\documentclass[onefignum,onetabnum]{siamart190516}

\usepackage{amsmath,amssymb,amsxtra}
\usepackage{mathrsfs}
\usepackage{graphicx}
\usepackage{float,verbatim}
\usepackage{threeparttable,booktabs}
\usepackage{underscore}
\usepackage{algorithm,algorithmic}
\usepackage{subfigure}
\usepackage{extarrows}
\usepackage{multirow}
\usepackage{makecell}
\usepackage{bm}
 \usepackage{color}
\graphicspath{{./InvPotenfigs/}}
\newtheorem{remark}{Remark}[section]
\newtheorem{condition}[theorem]{Condition}
\newtheorem{example}{Example}[section]
\newtheorem{assumption}{Assumption}[section]

\renewcommand{\d}{\mathrm{d}}

\def\II{(\Omega)}
\def\rbd{{\rm\bf d}}
\def\diam{{\rm diam}}
\def\dist{{\rm dist}}
\allowdisplaybreaks

\title{Stochastic Convergence Analysis of Inverse Potential Problem\thanks{The research of B. Jin is supported by Hong Kong RGC General Research Fund (Project 14306423 and 14306824) and ANR / RGC Joint Research Scheme (A-CUHK402/24), and a
start-up fund from The Chinese University of Hong Kong, and that of W. Zhang by the National Natural Science Foundation of China No. 12371423 and 12241104.}}

\author{Bangti Jin\thanks{Department of Mathematics, The Chinese University of Hong Kong, Shatin, New Territories, Hong Kong, P. R. China (\texttt{bangti.jin@gmail.com, b.jin@cuhk.edu.hk, qmquan@cuhk.edu.hk}).} \and Qimeng Quan\footnotemark[2] \and Wenlong Zhang\thanks{Department of Mathematics \& National Center for Applied Mathematics Shenzhen, Southern University of Science and Technology (SUSTech), 1088 Xueyuan
 Boulevard, University Town of Shenzhen, Xili, Nanshan, Shenzhen, Guangdong Province,
 P. R. China. (zhangwl@sustech.edu.cn)}}

\date{\today}

\begin{document}

\maketitle

\begin{abstract}
In this work, we investigate the inverse problem of recovering a potential coefficient in an elliptic partial differential equation from the observations at deterministic sampling points in the domain subject to random noise. We employ a least squares formulation with an $H^1\II$ penalty on the potential in order to obtain a numerical reconstruction, and the Galerkin finite element method for the spatial discretization. Under mild regularity assumptions on the problem data, we provide a stochastic $L^2\II$ convergence analysis on the regularized solution and the finite element approximation in a high probability sense. The obtained error bounds depend explicitly on the regularization parameter $\gamma$, the number $n$ of observation points and the mesh size $h$. These estimates provide a useful guideline for choosing relevant algorithmic parameters. Furthermore, we develop a monotonically convergent adaptive algorithm for determining a suitable regularization parameter in the absence of \textit{a priori} knowledge. Numerical experiments are also provided to complement the theoretical results.
\end{abstract}



\pagestyle{myheadings}
\thispagestyle{plain}

\section{Introduction}
In this work, we investigate the inverse problem of recovering a potential coefficient $q(x)$  in an elliptic partial differential equation (PDE). Let $\Omega\subset \mathbb{R}^d\ (d=1,2,3)$ be a simply connected open bounded domain with a boundary $\partial\Omega$. Consider the following boundary value problem with a zero Dirichlet boundary condition:
\begin{equation}\label{eqn:potential}
\left\{\begin{aligned}
		-\Delta u +qu &= f, \ &\mbox{in}&\ \Omega, \\
		u&=0, \ &\mbox{on}&\ \partial\Omega.
\end{aligned}\right.
\end{equation}
The function $f\in L^2(\Omega)$ in \eqref{eqn:potential} is a given source. The potential coefficient $q$ is assumed to belong to the admissible set
\begin{equation*}
	\mathcal{A} = \{q\in L^\infty\II: c_0\leq q(x)\leq c_1 \mbox{ a.e. in } \Omega\},
\end{equation*}
with lower and upper bounds $0\leq c_0<c_1<\infty$. (The coefficient $q$ may vanish identically in $\Omega$.) We use the notation $u(q)$  to indicate the dependence of the solution $u$ to problem \eqref{eqn:potential} on the potential $q$. The inverse potential problem is to recover the potential $q$ from the noisy observation of the exact data $u(q^\dagger)$ in the domain $\Omega$ (corresponding to the exact potential $q^\dagger$). It has several practical applications, including quantitative dynamic elastography (of recovering the density from one component of the displacement) \cite{choulli2021some}, the identification of the radiativity / reaction coefficient in heat conduction \cite{yamamoto2001simultaneous} and perfusion coefficient in Pennes' bio-heat equation \cite{pennes1948analysis,trucu2010space} etc.

In this work, the concerned inverse problem is to recover the potential $q^\dag\in\mathcal{A}$ from the pointwise measurement data $\bm m :=(m_i)_{i=1}^n$ subject to random noise:
\begin{equation}\label{eqn:rand-obs}
   m_i = u(q^\dag) (x_i) + \xi_i,\quad\mbox{with}\ (x_i)_{i=1}^n\subset\Omega,
\end{equation}
where the deterministic observational points $(x_i)_{i=1}^n$ are distributed quasi-uniformly over the domain $\Omega$, cf. Assumption \ref{ass:quasi-points}. The components of the random noise vector $\bm\xi=(\xi_i)_{i=1}^n$ are independent and identically distributed (i.i.d) sub-Gaussian random variables with a fixed parameter (noise level) $\sigma$ and zero expectation $\mathbb{E}[\bm \xi ]= \bm 0$; see  Definition \ref{def:sub-Gau} for the precise definition. The sub-Gaussian assumption holds for the popular Gaussian noise as well as uniform noise.

In practice, the observational data is measured at scattered points using sensors and contaminated by random measurement errors. Thus it is important to study the discrete observation model \eqref{eqn:rand-obs}. This model is challenging to analyze due to the following two facts. First, the nonlinearity of the parameter-to-state map $q\mapsto u(q)$ and the strong nonconvexity of the regularized functional render the existing argument \cite{chen2022stochastic} for linear inverse problems inapplicable. Chen et al \cite{chen2022stochastic} studied a regularized formulation of an inverse source problem with stochastic pointwise measurements and its discretization using the Galerkin FEM, and derived error estimates in a weaker norm for both regularization solution and finite element solution. The analysis employs the spectral decomposition of the associated linear forward operator and a noise separation strategy in the variational formulation. Unfortunately, these two crucial tools are no longer available for the inverse potential problem. Second, the scattered and stochastic nature of observational data poses big challenges to directly applying existing conditional stability and error analysis in \cite[Theorems 2.1 and 3.2]{jin2022convergence}. Indeed, one has to employ appropriate probabilistic tools to establish error estimates for observational data with random noise.

In this work, we employ Tikhonov regularization with an $H^1\II$ penalty for recovering the potential $q$ in the model \eqref{eqn:potential} and develop a numerical procedure based on the Galerkin finite element method (FEM). This work contributes in the following three aspects. First, we provide an \textit{a priori} bound in the $\psi_2$-Orlicz norm (with $\psi_2(t)=e^{t^2}-1$) on the regularized solution $q^*$ and the state approximation $u(q^\dag) - u(q^*)$ (i.e., prediction error) under a suitable choice of the regularization parameter $\gamma$  (with a fixed noise level $\sigma$):
\begin{equation}\label{ineq:pribound-q}
	\bigg\|\|u(q^*)-u(q^\dagger)\|_n\bigg\|_{\psi_2} + \gamma^{\frac12}\bigg\|\| q^*\|_{H^1\II}\bigg\|_{\psi_2} \leq c\gamma^{\frac12}\rho_0,
\end{equation}
with $\gamma^{{\frac12 + \frac{d}{12}}} = O(\sigma n^{-\frac12}(1+\rho_0)^{\frac d6}\rho_0^{-1})$, where the constant $\rho_0$ in \eqref{ineq:pribound-q} depends explicitly on  $\|q^\dagger\|_{H^1(\Omega)}$, the noise level $\sigma$ and the number $n$ of observational points.  The analysis crucially relies on a peeling argument \cite{geer2000empirical,chen2018stochastic,chen2022stochastic} and the sub-Gaussian nature of the noise vector $\bm \xi$, cf. Lemma \ref{lem:approx-sol-uppbound-q}.
Based on the estimate \eqref{ineq:pribound-q} and the conditional stability in \cite[Theorem 2.1]{jin2022convergence}, we establish \textit{a priori} stochastic convergence in the standard $L^2\II$ norm for the regularized solution $q^*$ in high probability. The error bound depends explicitly on the parameters $\gamma$, $\sigma$, $n$ and the probability tolerance $\tau$. See Theorem \ref{thm:con-err-potential} for the precise statement. The error bound is consistent with the conditional stability estimate, cf. Remark \ref{rmk:optimality}. Second, we discretize both Tikhonov functional and the governing PDE \eqref{eqn:potential} using the Galerkin FEM with continuous piecewise linear finite elements. We extend the theoretical findings for the continuous setting to the discrete one, and provide a stochastic $L^2(\Omega)$ error analysis of the discrete solution $q_h^*$ in high-probability; see Theorem \ref{thm:err-potential} for the precise statement. Due to the limited Sobolev regularity of the discrete state $u_h(q_h^*)$, the stochastic convergence analysis is far more technical and requires a markedly different proof strategy. Last, we propose an adaptive algorithm to choose a suitable $\gamma$ based on the error estimate, when no additional information is available, for each noise realization. The implementation of Algorithm \ref{alg:adpa-gamma} requires only an initial $\gamma_0$, amd then updates $\gamma$ iteratively by solving the discrete regularization problem. We prove a monotone convergence of the adaptive algorithm for each realization of $\bm m$. In addition, we present several numerical experiments to show the convergence behavior of the regularized solution and the adaptive algorithm. The stochastic convergence analysis of the continuous and discrete regularized problems in Theorems \ref{thm:con-err-potential} and \ref{thm:err-potential} represents the main contributions of the work.

Now we situate the present study in existing works on the concerned inverse potential problem in the context of variational regularization and statistical inversion, and discuss relevant results separately. In the context of variational regularization \cite{Engl1996Regularization,ItoJin:2015}, it has been extensively studied, and several reconstruction methods have been devised. These studies focus on the case that the data is available everywhere in the domain $\Omega$, and the noise $\xi(x)$ is a deterministic perturbation in the data space $L^2(\Omega)$ with a noise level $\delta:=\|\xi\|_{L^2\II}$, which is typically assumed to be small. The works \cite{yang2008inverse,yamamoto2001simultaneous} analyzed the well-posedness of the continuous regularized model for both elliptic and parabolic problems. In the literature, various conditions have been proposed to derive convergence rates in terms of the noise level $\delta$, including canonical source condition \cite{EnglKunischNeubauer:1989}. Hao and Quyen \cite{hao2010convergence} provided an approach that does not explicitly depend on the source condition. Chen et al \cite{Chen2020Convergence} investigated the inverse problem for elliptic and parabolic problems in negative Sobolev spaces, proved a novel conditional stability estimate and developed a variational inequality-type source condition, which also enables deriving convergence rates for the regularized solution.
Jin et al \cite{jin2022convergence} investigated a numerical procedure for the inverse problem in the elliptic/parabolic cases, which approximates both state $u$ and potential $q$ using the Galerkin FEM with continuous piecewise linear finite elements, and derived \textit{a priori} $L^2(\Omega)$ error estimates on the discrete solution $q_h^*$. The analysis relies on an energy argument and a novel H\"{o}lder-type conditional stability. See Section \ref{ssec:exist-reg} for further details about the results in \cite{Chen2020Convergence,jin2022convergence}. Moreover, Zhang et al \cite{zhang2022identification} presented a fixed-point iteration scheme for recovering the potential $q$ in the subdiffusion model from the terminal data, and proved an $L^2\II$ error analysis of the numerical solution.

In the context of statistical analysis of nonlinear PDE inverse problems, the discrete random measurement model and the Gaussian white noise model have received much recent attention, and the statistical theory has witnessed significant progress in the last few years; see the recent monograph \cite{Nickl:2023} for an up-to-date overview. The inverse potential problem for the stationary Schr\"{o}dinger model has been studied in the statistical setting in several  recent works  \cite{Nickl:2020,NicklWang:2020,NickWang:2024,Koers:2024}.
Nickl \cite{Nickl:2020} devised a nonparametric Bayesian prior on the potential coefficient $q$, and proved a minimax optimal contraction rate for the resulting posterior
distribution for both $q$ and $u(q)$ (up to log factors) \cite[Theorem 1 and Proposition 2]{Nickl:2020} and nonparametric Bernstein-von Mises theorems \cite[Theorem 8]{Nickl:2020}; see also \cite{NickWang:2024} for computational guarantees for both posterior mean vectors and \textit{maximum a posteriori} (MAP) estimates. Nickl et al \cite{NicklWang:2020} studied Tikhonov-type penalized least-square estimators of $q$ under squared Sobolev norm penalty, which can be viewed as a Bayesian ``MAP" estimator under some Gaussian process prior, and derived convergence rates for both $q$ and $u(q)$, including minimax-optimality for the prediction loss. In these works, convergence rates for the estimators of the potential $q$ were established by combining estimates on the prediction error and conditional stability estimates. Very recently, Koers et al \cite[Theorem 4.5]{Koers:2024} devised a linear method for the inverse potential problem, and proved optimal recovery rates for $u(q)$ and $q$. The present work differs from these interesting works in the following aspects: we investigate  recovering a potential $q^\dag\in H^2(\Omega)$ from pointwise observation at deterministic scattered points $(x_i)_{i=1}^n$ via the  estimator $q^*$, and treat explicitly the discretization error due to the Galerkin FEM approximation, whereas existing studies focus on estimating a potential $q$ with high Sobolev regularity (more precisely, $q\in H^\alpha(\Omega)$, with $\alpha>\max(d/2+2,2d-2)$) in the setting of random design and do not take into account the discretization error. We refer readers to Section \ref{ssec:exist-stat} for more detailed discussions.

The rest of the paper is organized as follows. In Section \ref{sec:stoc-conver}, we describe the continuous reconstruction model and its finite element approximation, and state \textit{a priori} $L^2\II$ error estimates in high probability. Moreover, we discuss relevant results in the literature. In Section \ref{sec:pre}, we recall basic theory on sub-Gaussian random variables. In Sections \ref{sec:error-cont} and \ref{sec:error-disc}, we provide the lengthy and technical proofs of Theorems \ref{thm:con-err-potential} and \ref{thm:err-potential}, respectively. In Section \ref{sec:numerics}, we give several numerical experiments to illustrate the theoretical analysis. In the appendix, we collect the proofs of two results on the existence of minimizers. We conclude the introduction with some useful notation.
For any $m\geq0$ and $p\geq1$, we denote by $W^{m,p}\II$ the standard Sobolev spaces of order $m$, equipped with the norm $\|\cdot\|_{W^{m,p}\II}$ and also write $H^{m}\II$ with the norm $\|\cdot\|_{H^m\II}$ when $p=2$ \cite{Adams2003Sobolev}. We denote the $L^2\II$ inner product by $(\cdot,\cdot)$. We denote by $C(\overline{\Omega})$ the space of continuous functions on $\overline{\Omega}$.
Throughout, we denote vectors by bold lowercase letters, and denote by $c$ a generic positive constant not necessarily the same at each occurrence but always independent of the sub-Gaussian parameter $\sigma$, the mesh size $h$, the number $n$ of measurement points, and the regularization parameter $\gamma$.

\section{Numerical procedure, main results and discussions} \label{sec:stoc-conver}
In this section, we describe the numerical procedure using a regularized output least-squares formulation to identify the potential $q$ and discretize the continuous formulation by the Galerkin FEM with continuous piecewise linear elements. Also we present the main results of the work and discuss several existing results about the inverse problem in detail.

\subsection{Regularized reconstruction and its discretization}\label{ssec:method}
To reconstruct the potential $q$, we
employ  Tikhonov regularization with an $H^1\II$ penalty \cite{Engl1996Regularization,ItoJin:2015}. This amounts to minimizing the following regularized functional:
\begin{equation}\label{eqn:conti-optim-prob-potential}
\begin{split}
    \min_{q\in \mathcal{A}} J_\gamma(q)=\|u(q)-\bm m\|_n^2 \!+\!\gamma\| q\|_{H^1\II}^2,
\end{split}
\end{equation}
where $\gamma>0$ is the regularization parameter balancing the two terms of the functional $J_\gamma(q)$, the function $u\equiv u(q)\in H^1_0\II$ satisfies
\begin{equation}\label{eqn:conti-weak-potential}	
   (\nabla u,\nabla\varphi)+(qu,\varphi)=(f,\varphi), \quad\forall \varphi\in H^1_0\II,
\end{equation}
and $\|\cdot\|_n$ is the discrete semi-norm defined in \eqref{eqn:semi-norm} below, i.e.,
$$\|u(q)-\bm m\|_n^2 \!=\frac1n\sum_{i=1}^n|u(q)(x_i)-m_i|^2.$$

By Proposition \ref{prop:exist-cont} below, the functional $J_\gamma(q)$ is well defined over the admissible set $\mathcal{A}$, and moreover, problem \eqref{eqn:conti-optim-prob-potential}--\eqref{eqn:conti-weak-potential} has at least one global minimizer $q^*\in \mathcal{A}\cap H^1(\Omega)$ $\mathbb{P}$ almost surely ($\mathbb{P}$ denotes the law of $\bm m$).

To approximate problem \eqref{eqn:conti-optim-prob-potential}--\eqref{eqn:conti-weak-potential} numerically, we employ the standard Galerkin FEM \cite{Brenner2002The}. Fix $h\in(0,h_0]$ for some $h_0>0$, and let $\mathcal{T}_h:=\cup\{K_j\}_{j=1}^{N_h}$ be a shape regular quasi-uniform simplicial triangulation of the domain $\Omega$ into mutually disjoint open face-to-face
subdomains $K_j$, such that $\Omega_h:= {\rm Int}(\cup_j\{\overline{K}_j\})\subset\Omega$ with all the boundary vertices of $\Omega_h$ locating on $\partial\Omega$ and ${\rm dist}(x,\partial\Omega)\leq ch^2$ for $x\in\partial\Omega_h$ \cite[Section 5.3]{larsson2003partial}. For any simplicial element $K\in\mathcal{T}_h$, we denote by $P_r(K)$ the space of polynomials of degree at most $r$ over $K$. On the triangulation $\mathcal{T}_h$, we define the conforming piecewise linear (P1) FEM spaces $V_h$ and $X_h$ (for the approximating polyhedron $\Omega_h$) by
\begin{equation*}
\begin{split}
&V_{h}:=\{v_{h}\in H^1(\Omega_h):v_{h}|_{K}\in P_1(K),\, \forall K\in\mathcal{T}_h\}, \\
& X_{h}:=\{v_{h}\in H^1_0(\Omega_h):v_{h}|_{K}\in P_1(K),\, \forall K\in\mathcal{T}_h\}.
\end{split}
\end{equation*}
 Similarly we can have conforming piecewise quadratic (P2) elements by taking $P_2(K)$ in the definitions.
Note that $V_h$ and $X_h$ can be naturally extended to the entire domain $\Omega$ by linear polynomials and zero function, respectively. We still denote the spaces of extended functions by $V_h$ and $X_h$. If the domain $\Omega$ is polyhedral, then we have $X_h = V_h\cap H_0^1\II$.
We use the spaces $V_h$ and $X_h$ to approximate the potential $q$ and the state $u$, respectively. Let $\mathcal{A}_h =\mathcal{A}\cap V_h$. Then the Galerkin FEM approximation of problem \eqref{eqn:conti-optim-prob-potential}--\eqref{eqn:conti-weak-potential} reads:
\begin{equation}\label{eqn:dis-optim-pure-dis}
	\min_{q_h\in \mathcal{A}_h} J_{\gamma,h}(q_h)=\|u_h(q_h)-\bm m\|^2_n + \gamma\| q_h\|_{H^1\II}^2,
\end{equation}
 where $u_h\equiv u_h(q_h)\in X_h$ satisfies
\begin{equation}\label{eqn:dis-weak-potential}	
	(\nabla u_h,\nabla\varphi_h)+(q_hu_h,\varphi_h)=(f,\varphi_h), \quad\forall \varphi_h\in X_h.
\end{equation}
 The objective functional $J_{\gamma,h}(q_h)$ is also well defined over the set $\mathcal{A}_h$, and moreover, problem \eqref{eqn:dis-optim-pure-dis}-\eqref{eqn:dis-weak-potential} has at least one minimizer $q_h^*\in\mathcal{A}_h$ $\mathbb{P}$ almost surely; see Proposition \ref{prop:exist-disc} below.

\subsection{Main results and discussions}
Next we state stochastic error bounds on the regularized approximation $q^*$ and the FEM approximation $q_h^*$ constructed via the numerical procedures in Section \ref{ssec:method} relative to the exact one $q^\dag$. Throughout we assume that the noise components $\bm \xi = (\xi_i)_{i=1}^d$ are i.i.d. sub-Gaussian random variables with parameter $\sigma$ and $\mathbb{E}\big[\xi_i\big] = 0$. For the convergence analysis, we need the following quasi-uniformity condition (introduced in \cite[Theorem 1.1]{utreras1988convergence}) on the scattered sampling points $(x_i)_{i=1}^n\subset \Omega$.

\begin{assumption}\label{ass:quasi-points}
The sampling points $(x_i)_{i=1}^n$ are scattered but quasi-uniformly distributed in the domain $\Omega$, i.e., there exists $b>0$ independent of $n$ such that $d_{\max} \leq bd_{\min}$ holds for all large $n$, where ${d_{\max}}$ and ${d_{\min}}$ are respectively defined by
\begin{equation*}
d_{\max}=\mathop {\rm sup}\limits_{x\in \Omega} \mathop {\rm inf}\limits_{1 \leq i \leq n} \|x-x_i\|_{\ell^2}
	\quad \mbox{and} \quad
		d_{\min}=\mathop {\rm inf}\limits_{1 \leq i \neq j \leq n} \|x_i-x_j\|_{\ell^2},
	\end{equation*}
	where $\|\cdot\|_{\ell^2}$ denotes the standard Euclidean norm of vectors.
\end{assumption}

\begin{figure}[hbt!]\label{fig:quasiuniform-mesh-points}
	\centering\setlength{\tabcolsep}{0pt}
	\begin{tabular}{ccc}
		 \includegraphics[width=0.32\textwidth,trim={0cm 0cm 0cm 0cm},clip]{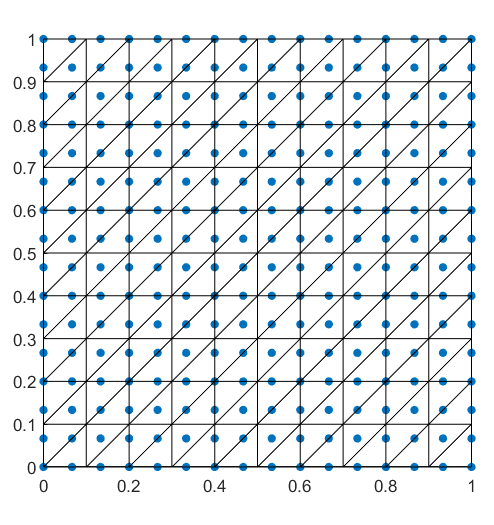} & \includegraphics[width=0.32\textwidth,trim={0cm 0cm 0cm 0cm},clip]{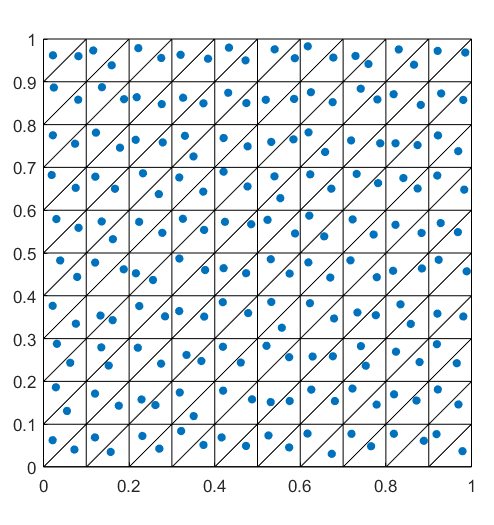} &
         \includegraphics[width=0.32\textwidth,trim={0cm 0cm 0cm 0cm},clip]{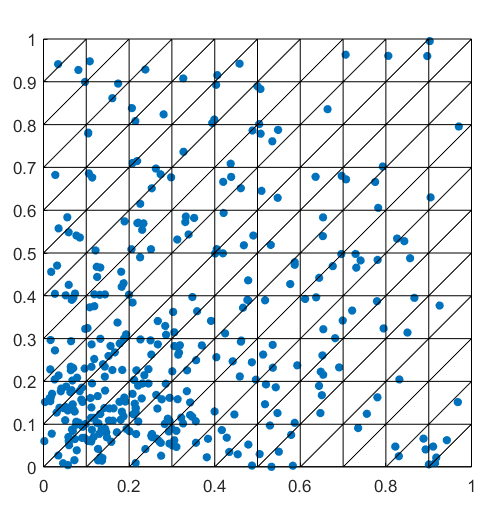}\\
         (a) uniform & (b) quasi-uniform & (c) non-quasi-uniform
	\end{tabular}
	\caption{A schematic illustration of  Assumption \ref{ass:quasi-points} on scattered sampling points $(x_i)_{i=1}^n$: (a) uniform, (b) quasi-uniform and (c) nonquasi-uniform.}
\end{figure}

Fig. \ref{fig:quasiuniform-mesh-points} shows two exemplary quasi-uniform distributions of sampling points $(x_i)_{i=1}^{n}$ over the unit square $\Omega=[0,1]^2$, as well as a non-quasi-uniform distribution of $(x_i)_{i=1}^{n}$. The quasi-uniform distribution in Fig. \ref{fig:quasiuniform-mesh-points} (b) is generated in two steps: we first determine the center of each triangule element in the uniform mesh, and then add a random  perturbation to the center. The density of the sampling points $(x_i)_{i=1}^{n}$ in Fig. \ref{fig:quasiuniform-mesh-points} (c) decays polynomially with the  radial distance from the origin (lower left corner). The sampling points  are not quasi-uniformly distributed since the minimum separation distance $d_{\min}$ is very small compared to $d_{\max}$, as indicated by the strong clustering around the origin.

We also make the  following regularity assumption on the problem data $q^\dag$ and $f$.
\begin{assumption}\label{ass:con-reg}
The domain $\Omega$ is smooth, $q^\dag\in H^1\II\cap\mathcal{A}$ and $f\in H^1\II$.
\end{assumption}

The next condition plays a crucial role in deriving the standard $L^2(\Omega)$ error bound. The distance function $\mathrm{dist}$ is defined by  $\mathrm{dist}(x,\partial\Omega)=\inf_{x'\in \partial\Omega}\|x-x'\|_{\ell^2}$.
\begin{condition}\label{Cond:positivity}
There exist $\beta,c>0$ such that $u^\dagger(x)\geq c\dist(x,\partial\Omega)^\beta~\mbox{a.e.}~ x\in \Omega$.
\end{condition}

Condition \ref{Cond:positivity} imposes a positivity condition on the exact state $u^\dag\equiv u(q^\dag)$ over the domain $\Omega$. The parameter $\beta>0$ quantifies the decay rate of $u^\dag$ near the domain boundary $\partial\Omega$: a larger $\beta$ indicates a faster decay, which implies a slower convergence rate of the regularized solution $q^*$; see Theorem \ref{thm:con-err-potential} and its proof for details. \cite[Proposition 2.1]{jin2022convergence} shows that the condition holds
with $\beta  = 2$ under Assumption \ref{ass:con-reg} with $f \geq  c_f > 0$.

Now we can state the first main result, i.e., weighted and standard $L^2\II$ error bounds on the regularized solution $q^*$. The proof proceeds in three steps: we first derive an \textit{a priori} bound on the prediction error $\|u^\dag-u(q^*)\|_n$ and the penalty $\|q^*\|_{H^1(\Omega)}$ under the given regularity assumption on the problem data, then establish a novel weighted $L^2\II$ error estimate (with the weight $u^\dag$) using the conditional stability argument (i.e., suitable test function in the weak formulation) under Assumptions \ref{ass:quasi-points} and \ref{ass:con-reg}; and last remove the weight function $u^\dag$ and obtain a standard $L^2\II$ error bound by utilizing the positivity condition on $u^\dag$ in Condition \ref{Cond:positivity}. The detailed proof of the theorem is given in Section \ref{sec:error-cont}. Note that the constant $c$ in the estimate is independent of $q^*$ and $n$ etc., and the estimate depends on $\sigma$ and $\|q^\dag\|_{H^1(\Omega)}$ via the factor $\rho_0$ only.

\begin{theorem}\label{thm:con-err-potential}
Let Assumptions \ref{ass:quasi-points} and \ref{ass:con-reg} be fulfilled, and let $q^*\in\mathcal{A}\cap H^1(\Omega)$ be a minimizer to problem \eqref{eqn:conti-optim-prob-potential}-\eqref{eqn:conti-weak-potential}.  Fix $\tau\in(0,\frac14)$,
$$\rho_0 = \|q^\dagger\|_{H^1\II} + \sigma n^{-\frac12}\quad\mbox{and} \quad \gamma^{{\frac12 + \frac{d}{12}}} = O(\sigma n^{-\frac12}(1+\rho_0)^{\frac d6}\rho_0^{-1}).$$ Then with probability at least $1-2\tau$, with $\ell_\tau = \sqrt{\log\frac{2}{\tau}}$, there holds
\begin{equation*}
\|(q^{\dag}-q^*)u^\dag\|^2_{L^2\II}\leq  c\big(\gamma^{\frac12} + n^{-\frac 3d} \big)^{\frac23}(1+\rho_0\ell_\tau)^3,
\end{equation*}
where the constant $c$ is independent of $q^*$, $\gamma$, $\rho_0$, $\tau$ and $n$, and dependent only on the domain $\Omega$ and the quasi-uniformity constant $b$.
Moreover, if Condition \ref{Cond:positivity} holds, then
\begin{equation*}
	\|q^{\dagger}-q^*\|_{L^2\II}\leq c\big[{\big(\gamma^{\frac12} + n^{-\frac 3d} \big)^{\frac23}{(1+\rho_0\ell_\tau)^3}}\big]^{\frac{1}{2(1+2\beta)}}.
\end{equation*}
\end{theorem}
\begin{remark}\label{rem:con-err-bound-match}
    Theorem \ref{thm:con-err-potential} provides  $L^2\II$ error bounds on the regularized approximation $q^*$ in the high probability sense. In the estimate, the factor $n^{-\frac{3}{d}}$ is due to the approximation error of the $L^2(\Omega)$ norm by the  discrete semi-norm $\|\cdot\|_n$, cf. Lemma \ref{lem:connect-seminorm-stdnorm} below. If the number $n$ of sampling points are large enough {\rm(}i.e., $\gamma\geq n^{-\frac6d}${\rm)}, the following error estimate holds with probability at least $1-2\tau$,    \begin{equation}\label{eqn:estimate-gamma}
        \|q^\dag-q^*\|_{L^2\II} \leq c\gamma^{\frac{1}{6(1+2\beta)}}[1+\rho_0\ell_\tau]^{\frac{3}{2(1+2\beta)}}.
    \end{equation}
The condition $\gamma\geq n^{-\frac 6d}$ is non-essential: To achieve a uniform bound with respect to $\gamma$, we use it to absorb the factor $n^{-\frac 3d}$. The estimate \eqref{eqn:estimate-gamma} is stated in terms of $\gamma$ only, and can be understood as the model error induced by the penalty term $\gamma \|q\|_{H^1(\Omega)}^2$ in $J_\gamma(q)$. Furthermore, if the condition that the number $n$ of sampling points is large, the \textit{a priori} choice $\gamma^{{\frac12 + \frac{d}{12}}} = O(\sigma n^{-\frac12}(1+\rho_0)^{\frac d6}\rho_0^{-1})$ implies the following estimate in terms of $n$:
\begin{equation}
    \|q^\dag-q^*\|_{L^2\II} \leq c(\sigma n^{-\frac12})^s \rho_0^r\ell_\tau^{\frac{3}{2(1+2\beta)}},
\end{equation}
with
$$s = \frac{2}{(6+d)(1+2\beta)}=\frac{12}{6+d}\!\cdot\frac{1}{3}\!\cdot\! \frac{1}{2(1+2\beta)} \quad\mbox{and}\quad r = \frac{11d+42}{6(6+d)(1+2\beta)},
$$
where the constant $c$ is independent of $\rho_0$, $\gamma$ and $n$. Note that for the three factors in the exponent $s$, the factor $\frac{12}{6+d}$ originates from the prediction error, which matches the minmax optimal rate for nonparametric regression {\rm(}cf. Remark \ref{rmk:rate-n} below{\rm)}, the factor $\frac{1}{3}$ stems from conditional stability {\rm(}cf. Remark \ref{rmk:optimality} below{\rm)}, and the factor $\frac{1}{2(1+2\beta)}$ is due to the positivity condition in Condition \ref{Cond:positivity}. Under Assumption \ref{ass:con-reg} and the condition $f\geq c_f>0$, we obtain the rate $\|q^\dag-q^*\|_{L^2(\Omega)}\leq cn^{-\frac{1}{5(6+d)}}$.
\end{remark}

\begin{remark}\label{rmk:optimality}
By Proposition \ref{prop:cond-stab}, the following conditional stability result holds 
\begin{equation*}
   \|(q-\widetilde q)u(q)\|_{L^2(\Omega)} \leq c(1+\max(\|q\|_{H^1(\Omega)},\|\widetilde q\|_{H^1(\Omega)}))^\frac56\|u(q)-u(\widetilde q)\|_{L^2(\Omega)}^{\frac13},
\end{equation*}
for any $q,\widetilde{q}\in \mathcal{A}\cap H^1(\Omega) $, where the constant $c$ depends only on $c_0$ and $c_1$.
The error bound in Remark \ref{rem:con-err-bound-match} is consistent with the stability estimate and the prediction error in Lemma \ref{lem:approx-sol-uppbound-q}. Note that this proof strategy does not lead to improved convergence rates for smoother potential coefficients, and hence the analysis focuses on the low-regularity regime. It is an interesting open problem to derive improved convergence rates for coefficients with better Sobolev regularity.
\end{remark}

To derive an error bound on the discrete approximation $q_h^*$, we need  the following slightly stronger regularity assumption on the problem data $q^\dag$ and $f$.
\begin{assumption}\label{ass:reg}
The domain $\Omega$ is smooth, $q^{\dagger}\in H^2\II\cap \mathcal{A}$ and $f\in H^1\II$.
\end{assumption}

Then we can state the second main result of the work. The analysis for $q_h^*$ is much more involved due to a lack of the necessary Sobolev regularity of FEM functions, and we overcome the challenge by using the finite-dimensionality of FEM approximations and the $H^2(\Omega)$-conforming interpolation operator. The technical proof of the theorem is given in Section \ref{sec:error-disc}.
\begin{theorem}\label{thm:err-potential}
Let Assumptions \ref{ass:quasi-points} and \ref{ass:reg} hold, and let $q_h^*\in\mathcal{A}_h$ be a minimizer to problem \eqref{eqn:dis-optim-pure-dis}-\eqref{eqn:dis-weak-potential}. Let
$$\eta:=\gamma^{\frac12} +h^2 + \gamma^{\frac{d}{12}}h^{4(\frac12-\frac{d}{8})}  + (\gamma^{\frac12+\frac{d}{12}}h^{4(\frac12-\frac{d}{8})})^{\frac12}.$$ Fix $\tau\in(0,\frac14)$ and $\rho_0 = \|q^\dagger\|_{H^1\II} + \sigma n^{-\frac12}$. Let $\gamma^{{\frac12 + \frac{d}{12}}} = O(\sigma n^{-\frac12}(1+\rho_0)^{\frac d6}\rho_0^{-1})$. Then with probability at least $1-2\tau$, with $\ell_\tau=\sqrt{\log\frac{2}{\tau}}$, there holds
\begin{equation*}
    \|(q^{\dag}-q_h^*)u^\dag\|^2_{L^2\II}\leq  c\big[\min(1, \big(\eta + h^2 + n^{-\frac3d} \big)^\frac23 + h ) \big]\big(1+\gamma^{-\frac12}\eta\rho_0\ell_\tau  \big)^3.
\end{equation*}
Moreover, if Condition \ref{Cond:positivity} holds, then
\begin{equation*}
    	\|q^{\dag}-q_h^*\|_{L^2\II}\leq c \Big(\big[\min(1, \big(\eta + h^2 + n^{-\frac3d} \big)^\frac23 + h ) \big]\big(1+\gamma^{-\frac12}\eta\rho_0\ell_\tau  \big)^3\Big)^{\frac{1}{2(1+2\beta)}}.
    \end{equation*}
\end{theorem}
\begin{remark}\label{rem:dis-err-bound-match}
Under the conditions in Theorem \ref{thm:err-potential}, the choice $h=O(\gamma^{\max(\frac{6-d}{6(4-d)},\frac13)})$ and $\gamma\geq n^{-\frac 6d}$, there holds with a probability at least $1-2\tau$,
\begin{equation*}
    	\|q^{\dag}-q_h^*\|_{L^2\II}\leq c\gamma^{\frac{1}{6(1+2\beta)}}[\rho_0\ell_\tau]^{\frac{3}{2(1+2\beta)}}.
\end{equation*}
This estimate is largely comparable with the continuous regularized approximation $q^*$, cf. Remark \ref{rem:con-err-bound-match}. Thus, the Galerkin discretization does not deteriorate the accuracy of the approximation.
This estimate exploits heavily the following two facts. First, the \textit{a priori} choice of the penalty parameter $\gamma^{{\frac12 + \frac {d} {12}}} = O(\sigma n^{-\frac12}(1+\rho_0)^{\frac d6} \rho_0^{-1})$ determines the model error of the regularized formulation  \eqref{eqn:conti-optim-prob-potential}-\eqref{eqn:conti-weak-potential}. Second, the condition $h=O(\gamma^{\max(\frac{6-d}{6(4-d)},\frac13)})$ ensures that the overall error is dominated by the model error or $\gamma$, and refining the mesh reduces only the discretization error but not the overall error. Note that the obtained rate is slightly higher than that in the deterministic case \cite[Corollary 3.2]{jin2022convergence}. This is attributed to the utilization of the a priori Lipschitz stability {\rm(}cf. Lemma \ref{lem:reg-u}{\rm)} when imposing the $H^1(\Omega)$ regularization.
\end{remark}

The proofs of Theorems \ref{thm:con-err-potential} and \ref{thm:err-potential} are given in Sections \ref{sec:error-cont} and \ref{sec:error-disc}, respectively.

\subsection{Review on existing studies}

In the literature, the inverse potential problem has been investigated in the frameworks of deterministic variational regularization and statistic analysis, under different assumptions on the problem setting. In this part, we discuss several closely related existing results in these two directions separately.

\subsubsection{Deterministic variational regularization}\label{ssec:exist-reg}

In this setting, one aims at recovering a potential $q^\dag\in \mathcal{A}$ from the \textit{deterministic} observational data $z^\delta$, which is a noisy version of the exact data $u(q^\dag)$, with a noise level $\delta = \|z^\delta - u(q^\dag)\|_{L^2(\Omega)}$. Note that the perturbation $\xi:=z^\delta-u(q^\dag)$ is assumed to belong to the data space $L^2(\Omega)$. Then the standard Tikhonov regularization with the Sobolev $H^s(\Omega)$ penalty (for $s=0,1$) for reconstructing the potential coefficient $q^\dag$ is to minimize
\begin{equation}\label{eqn:var-reg}
    J_\gamma (q) = \|u(q)-z^\delta\|_{L^2(\Omega)}^2 + \gamma \|q-q_0\|_{H^s(\Omega)}^2,
\end{equation}
where $q_0\in \mathcal{A}\cap H^s(\Omega)$ is an \textit{a priori} estimate of the exact potential  $q^\dag$, and $\gamma>0$ is the penalty parameter. Let $q_\gamma^\delta$ be a minimizer of $J_\gamma(q)$ over the admissible set $\mathcal{A}$. This procedure is well established for the inverse potential problem (see, e.g., \cite{EnglKunischNeubauer:1989}), and there are two recent studies \cite{Chen2020Convergence,jin2022convergence} on the error analysis of the approximations. Chen et al \cite{Chen2020Convergence} investigated the standard $L^2(\Omega)$ penalty (i.e., $s=0$ in \eqref{eqn:var-reg}). Then under the assumptions $|u(q^\dag)|\geq c_0>0$ in $\Omega$ and $q^\dag-q_0\in H_0^\kappa(\Omega)$ with $\kappa>0$ and $\kappa\neq1/2$, and $\alpha<\frac{\kappa}{1+\kappa}$, Chen et al \cite[Theorem 3.7]{Chen2020Convergence} proved that for the $L^2(\Omega)$ penalty, with the \textit{a priori} choice $\gamma=O(\delta^{2-\alpha})$, the following convergence rate holds
\begin{equation*}
\|q_\gamma^\delta-q^\dag\|_{L^2(\Omega)}\leq c\delta^{\frac{\alpha}{2}}.
\end{equation*}
The proof is based on a new source condition of variational inequality type for the problem, which is obtained using interpolation theory in Sobolev spaces and suitable \textit{a priori} estimates in negative Sobolev spaces. Note that the condition $|u(q^\dag)| \geq c_0 > 0$ in $\Omega$ can only hold for a nonzero Dirichlet boundary condition. In addition, Chen et al \cite{Chen2020Convergence} analyzed noisy gradient type measurement (i.e., $z^\delta$ is the noisy version of $\nabla u$) and the problem of recovering the space-dependent potential coefficient $q$ in the standard parabolic problem
\begin{equation} \label{eqn:parabolic}
\partial_t u - \Delta u + qu = f,\quad \mbox{in }\Omega\times(0,T],
\end{equation}
where $T>0$ is a given final time.

Jin et al \cite{jin2022convergence} investigated the Galerkin approximation of problem \eqref{eqn:var-reg} with an $H^1(\Omega)$ semi-norm penalty (and $q_0=0$). Then the Galerkin FEM approximation of problem \eqref{eqn:var-reg} reads
\begin{equation*}
    J_{\gamma,h}(q_h) = \|u_h(q_h)-z^\delta\|_{L^2(\Omega)}^2 + \gamma \|\nabla q_h\|_{L^2(\Omega)}^2,
\end{equation*}
with the FEM function $u_h(q_h)\in X_h$ satisfying \eqref{eqn:dis-weak-potential}.
Then for $q^\dag\in H^2(\Omega)\cap \mathcal{A}$ and $f\in L^2(\Omega)$, any global minimizer $q_h^*\in\mathcal{A}_h$ satisfies \cite[Theorem 3.2]{jin2022convergence}
\begin{equation*}
    \|(q^\dag-q_h^*)u(q^\dag)\|_{L^2(\Omega)} \leq c(h^\frac12 +\gamma^\frac14 + \min(h^\frac12+h^{-\frac12}\eta^\frac12,1))\gamma^{-\frac14}\eta^\frac12,\quad \mbox{with } \eta = h^2+\delta + \gamma^\frac12.
\end{equation*}
This estimate can also be converted into the standard $L^2(\Omega)$ norm under condition \eqref{Cond:positivity}.
The analysis is based on a conditional stability argument, and relies crucially on a suitable test function. The work \cite{jin2022convergence} also analyzed the inverse potential problem in the standard parabolic case \eqref{eqn:parabolic}. Note that the error analysis of the Galerkin scheme for problem \eqref{eqn:var-reg} with an $L^2(\Omega)$ penalty remains unavailable.

Both works \cite{Chen2020Convergence,jin2022convergence} do not directly employ standard range-type source conditions to derive convergence rates (as traditionally done for nonlinear PDE inverse problems \cite{Engl1996Regularization,ItoJin:2015}), and instead employ condition stability (using an energy argument) and positivity condition on the weight $u^\dag$. Due to the differences in the analysis strategy, the requisite regularity condition on the coefficient $q$ is different in the two works \cite{Chen2020Convergence,jin2022convergence}. This work follows closely the analysis strategy in \cite{jin2022convergence} but handles discrete measurements at deterministic sampling points subject to random noise. This requires treating random quantities suitably using probabilistic tools, especially maximal inequality for sub-Gaussian random processes in Lemma \ref{lem:max-ineq} and a technical lemma (cf. Lemma \ref{Geer-entropy}) due to van de Geer \cite[Lemma 8.4 and (10.6)]{geer2000empirical}.

\subsubsection{Statistical analysis}\label{ssec:exist-stat}
In the framework of statistical inversion, the observations are corrupted by statistical noise. One popular measurement model is the following random design regression model:
\begin{equation}\label{eqn:Gaussian-disc}
    Y_{q^\dag,i} = u(q^\dag)(X_i) + \varepsilon_i,\quad i = 1,2,\ldots,n,
\end{equation}
where the i.i.d. random design points $(X_i)_{i=1}^n$ follow the uniform distribution over the domain $\Omega$, and the noise components $(\varepsilon_i)_{i=1}^n$ are random variables following independent standard normal distribution. Following the argument in asymptotic statistics \cite[Section 1.2.3]{GineNickl:2016}, the discrete model \eqref{eqn:Gaussian-disc} is asymptotically equivalent to the following Gaussian white noise model:
\begin{equation}\label{eqn:Gaussian}
    Y_{q^\dag}^{(\epsilon)} = u(q^\dag) + \epsilon \mathbb{W}, \quad \epsilon=1/\sqrt{n},
\end{equation}
where $(\mathbb{W}(\psi),\psi\in L^2(\Omega))$ is a centered Gaussian white noise process indexed by the space $L^2(\Omega)$  \cite[pp. 19-20]{GineNickl:2016}.
The noise $\mathbb{W}$ does not belong to the standard data space $L^2(\Omega)$ (as commonly employed in variational regularization, cf. Section \ref{ssec:exist-reg}) but only $H^\eta(\Omega)$ with $\eta<-d/2$. Thus the deterministic analysis in Section \ref{ssec:exist-reg} does not apply directly to the model \eqref{eqn:Gaussian}.

The statistical analysis of the inverse potential problem for the stationary Schr\"{o}dinger equation under the measurement model \eqref{eqn:Gaussian} was first thoroughly studied in the work of Nickl \cite{Nickl:2020}. Nickl \cite{Nickl:2020} devised a nonparametric Bayesian prior on the coefficient $q$ using a wavelet basis, and proved a minimax optimal contraction rate for the resulting posterior distribution for both $u(q)$ and $q$ (up to log factors) \cite[Theorem 1 and Proposition 2]{Nickl:2020}. Moreover, \cite[Theorem 8]{Nickl:2020} addressed the important question of uncertainty quantification using the posterior distribution, by proving nonparametric Bernstein-von Mises theorems. Later Nickl and Wang \cite{NickWang:2024} established polynomial-time computational guarantees of the posteriori measure by Langevin type algorithms, for both posterior mean vectors as well as optimization based maximum a posteriori (MAP) estimates, provided that the parameter $q$ can be sufficiently well approximated by its orthogonal projection into a finite-dimensional subspace. (Note that due to the nonlinearity of the forward operator, the MAP estimate is different from the conditional mean of the posterior distribution.) Kekkonen \cite{Kekkonen:2022} investigated the inverse problem of recovering of a space-dependent potential coefficient $q(x)$ in the standard parabolic equation \eqref{eqn:parabolic} from discrete observations with random noise in the  space-time cylinder $\Omega\times(0,T)$ (similar to the model \eqref{eqn:Gaussian-disc}) under a rescaled Gaussian prior and the truncated Gaussian prior, and proved several posterior contraction rates.

Statistically, one may also consider the MAP estimators (i.e., penalized least squares estimators). In this aspect, the seminal work of Nickl et al \cite{NicklWang:2020} is closely related to the present study. Specifically, let $\mathbb{E}_{q^\dag}^\epsilon$ denote the expectation operator under the law $\mathbb{P}_{q^\dag}^\epsilon$ of $Y_{q^\dag}^{(\epsilon)}$ in \eqref{eqn:Gaussian}. Given a constant $c_{\min}>0$, let the sets $\mathcal{Q}$ and $\mathcal{Q}_{\alpha,r}(R)$ (for $r>c_{\min}$ and $R>0$) be defined respectively by
\begin{align*}
    &\mathcal{Q}:=\left\{q\in H^\alpha(\Omega): q>c_{\min} \mbox{ in }\Omega, q=1 \mbox{ on } \partial\Omega, \partial_n^j q =0, j=1,\ldots,\alpha-1\mbox{ on }\partial\Omega\right\}, \\
   & Q_{\alpha,r}(R):=\{q\in \mathcal{Q}: q>r\mbox{ in } \Omega, \|q\|_{H^\alpha(\Omega)}\leq R\},
\end{align*}
where $\partial_n$ denotes taking the normal derivative.
Then under the measurement model \eqref{eqn:Gaussian}, Nickl et al \cite{NicklWang:2020} investigated the following generalized Tikhonov regularized functional
\begin{equation*}
    J_{\lambda,\epsilon} (q)= -2\langle Y_{q^\dag}^{(\epsilon)},u(q)\rangle +\|u(q)\|_{L^2(\Omega)}^2 + \lambda \|\Phi^{-1}\circ q\|_{H^\alpha(\Omega)}^2,
\end{equation*}
where $\circ$ and $\langle\cdot,\cdot\rangle$ denote composition and duality product, respectively, and $\Phi$ is a regular link function \cite[Definition 7]{NicklWang:2020}. The loss $J_{\lambda,\epsilon}(q)$ can be interpreted as a generalized penalized least-squares estimator for $q$ under suitable Gaussian process priors.  Under the condition $\alpha>\max({d}/{2}+2,2d-2)$, Nickl et al \cite[Theorems 3.12 and 3.13]{NicklWang:2020} proved that for any $q^\dag\in\mathcal{Q}$, with the \textit{a priori} choice $\lambda_\epsilon = \epsilon^{\frac{4(\alpha+2)}{2(\alpha+2)+d}}$, there exists a minimizer to $\hat q_\epsilon \in \mathcal{Q}$ almost surely under $\mathbb{ P}_{q^\dag}^\epsilon$, and moreover, for each $R>0$ and $r>c_{\min}$, there hold
\begin{align}
    \sup_{q^\dag \in \mathcal{Q}_{\alpha,r}(R)} \mathbb{E}_{q^\dag}^\epsilon\|u(\widehat q_\epsilon)-u(q^\dag)\|_{H^\beta(\Omega)} &\leq c\epsilon^\frac{2(\alpha+2-\beta)}{2(\alpha+2)+d},\label{eqn:est-Nickl1}\\
    \sup_{q^\dag \in \mathcal{Q}_{\alpha,r}(R)} \mathbb{E}_{q^\dag}^\epsilon\|\widehat{q}_\epsilon-q^\dag\|_{L^2(\Omega)}&\leq c\epsilon^\frac{2\alpha}{2(\alpha+2)+d}.\label{eqn:est-Nickl2}
\end{align}
(For the estimate \eqref{eqn:est-Nickl2}, one has to assume also a positivity condition on the Dirichlet boundary condition.) These rates are minimax optimal in view of \cite[Proposition 2]{Nickl:2020}. The analysis is based on an abstract framework developed in \cite[Section 2]{Nickl:2020}, in which a stability estimate in negative Sobolev norms plays a crucial role, cf. the estimate \eqref{eqn:ell-neg} in Remark \ref{rmk:neg-stab-est} below for related discussions.

Note that these results \eqref{eqn:est-Nickl1}--\eqref{eqn:est-Nickl2} are derived under the condition $q^\dag \in \mathcal{Q}_{\alpha,r}$ with $\alpha >\max({d}/{2}+2,2d-2)$, whereas Theorem \ref{thm:con-err-potential} of this work focuses on the low regularity regime (i.e., $q\in H^1(\Omega)$ and only the exact potential $q^\dag\in H^2(\Omega)$, cf. Assumption \ref{ass:reg}). Thus these results are derived in different regularity regime, and this work complements existing results. Since the $H^\alpha(\Omega)$ regularity on $q$ implies $u(q)\in H^{\alpha+2}(\Omega)$, the prediction error $\|u^\dag-u(q^*)\|_n$ in Lemma \ref{lem:approx-sol-uppbound-q} below is comparable with the estimate \eqref{eqn:est-Nickl1} (with $\beta=0$ and $\alpha=1$), cf. Remark \ref{rmk:rate-n} below. Furthermore, Theorem \ref{thm:err-potential} analyzes the discretization error of the Galerkin FEM approximation $q_h^*$. Nonetheless, at a high level, this work shares the similarity with the work \cite{NicklWang:2020} in that it follows the paradigm of combining the consistency error (on the objective functional) with suitable conditional stability estimates.

Very recently Koers et al \cite{Koers:2024} proposed a linear approach to solve the inverse potential problem by first solving a linear inverse problem (of estimating $-\Delta u(q)$) and then recovering the coefficient $q$ via a pointwise nonlinear relation $q= (\Delta u(q)+f)/u(q)$, cf. \eqref{eqn:potential}. Koers et al derived optimal (adaptive) contraction rates for the induced posteriori distribution of $q$ and also provided frequentist coverage guarantees for the corresponding Bayesian credible sets \cite[Theorem 4.5]{Koers:2024}; see also \cite[Theorem 5.1]{Koers:2024} for analogous results in the parabolic case, cf. \eqref{eqn:parabolic}. In a different direction, Siebel \cite{Siebel:2024} further developed the approach in \cite{NicklWang:2020} for the measurement model \eqref{eqn:Gaussian-disc} for the conductivity problem, i.e., recovering the conductivity coefficient $a$ in the elliptic equation $-\nabla\cdot(a\nabla u)=f$ from the observation of the state $u$ over the whole domain $\Omega$.

\section{Preliminaries}
\label{sec:pre}
In this section we recall preliminary materials, including basic theory of sub-Gaussian random variables and discrete Sobolev semi-norms. Throughout, these tools will play a crucial role in the convergence analysis in Sections \ref{sec:error-cont} and \ref{sec:error-disc}.
\subsection{Basic theory of sub-Gaussian random variables}
First, we recall the definition of sub-Gaussian random variables and useful properties of empirical processes. The notation $\mathbb{E}$ denotes taking expectation with respect to the involved random variables.
\begin{definition}\label{def:sub-Gau}
	A random variable $Z$ is sub-Gaussian with parameter $\sigma$ if it satisfies
	\begin{equation*}
		\mathbb{E}\big[\exp(\lambda (Z - \mathbb E[Z]))\big] \leq \exp\big(\tfrac12\sigma^2\lambda^2\big),\quad \forall\lambda\in\mathbb{R}.
	\end{equation*}
The probability distribution function of a sub-Gaussian random variable $Z$ has an exponentially decaying tail, i.e.,
    \begin{equation*}
    	\mathbb{P}(|Z - \mathbb{E} [Z]|\ge z)\leq 2\exp\big(-\tfrac{z^2}{2\sigma^2}\big),\quad\forall z>0.
    \end{equation*}
\end{definition}

Next we introduce the Orlicz norm. Let $\psi$ be a monotonically increasing convex function on $\mathbb{R}_+\cup\{0\}$ with $\psi(0)=0$,
the Orlicz norm $\|Z\|_\psi$ of a random variable $Z$ is given by
\begin{equation*}
	\|Z\|_{\psi} = \inf\biggl\{c>0:\mathbb{E}\Big[\psi\Big(\frac{|X|}c\Big)\Big]\le 1\biggr\}.
\end{equation*}
Throughout, we choose the function $\psi_2(t):=\exp(t^2)-1$ for $t\geq0$. Then we have the following estimate \cite[(4.5)]{chen2018stochastic}.
\begin{lemma}\label{lem:bdd-sub-Gaussian}
    For any sub-Gaussian random variable $Z$ with $\|Z\|_{\psi_2} < \infty$, the following estimate holds:
\begin{equation}\label{ineq:proba-distrib-estimate}
	\mathbb{P}(|Z|\ge z)\le 2\exp\bigg(-\frac{z^2}{\|Z\|_{\psi_2}^2}\bigg),\quad \forall z>0.
\end{equation}
\end{lemma}

The next lemma provides a converse statement of the estimate \eqref{ineq:proba-distrib-estimate}.
\begin{lemma}\label{lem:Orlicz-bound}
Let $c_1$ and $c_2$ be two positive constants. If a random variable $Z$ satisfies
\begin{equation*}
	\mathbb{P}(|Z|>\alpha (1+z))\leq c_1\exp\big(-c_2^2z^2\big), \quad \forall \alpha>0, ~~z\ge 1,
\end{equation*}
then there exists a constant $c_3>0$ depending on $c_1$ and $c_2$ such that
$\|Z\|_{\psi_2}\le c_3\alpha$.
\end{lemma}

We will use the maximal inequality on sub-Gaussian random processes indexed by a semimetric space frequently. First we recall the concept of a semimetric space.
\begin{definition}
Fix a set $\mathbb{T}$. A function $\rbd$ defined on $\mathbb{T}$ is called a semimetric if it satisfies the following three properties:
    \begin{itemize}
        \item [{\rm(i)}] identity and non-negativity: $\rbd(t,t) =0$ and $\rbd(s,t) \geq0$ for any $s, t\in\mathbb{T}$;
        \item [{\rm(ii)}] symmetry: $\rbd(s,t) = \rbd(t,s)$ for any $s, t\in\mathbb{T}$;
        \item [{\rm(iii)}] triangle inequality: $\rbd(s,t) \leq \rbd(s,s') + \rbd(s',t)$ for all $s, s', t\in\mathbb{T}$.
    \end{itemize}
Moreover, the pair $(\mathbb{T},\rbd)$ is called a semimetric space.
\end{definition}

\begin{definition}\label{def:sub-Gau-process}
	Let $\mathbb{T}$ be a semimetric space with a semimetric $\rbd$ and $\{Z_t:t\in \mathbb{T}\}$ a random process indexed by $\mathbb{T}$. The random process $\{Z_t:t\in \mathbb{T}\}$ is called sub-Gaussian if
	\begin{equation*}
		\mathbb{P}(|Z_s-Z_t|>z)\le 2\,\exp\bigg( -\frac{z^2}{2\,\rbd(s,t)^2} \bigg),\quad \forall s,t\in \mathbb{T}, ~~z>0.
	\end{equation*}
\end{definition}
\begin{definition}\label{def:cov-num}
Let $\mathbb{T}$ be a semimetric space with a semimetric $\rbd$. A collection of points $(t_i)_{i=1}^n \subset \mathbb{T}$ is called an
    $\epsilon$-cover of $\big(\mathbb{T},\rbd\big)$ if for any $t\in \mathbb{T}$, there
    exists at least one $i \in \{1,\dots,n\}$ such that $\rbd(t, t_i) \leq \epsilon$. The
    $\epsilon$-covering number $\mathcal{N}(\epsilon,\mathbb{T},\rbd)$ is the minimal
    cardinality among all $\epsilon$-covers of $\mathbb{T}$.
\end{definition}

We shall use the following maximal inequality frequently \cite[Section 2.2.1]{vanDerVaart:1996}.
\begin{lemma}\label{lem:max-ineq}
Let $\mathbb{T}$ be a semimetric space with a semimetric $\rbd$, and let $\diam\, {\mathbb{T}}$ be the diameter of the set $\mathbb{T}$ with respect to the semimetric $\rbd$. If $\{Z_t:t\in \mathbb{T}\}$ is a separable sub-Gaussian random process, then
	\begin{equation*}
		\bigg\|\sup_{s,t\in \mathbb{T}}|Z_s-Z_t|\bigg\|_{\psi_2}\leq c\int^{\diam\, \mathbb{T}}_0\sqrt{\log N\big(\tfrac12\epsilon, \mathbb{T},\rbd\big)}~{\rm d}\epsilon.
	\end{equation*}
\end{lemma}

The following lemmas give useful bounds on the covering entropy of Sobolev subsets due to Birman and Solomyak \cite[p. 311, Theorem 5.2]{birman1967piecewise} and finite dimensional subsets \cite[Corollary 2.6]{geer2000empirical}.
Lemma \ref{lem:cover-en-bound-Sobolev} is stated for the case that the domain $\Omega$ is a unit cube in $\mathbb{R}^d$. For a general domain $\Omega$, one can first find a cube enclosing $\Omega$ and then use an extension argument and a scaling argument to derive an analogous result.
\begin{lemma}\label{lem:cover-en-bound-Sobolev}
Let $Q$ be the unit cube in $\mathbb{R}^d$ and $SW^{s,p}(Q)$ be the unit sphere of space $W^{s,p}(Q)$ for $s> 0$ and $p\ge 1$. Then for small $\epsilon>0$, there holds
\begin{equation*}
	\log N(\epsilon, SW^{s,p}(Q), \|\cdot\|_{L^q(Q)})\le c\epsilon^{-d/s},
\end{equation*}
where the integrability index $q$ satisfies the following condition with $q^*=p(1-sp/d)^{-1}$,
\begin{equation*}
\left\{\begin{aligned}
   1\le q\le\infty, &\quad \mbox{if } sp>d,\\
   1\leq q<\infty, &\quad \mbox{if }sp =d,\\
   1\le q < q^*, &\quad \mbox{if } sp< d.
\end{aligned}\right.
\end{equation*}
\end{lemma}

\begin{lemma}\label{lem:cover-en-bound-finite}
Let $G\subset L^2\II$ be a finite-dimensional subspace of dimensionality $N_G=\dim(G)>0$, and $G_R=\{g\in G: \|g\|_{L^2\II}\le R\}$. Then for small $\epsilon>0$, there holds
	\begin{equation*}
		\log N(\epsilon,G_R,\|\cdot\|_{L^2\II})\leq N_G\log(1+4R\epsilon^{-1}).
	\end{equation*}
\end{lemma}

\subsection{Discrete Sobolev semi-norm}

Pointwise measurements necessitate the use of discrete Sobolev semi-norms $\|\cdot\|_n$ in the regularized models. We first recall one classic Sobolev embedding result \cite[Theorem 4.12, p. 85]{Adams2003Sobolev}, which will be used frequently below.
\begin{lemma}\label{lem:Sob-embed}
For $d=1,2,3$,  the space   $H^2\II$ compactly embeds into $C(\overline{\Omega})$, and the space $H^1\II$ continuously embeds into $ L^4(\Omega)$.
\end{lemma}

For any $u,v\in C(\overline{\Omega})$ and $y\in \mathbb{R}^n$, we define
\begin{equation*}	(y,v)_n=n^{-1}\sum^n_{i=1}y_iv(x_i)\quad\mbox{and}\quad
(u,v)_n=n^{-1}\sum^n_{i=1}u(x_i)v(x_i),
\end{equation*}
and the discrete semi-norm \begin{equation}\label{eqn:semi-norm} \|u\|_n=\big(n^{-1}\sum_{i=1}^{n} |u(x_i)|^2 \big)^{1/2},\quad \forall u\in C(\overline{\Omega}).
\end{equation}

The next lemma connects the standard Sobolev norm with the discrete semi-norm \cite[Theorems 3.3 and 3.4]{utreras1988convergence}. The notation $|\cdot|_{H^3(\Omega)}$ denotes the $H^3(\Omega)$ semi-norm. Note that $\|\cdot\|_n$ is well-defined on the space $H^3(\Omega)$ due to the continuous embedding $H^3\II\hookrightarrow C(\overline{\Omega})$ (for $d=1,2,3$), cf. Lemma \ref{lem:Sob-embed}.
\begin{lemma}\label{lem:connect-seminorm-stdnorm}
Under Assumption \ref{ass:quasi-points}, there exists $c=c(\Omega,B)$ such that
	\begin{equation*}
		\|v\|_{L^2\II}\le c\big(\|v\|_n+ n^{-\frac3d}|v|_{H^3\II}\big)\quad \mbox{and}\quad \|v\|_n\le c\big(\|v\|_{L^2\II}+ n^{-\frac3d}|v|_{H^3\II}\big),\quad \forall v\in H^3(\Omega).
	\end{equation*}
\end{lemma}

\section{Proof of Theorem \ref{thm:con-err-potential}}\label{sec:error-cont}

In this section, we provide a stochastic convergence analysis of the minimizer $q^*$ to the regularized problem \eqref{eqn:conti-optim-prob-potential}-\eqref{eqn:conti-weak-potential} under Assumption \ref{ass:con-reg}.

First we give the existence of a global minimizer $q^*\in \mathcal{A}$ to the regularized problem \eqref{eqn:conti-optim-prob-potential}-\eqref{eqn:conti-weak-potential}. The proof of the lemma follows by a standard argument (i.e., the direct method in calculus in variation \cite{Dacorogna:2008}), and the complete proof is given in Appendix \ref{app:tech-proof}. The notation $\mathbb{P}$ denotes the law of the noisy measurement $\bm{m}$.
\begin{proposition}\label{prop:exist-cont}
The regularized functional $J_\gamma(q)$ is well defined on the admissible set $\mathcal{A}$, and moreover, there exists at least one global minimizer $q^*\in\mathcal{A}$ to problem \eqref{eqn:conti-optim-prob-potential}-\eqref{eqn:conti-weak-potential} $\mathbb{P}$ almost surely.
\end{proposition}

Next we give the Sobolev regularity and Lipschitz continuity of the state $u(q)$.
\begin{lemma}\label{lem:reg-u}
For any $q\in \mathcal{A}$, if the domain $\Omega$ is a convex polyhedron, there exists a constant $c$ independent of $q$ such that $$\|u(q)\|_{H^2(\Omega)}\leq  c.$$
Moreover, under Assumption \ref{ass:con-reg}, for any $q,\widetilde q\in\mathcal{A}\cap H^1(\Omega)$, there holds
\begin{equation}\label{eqn:Lipschtz-cont}
    \|u(q)-u(\widetilde{q})\|_{H^3(\Omega)}\leq c(1 + \|q\|_{H^1(\Omega)})\|q-\widetilde{q}\|_{H^1(\Omega)}.
\end{equation}
\end{lemma}
\begin{proof}
By the Lax-Milgram theorem, problem \eqref{eqn:conti-weak-potential} has a unique solution $u\in H_0^1(\Omega)$ such that  $\|u\|_{H^1(\Omega)}\leq c$ with $c$ independent of $q$. Now for a convex polyhedral domain $\Omega$, we rewrite problem \eqref{eqn:conti-weak-potential} into
\begin{equation*}
    \left\{\begin{aligned}
        - \Delta u(q)& = f-qu(q), \quad \mbox{in }\Omega,\\
        u(q) & = 0,\quad \mbox{on }\partial\Omega.
    \end{aligned}\right.
\end{equation*}
Then from the full regularity pickup of the elliptic operator $-\Delta$ (with a zero Dirichlet boundary) in $H^2\II$ on a convex polyhedral domain \cite[Theorem 3.1.2.1, p. 139]{Grisvard:2011}, the box constraint $q\in \mathcal{A}$ and the condition $f\in L^2(\Omega)$, it follows that
\begin{equation*}
    \|u(q)\|_{H^2\II} \leq c\|f-qu\|_{L^2\II} \leq c\big(\|f\|_{L^2\II} + \|q\|_{L^\infty\II} \|u\|_{L^2\II}\big) \leq c\|f\|_{L^2(\Omega)} \leq c.
\end{equation*}
This proves the first estimate. Next, the standard elliptic regularity theory \cite[Theorem 9.25]{Brezis2010FunctionalAS}, Lemma \ref{lem:Sob-embed} and Assumption \ref{ass:con-reg} imply $ u(q)\in H^2\II\hookrightarrow C(\overline{\Omega})$, for $d=1,2,3$. Let $w=u(q)-u(\widetilde{q})$. Then $w\in H_0^1(\Omega)$ satisfies
\begin{equation*}
    -\Delta w =  \widetilde q u(\widetilde q)-qu(q) = (\widetilde q-q)u(\widetilde q) + q ( u(\widetilde{q})-u(q) ) =: F_1 + F_2,\quad\mbox{in }\Omega.
\end{equation*}
Lemma \ref{lem:Sob-embed} implies $ u(q)\in H^2\II\hookrightarrow C(\overline{\Omega})$, and hence $F_1, F_2\in L^2(\Omega)$ with
\begin{align*}
\|F_1\|_{L^2(\Omega)} + \|F_2\|_{L^2(\Omega)} &\leq \|q-\widetilde q\|_{L^2(\Omega)}\|u(q)\|_{L^\infty(\Omega)} + \|q\|_{L^\infty(\Omega)}\|u(q) - u(\widetilde q)\|_{L^2(\Omega)} \\& \le c\|q-\widetilde q\|_{L^2(\Omega)}.
\end{align*}
This estimate and the governing equation fo $w$ immediately imply
\begin{equation}\label{eqn:est-q-H2L2}
    \|u(\widetilde q)-u(q)\|_{H^2(\Omega)}\leq c\|\widetilde q-q\|_{L^2(\Omega)},
\end{equation}
with the constant $c$ independent of $q$ and $\widetilde q$.
Meanwhile, the product rule yields the identity
\begin{align*}
    \nabla F_1 &= u(\widetilde q)\nabla (\widetilde q-q) + (\widetilde{q}-q)\nabla u(\widetilde q), \\
    \nabla F_2 & = ( u(\widetilde{q})-u(q) )\nabla q + q\nabla ( u(\widetilde{q})-u(q) ).
\end{align*}
Then the triangle inequality, H\"{o}lder's inequality and the continuous embeddings $H^1(\Omega)\hookrightarrow L^4(\Omega)$ and $H^2(\Omega)\hookrightarrow C(\overline{\Omega})$ from Lemma \ref{lem:Sob-embed} imply
\begin{align*}
   \|\nabla F_1\|_{L^2\II} & \leq \|\nabla (q-\widetilde q)u(\widetilde{q})\|_{L^2\II} + \|(q-\widetilde q)\nabla u(\widetilde{q})\|_{L^2\II} \\
    &\leq \|\nabla (q-\widetilde q)\|_{L^2\II}\|u(\widetilde{q})\|_{L^\infty(\Omega)} + \|q-\widetilde{q}\|_{L^4(\Omega)} \|\nabla u(\widetilde{ q})\|_{L^4\II}  \\
    & \leq c\|q-\widetilde{q}\|_{H^1(\Omega)}\|u(\widetilde{q})\|_{H^2(\Omega)}  \leq c \|q-\widetilde{q}\|_{H^1\II},
\end{align*}
with the constant $c$ independent of $\widetilde q$, in view of the first assertion. Moreover, by the estimate \eqref{eqn:est-q-H2L2}, we also have the following estimate
\begin{align*}
    \|\nabla F_2\|_{L^2(\Omega)} &\leq \|u(\widetilde{q})-u(q)\|_{L^\infty(\Omega)}\|\nabla q\|_{L^2(\Omega)} + \|q\|_{L^\infty(\Omega)} \|\nabla (u(\widetilde q) - u(q))\|_{L^2(\Omega)} \\&\leq c(\|\nabla q\|_{L^2(\Omega)} + 1) \|u(\widetilde q) - u(q)\|_{H^2\II} \leq c(\|\nabla q\|_{L^2(\Omega)} + 1)\|q-\widetilde{q}\|_{H^1\II}.
\end{align*}
The preceding estimates and the full regularity pickup of the elliptic operator $-\Delta$ (with a zero Dirichlet boundary) on a smooth domain (see, e.g., \cite[Theorem 9.25]{Brezis2010FunctionalAS}) imply the desired estimate \eqref{eqn:Lipschtz-cont}.
\end{proof}

\begin{remark}\label{rmk:neg-stab-est}The regularity estimates in Lemma \ref{lem:reg-u} play a crucial role in the subsequent convergence analysis, e.g., bounding the entropy integral \eqref{eqn:entropy} in the proof of Lemma \ref{lem:approx-sol-uppbound-q}. Note that Nickl et al \cite[p. 396, Section 4.4]{NicklWang:2020} proved the following stability estimate
\begin{equation}\label{eqn:ell-neg}
    \|u(q)-u(\widetilde{q})\|_{L^2(\Omega)}\leq \Big(1+\max\big(\|q\|_{H^\alpha(\Omega)}^4, \|\widetilde{q}\|_{H^\alpha(\Omega)}^4\big)\Big)\|q-\widetilde q\|_{(H^2(\Omega))^*},
\end{equation}
where $(H^2(\Omega))^*$ denotes the dual space of $H^2(\Omega)$. The estimate \eqref{eqn:ell-neg} was derived using a duality argument, and requires the condition $q,\widetilde q\in H^\alpha(\Omega)$ with $\alpha>\max({d}/{2}+2,2d-2)$. It is essential for applying the abstract framework in \cite[Section 2]{NicklWang:2020} {\rm(}see also \cite[eq. (2.4), p. 24]{Nickl:2023}{\rm)}. In contrast, the convergence result in Theorem \ref{thm:con-err-potential} only uses the estimate \eqref{eqn:Lipschtz-cont}, and hence can deal with a low-regularity potential coefficient $q\in \mathcal{A}\cap H^1(\Omega)$. Formally, the estimates \eqref{eqn:Lipschtz-cont} and \eqref{eqn:ell-neg} both indicate the degree of ill-posedness of the inverse potential problem: the forward operator is two-times smoothing.
\end{remark}

The next result gives a condition stability estimate for the inverse potential problem. This result improves \cite[Theorem 2.1]{jin2022convergence} with a more explicit dependence on $q$ and $\widetilde q$, and serves as the benchmark for the convergence rate of the regularized approximations.
\begin{proposition}\label{prop:cond-stab}
For any $q,\widetilde q\in \mathcal{A}\cap H^1(\Omega)$, the following estimate holds
\begin{equation*}
    \|(q-\widetilde q)u(q)\|_{L^2(\Omega)} \leq c(1+\max(\|q\|_{H^1(\Omega)},\|\widetilde q\|_{H^1(\Omega)}))^\frac{5}{6}\|u(q)-u(\widetilde{q})\|_{L^2(\Omega)}^\frac13.
\end{equation*}
\end{proposition}
\begin{proof}
By the weak formulations of $u\equiv u(q)$ and $\widetilde u\equiv u(\widetilde q)$, it follows that for any $\varphi\in H_0^1(\Omega)$,
\begin{equation*}
		((q-\widetilde q)u,\varphi) = -(\nabla(u-\widetilde u),\nabla\varphi)-(\widetilde q(u-\widetilde u(q)),\varphi)=:{\rm I} + {\rm II}.
	\end{equation*}
Let $\varphi=(q-\widetilde q)u$. The box constraint of $q$ and $\widetilde q$ yields $\|\varphi\|_{L^2(\Omega)}\leq c$. Moreover, we have
$\nabla\varphi=(\nabla q-\nabla \widetilde q)u+(q-\widetilde q)\nabla u$. This, the regularity estimate $\|u\|_{H^2(\Omega)}\leq c$ from Lemma \ref{lem:reg-u}  and the Sobolev embedding $H^2(\Omega)\hookrightarrow L^\infty(\Omega)$ from Lemma \ref{lem:Sob-embed} imply	
\begin{align*}
\|\nabla\varphi\|_{L^2\II}&\leq \|\nabla (q- \widetilde q)\|_{L^2\II}\|u\|_{L^\infty\II}+\|q-\widetilde q\|_{L^\infty\II}\|\nabla u\|_{L^2\II} \\& \leq c(1+\max(\|q\|_{H^1(\Omega)}, \|\widetilde q\|_{H^1(\Omega)}) ).
\end{align*}
Hence $\varphi\in H_0^1\II$. By the Cauchy-Schwarz inequality, the box constraint on $q$ and $\widetilde q$, and Lemmas \ref{lem:connect-seminorm-stdnorm} and \ref{lem:reg-u}, we obtain
\begin{align*}
|{\rm II}|&\leq c\|u-\widetilde u\|_{L^2\II}\|\widetilde q\|_{L^\infty\II}\|\varphi\|_{L^2\II} \leq c\|u-\widetilde u\|_{L^2\II}.
\end{align*}
The Gagliardo--Nirenberg interpolation inequality \cite{BrezisMironescu:2018},
\begin{equation}\label{eqn:Gargliardo}
\|v\|_{H^1(\Omega)}\leq c\|v\|_{L^2(\Omega)}^{\frac23}\|v\|_{H^3(\Omega)}^{\frac13},\quad\forall v\in H^3(\Omega),
\end{equation}
 Lemmas \ref{lem:connect-seminorm-stdnorm} and \ref{lem:reg-u}, we get
 \begin{equation*}
    \begin{split}
       \|\nabla\big(u-\widetilde u \big)\|_{L^2\II} &\leq c\|u-\widetilde u\|^\frac23_{L^2\II}\|u-\widetilde u\|^{\frac13}_{H^3\II}\\& \leq c\|u-\widetilde u\|^\frac23_{L^2\II}((1+\|q\|_{H^1(\Omega)})\|q-\widetilde q\|_{H^1\II})^{\frac13}\\
       & \leq c(1+\max(\|q\|_{H^1(\Omega)},\|\widetilde q\|_{H^1(\Omega)}))^{\frac23}\|u-\widetilde u\|^\frac23_{L^2\II}.
    \end{split}	
    \end{equation*}
Together with the Cauchy-Schwarz inequality, we arrive at
    \begin{align*}
    	|{\rm I}| &\leq \|\nabla\big(u^\dag-u(q^*)\big)\|_{L^2\II}\|\nabla\varphi\|_{L^2\II}\\ &
        \leq c(1+\max(\|q\|_{H^1(\Omega)},\|\widetilde q\|_{H^1(\Omega)}))^{\frac53}\|u-\widetilde u\|^\frac23_{L^2\II}.
\end{align*}
Combining the preceding two estimates yields the desired assertion.
\end{proof}

The next lemma gives the sub-Gaussianity of $\{\big(\xi,u(q)-u^\dagger\big)_n, q\in \mathcal{A}\}$.  Let $u^\dag=u(q^\dag)$.
\begin{lemma}\label{lem:con-subGaussian}
Let $\rbd(q,\widetilde{q}) = \sigma n^{-\frac12}\|u(q)-u(\widetilde{q})\|_n$. Then $\{\big(\xi,u(q)-u^\dagger\big)_n, q\in \mathcal{A}\}$ is a sub-Gaussian random process equipped with the semimetric $\rbd$.
\end{lemma}
\begin{proof}
The proof follows similarly as \cite[Lemma 2.2.7]{vanDerVaart:1996} or \cite[Lemma 4.6]{chen2018stochastic}, and we only include it for completeness. First, the well-definedness of $\rbd$ follows from the regularity $u(q),u(\widetilde{q})\in H^2\II$ in Lemma \ref{lem:reg-u} and the continuous embedding $H^2(\Omega)\hookrightarrow C(\overline{\Omega})$ from Lemma \ref{lem:Sob-embed}. The non-negativity, symmetry and identity hold by the definition of $\rbd$. Since $\|\cdot\|_n$ is a semi-norm over $C(\overline{\Omega})$, the triangle inequality holds for any $q,\widetilde{q},q'\in\mathcal{A}$:
\begin{align*}
\rbd(q,\widetilde{q}) &= \sigma n^{-\frac12}\|u(q)-u(\widetilde{q})\|_n  \leq \sigma n^{-\frac12}\|u(q)-u(q')\|_n + \sigma n^{-\frac12}\|u(q')-u(\widetilde{q})\|_n\\
&=\rbd(q,q') + \rbd(q',\widetilde{q}).
\end{align*}
Thus, $\rbd$ is a semimetric. To prove the sub-Gaussianity of $\{\big(\xi,u(q)-u^\dagger\big)_n, q\in \mathcal{A}\}$, let $\mathbb{Z}_q:= \big(\xi,u(q)-u^\dagger\big)_n$. Then
    \begin{equation*}
        \mathbb{Z}_q - \mathbb{Z}_{\widetilde{q}} = \big(\xi,u(q)-u(\widetilde{q})\big)_n = n^{-1}\sum_{j=1}^n\xi_j\big( u(q)-u(\widetilde{q})\big)(x_j).
    \end{equation*}
By Definition \ref{def:sub-Gau}, we deduce
\begin{align*} \mathbb{E}\big[\exp\big(\lambda(\mathbb{Z}_q - \mathbb{Z}_{\widetilde{q}})  \big)\big] & = \prod_{j=1}^n \mathbb{E}\big[\exp(\lambda n^{-1} \xi_j\big( u(q)-u(\widetilde{q})\big)(x_j))\big] \\
& \leq \exp\Big(\tfrac12\sigma^2\lambda^2n^{-2}\sum_{j=1}^n |\big( u(q)-u(\widetilde{q})\big)(x_j))|^2\Big)\\&  = \exp\Big(\tfrac12\sigma^2\lambda^2n^{-1}\|u(q)-u(\widetilde{q})\|^2_n\Big),
\end{align*}
since $(\xi_j)_{j=1}^n$ are i.i.d. sub-Gaussian variables with parameter $\sigma$ and zero expectation. Thus, by definition, the desired assertion follows.
\end{proof}

Next we give an important \textit{a priori} $\psi_2$-Orlicz bound on  $\|u^\dagger-u(q^*)\|_n$ and $\| q^*\|_{H^1(\Omega)}$.
\begin{lemma}\label{lem:approx-sol-uppbound-q} Let Assumptions \ref{ass:quasi-points} and \ref{ass:con-reg} be fulfilled, and let $q^*$ be a minimizer to problem \eqref{eqn:conti-optim-prob-potential}-\eqref{eqn:conti-weak-potential}. Then with $\rho_0 = \|q^\dagger\|_{H^1\II} + \sigma n^{-\frac12}$ and $\gamma^{{\frac12 + \frac{d}{12}}} = O(\sigma n^{-\frac12}(1+\rho_0)^{\frac d6} \rho_0^{-1})$, the following $\psi_2$-Orlicz bound holds
	\begin{equation*}
		\bigg\|\|u(q^*)-u^\dagger\|_n\bigg\|_{\psi_2} + \gamma^{\frac12}\bigg\|\| q^*\|_{H^1\II}\bigg\|_{\psi_2} \leq c\gamma^{\frac12}\rho_0.
	\end{equation*}
\end{lemma}
\begin{proof}
By the minimizing property of $q^* \in \mathcal{A}$ to the functional $J_\gamma(q) $, it follows that
	\begin{equation*}
		\|u(q^*)-\bm m\|_{n}^2 + \gamma \|q^*\|_{H^1\II}^2\leq \|u^\dagger-\bm m\|_{n}^2 + \gamma\| q^\dagger\|_{H^1\II}^2.
	\end{equation*}
	Then the definition of $\rho_0$ gives the basic inequality
	\begin{align}\label{eqn:min-prop}
		\|u(q^*)-u^\dagger\|_n^2 + \gamma \|q^*\|_{H^1\II}^2 & \leq 2\big(\bm\xi,u(q^*) - u^\dagger\big)_n +  \gamma\rho_0^2.
	\end{align}
Now we prove the desired estimate by a peeling argument \cite{vanDerVaart:1996}. Let $\delta>0$ and $\rho>0$ be two constants to be determined later, and for $i,j\ge 1$, let
\begin{equation*}
  A_0=[0,\delta), ~~~ A_i=[2^{i-1}\delta,2 ^i\delta), ~~~ B_0=[0,\rho), ~~~ B_j=[2^{j-1}\rho,2^j\rho)\,.
\end{equation*}
For $i,j\ge 0$, we define the set
\begin{equation*}
		F_{ij}: = \{q \in \mathcal{A}: \|u(q)-u^\dagger\|_n \in A_i, ~~~ {\|q\|_{H^1\II}} \in B_j \}.
	\end{equation*}
	Hence, the following disjoint decompositions hold
	\begin{equation*}
		\big\{\|u(q)-u^\dag\|_n>\delta\big\} = \bigcup_{i=1}^\infty\bigcup_{j=0}^\infty F_{ij}\quad\mbox{and}\quad\big\{\| q\|_{H^1\II}>\rho\big\} = \bigcup_{i=0}^\infty\bigcup_{j=1}^\infty F_{ij}.
	\end{equation*}
Next we prove the first estimate of the lemma. By the subadditivity of probability, we have
\begin{equation*}
	\mathbb{P}\big(\|u(q^*)-u^\dag\|_n>\delta\big)\leq\sum_{i=1}^\infty\sum_{j=0}^\infty \mathbb{P}(q^*\in F_{ij}) = \sum_{i=1}^\infty\sum_{j=1}^\infty \mathbb{P}(q^*\in F_{ij}) + \sum_{i=1}^\infty\mathbb{P}(q^*\in F_{i0}).
\end{equation*}
It suffices to estimate the probability $\mathbb{P}(q^*\in F_{ij})$ of the event $\{q^*\in F_{ij}\}$. First consider the case  $i,j\geq1$. The condition $q^*\in F_{ij}$ and the inequality \eqref{eqn:min-prop} imply that the following event
\begin{equation*}
		2^{2(i-1)}\delta^2 + 2^{2(j-1)}\gamma\rho^2 \leq 2\sup_{q\in F_{ij}}|(\bm \xi, u(q)-u^\dagger)_n| + \gamma\rho_0^2
\end{equation*}
holds. Let $G_{i,j}=\cup_{k,l=0}^{i,j} F_{kl}$. Then using the semimetric $\rbd$ defined in Lemma \ref{lem:con-subGaussian} and the triangle inequality, we get
    \begin{align*}
    \diam(G_{ij}) := \sup_{q,\widetilde{q}\in G_{ij}}\rbd(q,\widetilde{q}) = \sup_{q,\widetilde{q}\in G_{ij}}\sigma n^{-\frac12}\|u(q)-u(\widetilde{q})\|_n  \leq \sigma n^{-\frac12}2^{i+1}\delta.
    \end{align*}
Define the set $u(G_{ij})=\{u(q)-u(q^\dagger): q\in G_{ij}\}$. Then we have the following identity
\begin{equation}\label{eqn:cover-id}
    N\left(\tfrac\epsilon2, G_{ij},\rbd\right) = N\left(\tfrac{\epsilon}{2\sigma n^{-\frac12}}, u(G_{ij}),\|\cdot\|_n\right).
\end{equation}
To see the identity, let $B_r(q,\rbd)$ be a ball centered at $q$ with a radius $r$ in the semi-distance $\rbd$. Then by definition of $N\left(\tfrac\epsilon2, G_{ij},\rbd\right)$ (cf. Definition \ref{def:cov-num}), there exists a sequence $(q_i)_{i=1}^{N\left(\tfrac\epsilon2, G_{ij},\rbd\right)}$ such that
\begin{equation*}
    G_{ij}\subset \bigcup_{i=1}^{N\left(\epsilon/2, G_{ij},\rbd\right)} B_{\tfrac\epsilon 2}(q_i,\rbd).
\end{equation*}
By the definitions of the semimetrics $\rbd$ and $\|\cdot\|_n$, we have
\begin{equation*}
    u(G_{ij}) \subset \bigcup_{i=1}^{N\left(\epsilon/2, G_{ij},\rbd\right)} B_{\tfrac{\epsilon}{2\sigma n^{-\frac12}}}(u(q_i),\|\cdot\|_n).
\end{equation*}
This implies that the inequality $$ N\left(\tfrac\epsilon2, G_{ij},\rbd\right)\geq N\left(\tfrac{\epsilon}{2\sigma n^{-\frac12}}, u(G_{ij}),\|\cdot\|_n\right).$$ Similarly, we can prove the converse inequality and thus the desired identity \eqref{eqn:cover-id} follows.
Let $Z_q = (\bm \xi, u(q)-u^\dagger)_n$. Then by definition \ref{def:sub-Gau-process}, Lemmas \ref{lem:max-ineq} and \ref{lem:con-subGaussian}, we get
\begin{align}
&\bigg\|\sup_{q\in F_{ij}}|(\bm \xi, u(q)-u(q^\dagger))_n|\bigg\|_{\psi_2} \leq \bigg\|\sup_{q\in G_{ij}}|(\bm \xi, u(q)-u(q^\dagger))_n|\bigg\|_{\psi_2}\nonumber\\
\leq& \bigg\|\sup_{q,\Tilde{q}\in G_{ij}}|(\bm \xi, u(q)-u(\Tilde{q}))_n|\bigg\|_{\psi_2}
 \leq c\int_0^{\sigma n^{-\frac12}2^{i+1}\delta} \sqrt{\log N\left(\tfrac\epsilon2, G_{ij},\rbd\right)}~{\rm d}\epsilon\nonumber \\
	 =& c\int_0^{\sigma n^{-\frac12}2^{i+1}\delta} \sqrt{\log N\left(\tfrac{\epsilon}{2\sigma n^{-\frac12}}, u(G_{ij}),\|\cdot\|_n\right)}~{\rm d}\epsilon\nonumber \\
	 \leq& c\int_0^{\sigma n^{-\frac12}2^{i+1}\delta} \sqrt{\log N\left(\tfrac{\epsilon}{2c\sigma n^{-\frac12}}, u(G_{ij}),\|\cdot\|_{L^\infty\II}\right)}~{\rm d}\epsilon.\label{eqn:maximal-cont}
\end{align}
The continuity estimate in Lemma \ref{lem:reg-u} gives
\begin{equation*}
\|u(q)-u(q^\dagger)\|_{H^3\II} \leq c(1+\|q^\dagger\|_{H^1(\Omega)})\|q-q^\dagger\|_{H^1(\Omega)} \leq c2^{j+1}\rho(1+\rho_0).
\end{equation*}
Thus the set $u(G_{ij})$ is contained in an $H^3\II$ ball of radius $r=c2^j\rho(1+\rho_0)$ (by absorbing one factor of $2$ into $c$). By Lemma \ref{lem:cover-en-bound-Sobolev} (with the choice $(s,p,q)=(3,2,\infty)$), we deduce
\begin{equation}\label{eqn:entropy}
 	\log N\left(\tfrac{\epsilon}{2c\sigma n^{-\frac12}}, u(F_{ij}),\|\cdot\|_{L^\infty\II}\right) \leq c\bigg(\frac{\sigma n^{-\frac12}2^j\rho(1+\rho_0)}{\epsilon}\bigg)^{\frac d3}.
\end{equation}
Combining the estimates \eqref{eqn:maximal-cont} and \eqref{eqn:entropy} gives
\begin{align*}
	\bigg\|\sup_{q\in F_{ij}}|(\bm \xi, u(q)-u^\dagger)_n|\bigg\|_{\psi_2} & \leq c\int_0^{\sigma n^{-\frac12}2^{i+1}\delta}\bigg(\frac{\sigma n^{-\frac12}2^j\rho(1+\rho_0)}{\epsilon}\bigg)^{\frac d6}~{\rm d}\epsilon \\
    & \leq c\sigma n^{-\frac12}(2^i\delta)^{1-\frac d6}(2^j\rho)^{\frac d6}{(1+\rho_0)}^{\frac d6}.
 \end{align*}
By Lemmas \ref{lem:con-subGaussian} and \ref{lem:bdd-sub-Gaussian}, we have for $i,j\geq 1$,
\begin{align*}
	\mathbb{P}(q^*\in F_{ij}) &\leq\mathbb{P}\bigg(2\sup_{q\in F_{ij}}|(\bm \xi, u(q)-u^\dagger)_n| \geq 2^{2(i-1)}\delta^2 + 2^{2(j-1)}\gamma\rho^2 - \gamma\rho_0^2 \bigg) \\
		& \leq 2\exp\bigg(-c\sigma^{-2} n \bigg|\frac{2^{2(i-1)}\delta^2 + 2^{2(j-1)}\gamma\rho^2 - \gamma\rho_0^2}{(2^i\delta)^{1-\frac d6}(2^j\rho)^{\frac d6}{(1+\rho_0)}^{\frac d6}}\bigg|^2\bigg) := {\rm I}_{ij}.
\end{align*}
Now we choose $\delta^2=\gamma\rho_0^2 (1+z)^2$ and $\rho = \rho_0$ for $z\geq1$. By Young's inequality $ab\leq \frac{a^p}{p}+\frac{b^q}{q}$ for $a,b\geq0$ and $p,q> 1$ with $\frac1p+\frac1q=1$, we have $$(2^i(1+z))^{1-\frac d6}(2^j)^{\frac d6} \leq c\big(2^i(1+z) + 2^j\big).$$
Using this inequality  and the \textit{a priori} choice $\gamma^{\frac12 + \frac{d}{12}} = O(\sigma n^{-\frac12}(1+\rho_0)^{\frac d6}\rho_0^{-1})$, we obtain
\begin{equation*}
   {\rm I}_{ij} = 2\exp\bigg(-c\bigg|\frac{2^{2i}z(1+z) + 2^{2j}}{(2^i(1+z))^{1-\frac d6}(2^j)^{\frac d6}}\bigg|^2\bigg) \leq 2\exp\bigg(-c\big(2^{2i}z^2 + 2^{2j}\big)\bigg).
\end{equation*}
Similarly, the following estimate holds for $j=0$:
    \begin{align*}
    	\mathbb{P}(q^*\in F_{i0}) \leq 2\exp ({- c2^{2i} z^2}).
    \end{align*}
The estimates for the cases $i,j\geq1$ and $j=0$ and the inequality $\exp(-c)<1$ yield
\begin{align*}
	&\mathbb{P}\big(\|u(q^*)-u^\dagger\|_n>\gamma^{\frac12}\rho_0(1+z)\big)\\
 \leq& \sum_{i=1}^\infty\sum_{j=1}^\infty2\exp\big(-c\big(2^{2i}z^2 + 2^{2j}\big)\big) + \sum_{i=1}^\infty2\exp(- c2^{2i} z^2) \leq 4\exp (- cz^2),
\end{align*}
This and Lemma \ref{lem:Orlicz-bound} complete the proof of the $\psi_2$-Orlicz bound on the prediction error $\|u(q^*)-u^\dag\|_n$. By interchanging the summation range of $i$ and $j$, we can get the bound on $\|q^*\|_{H^1\II}$ similarly.
\end{proof}

\begin{remark}\label{rmk:rate-n}
By combining Lemma \ref{lem:approx-sol-uppbound-q} with Lemmas \ref{lem:bdd-sub-Gaussian} and \ref{lem:Orlicz-bound}, we deduce that the following estimates in high probability hold for $z>0$,
\begin{equation*}
\mathbb{P}(\|u(q^*)-u^\dagger\|_n>\gamma^{\frac12}\rho_0z) \leq 2\exp(-cz^2)\quad\mbox{and}\quad \mathbb{P}(\| q^*\|_{H^1\II}>\rho_0z) \leq 2\exp(-cz^2).
\end{equation*}
This result implies a prediction error $\|u(q^*)-u^\dag\|_n$ of $O(\gamma^\frac12\rho_0)$ in high probability. {\color{blue}Furthermore, we also have $\mathbb{E} [\|u(q^*)-u^\dagger\|_n^2]\leq c\gamma\rho_0^2$ and $\mathbb{E}[\| q^*\|_{H^1\II}^2]\leq c\rho_0^2$.} With the a priori choice $\gamma^{\frac12+\frac{d}{12}} = O(\sigma n^{-\frac12}{(1+\rho_0)}^{\frac d6}\rho_0^{-1})$, we deduce a bound $O((\sigma n^{-\frac12})^{\frac{6}{6+d}}(1+\rho_0)^\frac{2d}{6+d})$ on the prediction error $\|u(q^*)-u^\dag\|_n$ in terms of the number $n$ of sampling points. Under the a priori Sobolev regularity $u^\dag \in H^3(\Omega)$ from Lemma \ref{lem:reg-u}, the obtained rate coincides with the rate $O(n^{-\frac{3}{6+d}})$ for the nonparametric regression \cite[p. 188]{geer2000empirical}. This rate is known to be minmax optimal in the class of H\"{o}lder continuous regression functions \cite{Stone:1982}.
\end{remark}
\begin{remark}
In the proof of Lemma \ref{lem:approx-sol-uppbound-q}, we have made full use of the elliptic regularity $u(q)\in H^3(\Omega)$ of the state $u(q)$ under Assumption \ref{ass:con-reg} to bound the covering number, cf. \eqref{eqn:entropy}. This regularity holds only for the continuous problem \eqref{eqn:potential} on a smooth domain $\Omega$. On a convex polyhedral domain, we can only expect the regularity $u(q)\in H^2(\Omega)$. Then by repeating the argument of Lemma \ref{lem:approx-sol-uppbound-q}, with $\rho_0 = \|q^\dagger\|_{H^1\II} + \sigma n^{-\frac12}$ and $\gamma^{{\frac12 + \frac{d}{8}}} = O(\sigma n^{-\frac12}\rho_0^{-1})$,
the following  $\psi_2$-Orlicz bound holds
\begin{equation*}
\bigg\|\|u(q^*)-u^\dagger\|_n\bigg\|_{\psi_2} + \gamma^{\frac12}\bigg\|\| q^*\|_{H^1\II}\bigg\|_{\psi_2} \leq c\gamma^{\frac12}\rho_0,
\end{equation*}
or by Lemmas \ref{lem:bdd-sub-Gaussian} and \ref{lem:Orlicz-bound} equivalently the following estimates hold for $z>0$,
$\mathbb{P}(\|u(q^*)-u^\dagger\|_n>\gamma^{\frac12}\rho_0z) \leq 2\exp(-cz^2)$ and $\mathbb{P}(\|q^*\|_{H^1\II}>\rho_0z) \leq 2\exp(-cz^2)$.
\end{remark}

 Now we can give the proof of Theorem \ref{thm:con-err-potential}
\begin{proof}[Proof of Theorem \ref{thm:con-err-potential}]
By Lemma \ref{lem:approx-sol-uppbound-q} and the estimate \eqref{ineq:proba-distrib-estimate}, there hold with probability at least $1-2\tau$
\begin{equation*}
\|u(q^*)-u^\dagger\|_n \leq c\gamma^{\frac12}\rho_0\ell_\tau \quad \mbox{and}\quad \|\nabla q^*\|_{L^2\II} \leq c\rho_0\ell_\tau.
\end{equation*}
Next from the weak formulations of $u^\dag\equiv u(q^\dag)$ and $u(q^*)$, it follows that for any $\varphi\in H_0^1(\Omega)$,
\begin{equation*}
		((q^\dagger-q^*)u^\dagger,\varphi) = -(\nabla(u^\dagger-u(q^*)),\nabla\varphi)-(q^*(u^\dagger-u(q^*)),\varphi)=:{\rm I} + {\rm II}.
	\end{equation*}
Let $\varphi=(q^\dagger-q^*)u^\dagger$. Assumption \ref{ass:con-reg} and the box constraint of $q^\dagger$ and $q^*$ yield $\|\varphi\|_{L^2(\Omega)}\leq c$. Moreover, we have
$\nabla\varphi=(\nabla q^\dag-\nabla q^*)u^\dag+(q^\dag-q^*)\nabla u^\dag$.
Consequently,  	
\begin{align*}
\|\nabla\varphi\|_{L^2\II}&\leq \|\nabla (q^\dag- q^*)\|_{L^2\II}\|u^\dag\|_{L^\infty\II}+\|q^\dag-q^*\|_{L^\infty\II}\|\nabla u^\dag\|_{L^2\II} \\& \leq c(\|\nabla q^\dag\|_{L^2\II} + \|\nabla q^*\|_{L^2\II}  + 1)\leq c\big(1 + \rho_0\ell_\tau \big).
\end{align*}
Hence $\varphi\in H_0^1\II$. By the Cauchy-Schwarz inequality, the box constraint on $q^\dag$ and $q^*$, and Lemmas \ref{lem:connect-seminorm-stdnorm} and \ref{lem:reg-u}, we obtain
\begin{align*}
|{\rm II}|&\leq c\|u^\dag-u(q^*)\|_{L^2\II}\|q^*\|_{L^\infty\II}\|\varphi\|_{L^2\II} \\& \leq c(\|u^\dag-u(q^*)\|_n + n^{-\frac3d} \|u^\dag-u(q^*)\|_{H^3(\Omega)}) \\&
\leq c(\|u^\dag-u(q^*)\|_n + n^{-\frac3d}(1+\rho_0) \|q^\dag-q^*\|_{H^1(\Omega)})
\\&\leq c \big(\gamma^{\frac12} + n^{-\frac 3d}(1+\rho_0) \big)\big(1+\rho_0\ell_\tau).
\end{align*}
To bound the term ${\rm I}$, by the Gagliardo--Nirenberg interpolation inequality \cite{BrezisMironescu:2018} (cf. the inequality \eqref{eqn:Gargliardo}), Lemmas \ref{lem:connect-seminorm-stdnorm} and \ref{lem:reg-u}, we get
    \begin{equation*}
    \begin{split}
       \|\nabla\big(u^\dag-u(q^*)\big)\|_{L^2\II} &\leq \|u^\dag-u(q^*)\|^\frac23_{L^2\II}\|u^\dag-u(q^*)\|^{\frac13}_{H^3\II}\\& \leq c\big[\gamma^{\frac12} + n^{-\frac 3d}(1+\rho_0) \big(1+\rho_0\ell_\tau)\big]^{\frac23}(1+\rho_0)^{\frac 13} \|q^\dag - q^*\|_{H^1(\Omega)}^{\frac13} \\& \leq c\big(\gamma^{\frac12} + n^{-\frac 3d} \big)^{\frac23}(1+\rho_0\ell_\tau)^2.
    \end{split}	
    \end{equation*}
Together with the Cauchy-Schwarz inequality, we arrive at
    \begin{align*}
    	|{\rm I}| \leq \|\nabla\big(u^\dag-u(q^*)\big)\|_{L^2\II}\|\nabla\varphi\|_{L^2\II} \leq c\big(\gamma^{\frac12} + n^{-\frac 3d} \big)^{\frac23}(1+\rho_0\ell_\tau)^3.
\end{align*}
This proves the first assertion.  To prove the second estimate, we repeat the argument of \cite[Theorem 2.1]{jin2022convergence}. Specifically, we decompose the domain $\Omega$ into two disjoint sets $\Omega=\Omega_\rho\cup \Omega_\rho^c$, with $\Omega_\rho=\{x\in\Omega:{\rm dist}(x,\partial\Omega)\geq\rho\}$ and $\Omega_\rho^c=\Omega\setminus\Omega_\rho$, with the constant $\rho>0$ to be chosen. On the subdomain $\Omega_\rho$, Condition \ref{Cond:positivity} leads to
   \begin{align*}
     \int_{\Omega_\rho}(q^\dag-q^*)^2{\rm d}x &\leq \rho^{-2\beta}\int_{\Omega_\rho }(q^\dag-q^*)^2\mathrm{dist}(x,\partial\Omega)^{2\beta}{\rm d}x \leq c\rho^{-2\beta}\|(q^{\dag}-q^*)u^\dag\|^2_{L^2\II},
   \end{align*}
Meanwhile, since $q^\dag,q^*\in\mathcal{A}$, we have
$$\int_{\Omega_\rho^c}(q^{\dag}-q^*)^2{\rm d}x \leq c|\Omega_\rho^c|\leq c\rho,$$
where $|\Omega_\rho^c|$ denotes the Lebesgue measure of the set $\Omega_\rho^c$. By setting
$\rho = c\rho^{-2\beta}\|(q^{\dag}-q^*)u^\dag\|^2_{L^2\II}$, i.e., $\rho=c\|(q^*-q^\dag)u^\dag\|_{L^2(\Omega)}^\frac{2}{1+2\beta}$, we can balance the last two estimates and obtain the following estimate
$$\|q^\dag-q^*\|_{L^2\II} \leq c\|(q^{\dag}-q^*)u^\dag\|_{L^2\II}^{\frac{1}{1+2\beta}}.$$
Now using the first assertion, we arrive at the second assertion.
\end{proof}

\begin{remark}
The overall proof of Theorem \ref{thm:con-err-potential} heavily relies on the test function $\varphi = (q^\dag - q^*)u^\dag \in H_0^1(\Omega)$. This choice necessitates an a priori $H^1(\Omega)$ regularity on the regularized solution $q^*$, which motivates the construction of the regularized model \eqref{eqn:conti-optim-prob-potential}. Note that in principle one may employ an $H^s(\Omega)$ penalty, with $s>1$, which however is numerically inconvenient to approximate using the Galerkin FEM, at least when $s$ is not an integer.
\end{remark}

\section{Proof of Theorem \ref{thm:err-potential}}\label{sec:error-disc}

In this section, we provide a stochastic convergence rate analysis of the discrete approximation $q_h^*$, which differs markedly from that for the continuous problem in Section \ref{sec:error-cont}, due to limited Sobolev regularity of the FEM approximation. First we state the existence of a minimizer $q_h^*$, and give the proof in Appendix \ref{app:tech-proof}.

\begin{proposition}\label{prop:exist-disc}
The functional $J_{\gamma,h}(q_h)$ is well defined over the admissible set $\mathcal{A}_h$, and there exists at least one global minimizer $q_h^*\in \mathcal{A}_h$ to problem \eqref{eqn:dis-optim-pure-dis}-\eqref{eqn:dis-weak-potential} $\mathbb{P}$ almost surely.
\end{proposition}

We first recall several estimates from the finite element analysis \cite{Brenner2002The}.
Let $\Pi_h$ be the Lagrange nodal interpolation operator on the space $V_h$. It satisfies the following error estimates \cite[Corollary 4.4.24]{Brenner2002The}:
\begin{align}\label{ineq:Pih-approx-inf}
\|v-\Pi_hv\|_{L^\infty(K)} &\leq ch^{2-\frac d2} \|v\|_{H^2(K)}, \quad \forall v\in H^2(K), \forall K\in\mathcal{T}_h,\\
\|v-\Pi_hv\|_{L^2\II} &\leq ch^2 \|v\|_{H^2\II}, \quad \forall v\in H^2\II. \label{ineq:Pih-approx-2}
\end{align}

Moreover, we define the $L^2\II$ projection operator $P_h:L^2\II\mapsto X_h$ by
\begin{equation*}
	(P_hv,v_h)=(v,v_h), \quad \forall v\in L^2\II,  v_h\in X_h.
\end{equation*}
Then for $1\leq p\leq \infty$ and $s=0,1,2$, $k=0,1$ with $k\leq s$ \cite{thomee2007galerkin}:
\begin{equation}\label{ineq:Ph-approx}
	\big\|v-P_hv\big\|_{W^{k,p}\II}\leq  ch^{s-k}\big\|v\big\|_{W^{s,p}\II}, \quad \forall v\in W^{s,p}\II \cap H_0^1\II.
\end{equation}

 Let $\mathbb{T}_K:=(x_i)_{i=1}^n\cap K$ for each $K\in\mathcal{T}_h$ and $\#\mathbb{T}_{K}$ be the cardinality of $\mathbb{T}_K$. The notation $|K|$ denotes the Lebesgue measure of the set $K$. The next lemma gives a bound on the cardinality $\#\mathbb{T}_K$.
\begin{lemma}\label{lem:cardinality-estimate}
    Let Assumption \ref{ass:quasi-points} hold and $\mathcal{T}_h$ be quasi-uniform. Then $\#\mathbb{T}_{K} = O(nh^d)$.
\end{lemma}
\begin{proof}
Since the sampling points $(x_i)_{i=1}^n$  satisfy Assumption \ref{ass:quasi-points}, we can connect them suitably and generate a new quasi-uniform triangulation $\mathcal{T}_n = \cup {K_n}$ of the domain $\Omega$, with a mesh size $O(n^{-\frac1d})$. Hence, for any $K_n\in\mathcal{T}_n$ and $K\in\mathcal{T}_h$, we have
    \begin{equation*}
     |K|/|K_n| = O(h^d/ (n^{-\frac1d})^d) = O(nh^d).
    \end{equation*}
That is, each $K$ can contain $O(nh^d)$ element of $K_n$.
Since each element $K_n\in\mathcal{T}_n$ has only $d+1$ nodes (in $\mathbb{R}^d$), the definition $\mathbb{T}_K:=(x_i)_{i=1}^n\cap K$ implies $\#\mathbb{T}_K = O(nh^d)$.
\end{proof}

\begin{lemma}\label{lem:pointwise-L2-stab}
Let Assumption \ref{ass:quasi-points} be fulfilled. Then there exist $c,c'>0$, independent of $n$ and $h$, such that
	\begin{equation*}
		c'\|v_h\|_{L^2\II}\leq\|v_h\|_{n}\leq c\|v_h\|_{L^2\II},\quad\forall v_h\in V_h.
	\end{equation*}
\end{lemma}
\begin{proof}
 By Lemma \ref{lem:cardinality-estimate} and the definition of $\|\cdot\|_n$ and the inverse estimate $\|v_h\|_{L^\infty(K)}\leq ch^{-\frac{d}{2}}\|v_h\|_{L^2(K)}$ on the FEM space $V_h$ \cite[Theorem 4.5.11]{Brenner2002The}, there holds
\begin{align*}
\|v_h\|^2_{n} &= n^{-1}\sum_{i=1}^n |v_h(x_i)|^2 \leq n^{-1}\sum_{K\in\mathcal{T}_h} \#\mathbb{T}_K\|v_h\|^2_{L^\infty(K)}\\ &\leq c\sum_{K\in\mathcal{T}_h}\|v_h\|^2_{L^2(K)} = c\|v_h\|^2_{L^2\II}.
\end{align*}
Next, since $\#\mathbb{T}_{K} = O(nh^d)$, there holds
    \begin{align*}
        \|v_h\|_n^2 & = n^{-1}\sum_{i=1}^n |v_h(x_i)|^2 = \sum_{K\in\mathcal{T}_h} n^{-1}\sum_{x_i\in\mathbb{T}_K} |v_h(x_i)|^2 \\& = c\sum_{K\in\mathcal{T}_h} \Big[h^d (\#\mathbb{T}_{K})^{-1}\sum_{x_i\in\mathbb{T}_K} |v_h(x_i)|^2 \Big] =:  c\sum_{K\in\mathcal{T}_h} {\rm I}_K.
    \end{align*}
Meanwhile, let $\widehat{K}$ be the reference element and the set $(\widehat{x}_i)_{i=1}^{n_{\widehat{K}}}\subset \widehat{K}$ (with $n_{\widehat{K}}$ being the number of sampling points in the element $\widehat{K}$). Then by Lemma \ref{lem:connect-seminorm-stdnorm} and noting that $v_h$ is linear over $\widehat{K}$ (and $|v_h|_{H^2(\widehat{K})}=0$), we have
\begin{equation*}
c_1\|v_h\|_{L^2(\widehat{K})} \leq \|v_h\|_{n_{\widehat{K}}} \leq c_2\|v_h\|_{L^2(\widehat{K})}.
\end{equation*}
This estimate, the scaling argument and the fact $|K| = O(h^d)$ yield the desired estimate ${\rm I}_K \geq c\|v_h\|_{L^2(K)}^2$.
\end{proof}

The next result gives a crucial \textit{a priori} estimate on the FEM approximations $u_h(q)$ and $u_h(\Pi_hq^\dag)$ in the $\|\cdot\|_n$.
\begin{lemma}\label{lem:point-err}
Let Assumption \ref{ass:quasi-points} be fulfilled. Then for any $q\in\mathcal{A}$, there holds
	\begin{align*}
		\|u(q) - u_h(q)\|_n  \leq ch^2\|u(q)\|_{H^2(\Omega)}.
	\end{align*}
Additionally, if Assumption \ref{ass:reg} holds, then there holds
\begin{align*}
 \|u^\dag-u_h(\Pi_hq^{\dag})\|_{n} \leq ch^2,
\end{align*}
where the constant $c$ depends on $\|q^\dag\|_{H^2(\Omega)}$.
\end{lemma}
\begin{proof}
First we bound $\|u(q) - u_h(q)\|_n$.
By C\'{e}a's lemma and the standard duality argument \cite[Chapter 5]{Brenner2002The}, we obtain
\begin{equation}\label{eqn:approx-L2}
	\|u(q)-u_h(q)\|_{L^2\II} + h \|\nabla\big(u(q)-u_h(q)\big)\|_{L^2\II}  \leq ch^2\|u(q)\|_{H^2(\Omega)}.
\end{equation}
By the estimate \eqref{ineq:Pih-approx-inf} and $\#\mathbb{T}_{K} = O(nh^d)$ from Lemma \ref{lem:cardinality-estimate}, we have
\begin{align}
	\|u(q)-\Pi_hu(q)\|_n^2 & \leq n^{-1}\sum_{{K}\in\mathcal{T}_h}\#\mathbb{T}_{K}\|u(q)-\Pi_hu(q)\|_{L^\infty({K})}^2 \nonumber\\ & \leq cn^{-1}\sum_{{K}\in\mathcal{T}_h}\#\mathbb{T}_{K}h^{4-d}\|u(q)\|_{H^2({K})}^2\nonumber \\
    &\leq ch^4\|u(q)\|_{H^2\II}^2 . \label{eqn:approx-Pi-n}
\end{align}
Since $\Pi_hu(q) - u_h(q)\in X_h$,  Lemma \ref{lem:pointwise-L2-stab}, the triangle inequality and the estimate \eqref{ineq:Pih-approx-2} give
\begin{align*}
   & \|\Pi_hu(q) - u_h(q)\|_n \leq c\|\Pi_hu(q) - u_h(q)\|_{L^2\II}\\
    \leq &c\big(\|\Pi_hu(q) - u(q)\|_{L^2\II}+ \|u(q) - u_h(q)\|_{L^2\II}\big) \leq ch^2\|u(q)\|_{H^2(\Omega)}.
\end{align*}
Combining these estimates with the triangle inequality gives $$\|u(q)-u_h(q)\|_n\leq ch^2\|u(q)\|_{H^2(\Omega)}.$$ It remains to bound  $\|u^\dagger-u_h(\Pi_hq^{\dagger})\|_{n}$. By Assumption \ref{ass:reg} and \cite[Lemma 3.1]{jin2022convergence}, there holds
\begin{equation*}
	\|u^\dagger-u_h(\Pi_hq^{\dagger})\|_{L^2\II} + h \|\nabla (u^\dagger-u_h(\Pi_hq^{\dagger}))\|_{L^2\II}\leq  ch^2.
\end{equation*}
Repeating the proof for $\|u(q) - u_h(q)\|_n$ with Assumption \ref{ass:reg} yields $$\|u^\dagger-u_h(\Pi_hq^{\dagger})\|_{n} \leq ch^2.$$
This completes the proof of the lemma.
\end{proof}

\begin{remark}\label{h3-converge}
To  take the full advantage of the $H^3(\Omega)$ regularity in Lemma \ref{lem:reg-u}, one should use P2 finite elements instead of  P1 finite elements. Then the error bound in Lemma \ref{lem:point-err} can be improved to
$\|u(q) - u_h(q)\|_n  + \|u^\dagger-u_h(\Pi_hq^{\dagger})\|_{n} \leq ch^3$.
\end{remark}

To analyze the discrete approximation $q_h^*$, we recall a useful lemma due to van de Geer on stochastic convergence \cite[Lemma 8.4 and (10.6)]{geer2000empirical}. Using the terminology of the stochastic convergence order, we denote a random variable $X=\mathcal{O}_{p}(z)$ if $X$ is a sub-Gaussian random variable with zero expectation and parameter $z$.
 \begin{lemma}
\label{Geer-entropy}
Let $G$ be a function space with $\log N(\varepsilon, B_G,\|\cdot\|_{n}) \leqslant  A \varepsilon^{-r}$, with $N(\varepsilon, B_G,\|\cdot\|_{n})$ being the local covering number of the unit ball $B_G$. Then there holds
\begin{equation*}
    \sup_{g \in G} \frac{\left|(\bm{\xi}, g)_{n}\right|}{\|g\|_{n}^{1-\frac r2}\|g\|_{G}^{\frac r2}}=\mathcal{O}_{p}\left( \sqrt{A} \sigma n^{-\frac12}\right).
\end{equation*}
\end{lemma}

The next result gives a $\psi_2$-Orlicz bound, which is a discrete analogue of Lemma \ref{lem:approx-sol-uppbound-q}.
\begin{lemma}\label{lem:approx-sol-uppbound-qh}
Let Assumptions \ref{ass:quasi-points} and \ref{ass:reg} hold, and let $q_h^*\in\mathcal{A}_h$ be a minimizer to problem \eqref{eqn:dis-optim-pure-dis}-\eqref{eqn:dis-weak-potential}. Let $\rho_0 = \|q^\dagger\|_{H^1\II} + \sigma n^{-\frac12}$. Then with the a priori choice $\gamma^{{\frac12 + \frac {d} {12}}} = O(\sigma n^{-\frac12}{(1+\rho_0)^{\frac d6}}\rho_0^{-1})$, the following $\psi_2$-Orlicz bound holds
\begin{align*}
	&\bigg\|\big\|u_h(q_h^*) - u_h(\Pi_hq^\dagger)\big\|_n\bigg\|_{\psi_2} + \gamma^{\frac12}\bigg\|\big\|q_h^*\big\|_{H^1\II}\bigg\|_{\psi_2} \\
   \leq& c\big( \gamma^{\frac12} +h^2 + \gamma^{\frac{d}{12}}h^{4(\frac12-\frac{d}{8})}  + (\gamma^{\frac12+\frac{d}{12}}h^{4(\frac12-\frac{d}{8})})^{\frac12}
   \big) \rho_0.
\end{align*}
In particular, with $h=O(\gamma^{\frac{6-d}{6(4-d)}})$, the following estimate holds
  \begin{align*}
		\bigg\|\big\|u_h(q_h^*) - u_h(\Pi_hq^\dagger)\big\|_n\bigg\|_{\psi_2} + \gamma^{\frac12}\bigg\|\big\|q_h^*\big\|_{H^1\II}\bigg\|_{\psi_2}
   \leq c \gamma^{\frac12}  \rho_0.
		\end{align*}
\end{lemma}
\begin{proof}
It suffices to prove the first assertion, and we denote the $\psi_2$-Orlicz norm by $A$ below. By the minimizing property of $q_h^*\in\mathcal{A}_h$ and $\Pi_hq^\dagger\in\mathcal{A}_h$, there holds
\begin{equation*}
	\|u_h(q_h^*)-\bm m\|_{n}^2 + \gamma \|q_h^*\|_{H^1\II}^2\leq \|u_h(\Pi_hq^\dagger)-\bm m\|_{n}^2 + \gamma\|\Pi_hq^\dagger\|_{H^1\II}^2.
\end{equation*}
Let $u^0\equiv u(0)\in H^3\II\cap H^1_0\II$ (cf. Lemma \ref{lem:reg-u}) and $u_h^0 \equiv u_h(0)\in X_h$, and for any $q\in\mathcal{A}$, let $\overline{u}(q)=u(q)-u^0$ and $\overline{u}_h(q) = u_h(q)-u^0_h$. Direct computation leads to
\begin{align*}
	\|\overline{u}_h(q_h^*)-\overline{u}_h(\Pi_hq^\dagger)\|_n^2 + \gamma \| q_h^*\|_{H^1\II}^2 & \leq {\rm I} +\gamma \|\Pi_hq^\dag\|_{H^1(\Omega)}^2,
 \end{align*}
with ${\rm I}= 2(\bm{m}-u_h(\Pi_hq^\dag),\overline{u}_h(q_h^*)-\overline{u}_h(\Pi_hq^\dag))_n$.
This and the definition of the noisy data $\bm{m}$ in \eqref{eqn:rand-obs} give the following decomposition
 \begin{align*}
{\rm I} &=  2\big(\bm\xi,\overline{u}_h(q_h^*) - \overline{u}_h(\Pi_hq^\dagger)\big)_n +  2\big(\overline{u}^\dagger-\overline{u}_h(\Pi_hq^{\dagger}),\overline{u}_h(q_h^*) - \overline{u}_h(\Pi_hq^\dagger)\big)_n \\
	& \leq 2\big(\bm\xi,\overline{u}_h(q_h^*) - \overline{u}(q_h^*)\big)_n  + 2\big(\bm\xi,\overline{u}(q_h^*) - \overline{u}(\Pi_hq^\dagger)\big)_n+ 2\big(\bm\xi,\overline{u}(\Pi_hq^\dagger) - \overline{u}_h(\Pi_hq^\dagger)\big)_n  \\&\quad + 2\|\overline{u}^\dagger-\overline{u}_h(\Pi_hq^{\dagger})\|_n\|\overline{u}_h(q_h^*) - \overline{u}_h(\Pi_hq^\dagger)\|_n
 =:F_1 +F_2 +F_3 +F_4.
\end{align*}
Moreover, let $F_5=\gamma \|\Pi_hq^\dag\|_{H^1(\Omega)}^2$.
Then clearly, one of the five events $$\|\overline{u}_h(q_h^*)-\overline{u}_h(\Pi_hq^\dagger)\|_n^2 + \gamma \| q_h^*\|_{H^1\II}^2 \leq c F_i, \quad i=1,\ldots,5$$ must happen. We analyze the five cases separately. The proof is lengthy and technical, and it is divided into three steps. \\
\textbf{Step 1. Bound $F_1$ and $F_3$.} At this step, we crucially use an $H^2(\Omega)$ conforming FEM space $W_h$ (e.g., Hermite elements and Argyris elements for $d=1,2$ \cite{Brenner2002The} and Zhang elements for $d=3$ \cite{zhang2009family}). Specifically, we define an interpolation operator $\widetilde{\Pi}_h: H^2\II\mapsto W_h \subset H^2(\Omega)$ such that there holds for any $0\leq s\leq r\leq 2$:
$$\|v-\widetilde{\Pi}_h v\|_{H^s\II} \leq ch^{r-s}\|v\|_{H^r\II},\quad \forall v\in H^r(\Omega).$$
Moreover, the triangle inequality, Lemma \ref{lem:pointwise-L2-stab}, the $\|\cdot\|_n$ (cf. \eqref{eqn:approx-Pi-n}) and $L^2\II$ approximation properties of $\Pi_h$, and the approximation property and the $H^2\II$ stability of $\widetilde{\Pi}_h$ imply that for any $v\in H^2(\Omega)$, there holds
\begin{align}
    \|v-\widetilde{\Pi}_h v\|_n & \leq \|\big(v-\widetilde{\Pi}_h v) - \Pi_h\big(v-\widetilde{\Pi}_h v\big)\|_n + \|\Pi_h\big(v-\widetilde{\Pi}_h v\big)\|_n\nonumber\\
    & \leq ch^2\|v-\widetilde{\Pi}_h v\|_{H^2\II} + c\|\Pi_h\big(v-\widetilde{\Pi}_h v)\|_{L^2\II} \nonumber \\& \leq ch^2\|v-\widetilde{\Pi}_h v\|_{H^2\II} + c\|v-\widetilde{\Pi}_h v\|_{L^2\II}
    \leq ch^2\|v\|_{H^2\II}.\label{eqn:approx-Pih}
\end{align}
Obviously, we have the following decomposition
\begin{equation*}
    F_1 = 2\big(\bm\xi,\overline{u}_h(q_h^*) - \widetilde{\Pi}_h\overline{u}(q_h^*)\big)_n + 2\big(\bm\xi,\widetilde{\Pi}_h\overline{u}(q_h^*) - \overline{u}(q_h^*)\big)_n.
\end{equation*}
Now we define the discrete part $F_1^d(q_h)$ and regular part $F_1^r(q_h)$ of $F_1$ respectively by
\begin{align*}
    &F_1^d(q_h):= \frac{\overline{u}_h(q_h) - \widetilde{\Pi}_h\overline{u}(q_h)}{\|\overline{u}_h(q_h) - \overline{u}_h(\Pi_hq^\dag)\|_n+\gamma^{\frac12}\|q_h\|_{H^1\II}} \in X_h + W_h, \\&
    F_1^r(q_h):= \frac{\widetilde{\Pi}_h\overline{u}(q_h) - \overline{u}(q_h)}{\|\overline{u}_h(q_h) - \overline{u}_h(\Pi_hq^\dag)\|_n+\gamma^{\frac12}\|q_h\|_{H^1\II}}\in H^2\II.
\end{align*}
Then it suffices to bound $\sup_{q_h\in \mathcal{A}_h}|\big(\bm\xi,F_1^i\big)_n|$ with the index $i\in\{d,r\}$. Similar to the proof of Lemma \ref{lem:con-subGaussian}, one can show that $\{\big(\bm{\xi}, F_1^d(q_h) \big)_n, q_h\in\mathcal{A}_h \}$ is a sub-Gaussian random process equipped with the semimetric
\begin{equation*}
    \widehat{\rbd}(q_h,\widetilde{q}_h):= \sigma n^{-\frac12}\| F_1^d(q_h) -  F_1^d(\widetilde{q}_h)\|_n, \quad \forall q_h,\widetilde{q}_h \in \mathcal{A}_h.
\end{equation*}
Now by the triangle inequality, Lemma \ref{lem:point-err} and \eqref{eqn:approx-Pih}, we have
\begin{align*}
    \|\overline{u}_h(q_h)-\widetilde\Pi_h\overline{u}(q_h)\|_n & \leq \|\overline{u}_h(q_h)-\overline{u}(q_h)\|_n + \|\overline{u}(q_h)-\widetilde\Pi_h\overline{u}(q_h)\|_n\leq ch^2\|\overline{u}(q_h)\|_{H^2(\Omega)}.
\end{align*}
This and the stability estimate $\|\overline{u}(q_h)\|_{H^2\II} \leq c\|q_h\|_{H^1\II}$ from Lemma \ref{lem:reg-u} yield
\begin{align*}
    \diam(\mathcal{A}_h) & = \sup_{q_h,\widetilde{q}_h\in\mathcal{A}_h}\widehat{\rbd}(q_h,\widetilde{q}_h) \leq 2\sigma n^{-\frac12}\sup_{q_h\in\mathcal{A}_h}\| F_1^d(q_h)\|_n\\
    & \leq c\sigma n^{-\frac12}h^2\|\overline{u}(q_h)\|_{H^2\II}\Big(\gamma^{\frac12}\|q_h\|_{H^1\II}\Big)^{-1} \leq c\sigma n^{-\frac12}h^2\gamma^{-\frac 12}.
\end{align*}
By the triangle inequality, the $L^2(\Omega)$ approximation property of the operator $\widetilde{\Pi}_h$, the estimate \eqref{eqn:approx-L2}, and the stability estimate $\|\bar u(q_h)\|_{H^2(\Omega)}\leq c\|q_h\|_{H^1(\Omega)}$, we get
\begin{align}
    \sup_{q_h\in \mathcal{A}_h}\|F_1^d(q_h)\|_{L^2\II} & \leq \sup_{q_h\in \mathcal{A}_h}\frac{\|\overline{u}_h(q_h) - \overline{u}(q_h)\|_{L^2\II} + \|\overline{u}(q_h) - \widetilde{\Pi}_h\overline{u}(q_h)\|_{L^2\II}}{\|\overline{u}_h(q_h) - \overline{u}_h(\Pi_hq^\dag)\|_n+\gamma^{\frac12}\|q_h\|_{H^1\II}}\nonumber
    \\& \leq \sup_{q_h\in \mathcal{A}_h}\frac{\|\overline{u}_h(q_h) - \overline{u}(q_h)\|_{L^2\II}}{\|\overline{u}_h(q_h) - \overline{u}_h(\Pi_hq^\dag)\|_n+\gamma^{\frac12}\|q_h\|_{H^1\II}} + \sup_{q_h\in \mathcal{A}_h}\|F_1^r(q_h)\|_{L^2\II} \nonumber\\&
    \leq ch^2\|\overline{u}(q_h)\|_{H^2\II}\Big(\gamma^{\frac12}\|q_h\|_{H^1\II}\Big)^{-1}\leq ch^2\gamma^{-\frac 12}.\label{eqn:apriori-F1}
\end{align}
By the argument in Lemma \ref{lem:pointwise-L2-stab}, the following Lipschitz continuity property holds
\begin{equation*}
    \widehat{\rbd}(q_h,\widetilde{q}_h) \leq c\sigma n^{-\frac12}\|  F_1^d(q_h) -  F_1^d(\widetilde{q}_h)\|_{L^2\II},\quad \forall q_h,\widetilde{q}_h \in\mathcal{A}_h.
\end{equation*}
By Lemmas \ref{lem:max-ineq} and \ref{lem:cover-en-bound-finite}, the inequality $\log(1+x) \leq x$ over $\mathbb{R}_+$ and the estimate $\mathrm{dim}(X_h + W_h)=O(h^{-d})$, we get
\begin{align*}
 &\quad \bigg\|\,\sup_{q_h\in \mathcal{A}_h}|(\bm{\xi}, F_1^d(q_h) )_n|\,\bigg\|_{\psi_2} \leq c\int_0^{c\sigma n^{-\frac12}h^2\gamma^{-\frac12}}\sqrt{\log N\big(\tfrac{\epsilon}{2}, \mathcal{A}_h, \widehat{\rbd}\big)} \ {\rm d}\epsilon
    \\& \leq c\int_0^{c\sigma n^{-\frac12}h^2\gamma^{-\frac12}}\sqrt{\log N\big(\tfrac{\epsilon}{c\sigma n^{-\frac12}}, X_h + W_h, \|\cdot\|_{L^2\II}\big)} \ {\rm d}\epsilon \\
    &\leq c\int_0^{c\sigma n^{-\frac12}h^2\gamma^{-\frac12}} \dim(X_h + W_h)^{\frac12}\bigg(\log\Big( 1+ c\sigma n^{-\frac12}\sup_{q_h\in \mathcal{A}_h}\|F_1^d(q_h)\|_{L^2\II}\epsilon^{-1}\Big)\bigg)^{\frac12}\ {\rm d}\epsilon \\
    &\leq c\int_0^{c\sigma n^{-\frac12}h^2\gamma^{-\frac12}} \dim(X_h + W_h)^{\frac12}\Big(c\sigma n^{-\frac12}h^2\gamma^{-\frac12}\epsilon^{-1}\Big)^{\frac12}\ {\rm d}\epsilon \\& \leq c\sigma n^{-\frac12}h^{-\frac{d}{2}}h^2\gamma^{-\frac12} \leq c\gamma^{\frac{d}{12}}h^{2-\frac d2}\rho_0(1+\rho_0) ^{-\frac d6},
\end{align*}
using the choice $\gamma^{\frac12 + \frac{d}{12}} = O(\sigma n^{-\frac12} (1+\rho_0) ^{\frac d6}\rho_0^{-1})$.
Next, for the regular part $F_1^r$, Lemma \ref{Geer-entropy} for $r=\frac d2$, the approximation in $\|\cdot\|_n$ and the $H^2\II$ stability of $\widetilde{\Pi}_h$ and the estimate \eqref{eqn:approx-Pih} yield
\begin{align*}
 \sup_{q_h\in \mathcal{A}_h}\big(\bm\xi,F_1^r(q_h)\big)_n & = \mathcal{O}_p\Big(\sigma n^{-\frac 12}\| F_1^r(q_h)\|_n^{1-\frac d4}\|F_1^r(q_h)\|_{H^2\II}^{\frac d4}\Big) \\
 & \leq\mathcal{O}_p\Big(\sigma n^{-\frac 12}\Big( h^2\gamma^{-\frac12}\Big)^{1-\frac d4}(\gamma^{-\frac12})^{\frac d4}\Big) = \mathcal{O}_p\Big(\sigma n^{-\frac 12} h^{2-\frac d2}\gamma^{-\frac12}\Big).
\end{align*}
Similarly, from the choice $\gamma^{\frac12 + \frac{d}{12}} = O(\sigma n^{-\frac12}(1+\rho_0) ^{\frac d6}\rho_0^{-1})$, we deduce
\begin{equation*}
    \bigg\|\sup_{q_h\in \mathcal{A}_h}\big(\bm\xi,F_1^r(q_h)\big)_n\bigg\|_{\psi_2} \leq c\sigma n^{-\frac 12}h^{2-\frac d2} \gamma^{-\frac12} \leq c\gamma^{\frac{d}{12}}h^{2-\frac d2}\rho_0(1+\rho_0) ^{-\frac d6}.
\end{equation*}
Repeating the preceding argument yields the following bound on the term $F_3$:
\begin{align*}
\|\overline{u}_h(q_h^*)-\overline{u}_h(\Pi_hq^\dagger)\|_n^2 + \gamma \| q_h^*\|_{H^1\II}^2 &\leq  \mathcal{O}_p\Big(\sigma n^{-\frac12} h^{2-\frac d2}\|\Pi_hq^\dag\|_{H^1\II}\Big)\\ &\leq\mathcal{O}_p\Big(\gamma^{\frac12 + \frac{d}{12}} h^{2-\frac d2}\rho_0^2(1+\rho_0)^{-\frac{d}{6}}\Big).
\end{align*}
Thus it follows directly
\begin{equation*}
    A\leq c\gamma^{\frac{d}{12}}h^{4(\frac12-\frac{d}{8})}\rho_0(1+\rho_0) ^{-\frac d6} +  c\Big(\gamma^{\frac12+\frac{d}{12}}h^{4(\frac12-\frac{d}{8})}(1+\rho_0)^{-\frac{d}{6}}\Big)^{\frac 12}\rho_0.
\end{equation*}
\textbf{Step 2. Bound $F_2$.}
Upon letting $G=H^3(\Omega)$ and applying Lemma \ref{Geer-entropy} ($\gamma=\frac{d}{3}$, cf. Lemma \ref{lem:cover-en-bound-Sobolev}), we obtain
\begin{align*}
  &\big(\bm\xi,\overline{u}(q_h^*) - \overline{u}(\Pi_hq^\dag)\big)_n  = \mathcal{O}_{p}\left(\sigma n^{-\frac12} \|\overline{u}(q_h^*) - \overline{u}(\Pi_hq^\dag)\|_{n}^{1- \frac{d}{6}}{\big[(1+\rho_0)\|q_h^* - \Pi_hq^\dag\|_{H^1(\Omega)}\big]}^{\frac{d}{6}}\right).
\end{align*}
By the triangle inequality, Lemma \ref{lem:point-err} and the \textit{a priori} estimate $\|\overline{u}(q_h)\|_{H^2(\Omega)}\leq c\|q_h\|_{H^1(\Omega)}$ (cf. Lemma \ref{lem:reg-u}), we deuce
\begin{equation}\label{ineq:boundF2-pre-estimate}
    \begin{split}
        &\|\overline{u}(q_h^*) - \overline{u}(\Pi_hq^\dag)\|_{n}\\
 \leq & \|\overline{u}_h(q_h^*) - \overline{u}_h(\Pi_hq^\dag)\|_{n}+\|\overline{u}(q_h^*)-\overline{u}_h(q_h^*)\|_n + \|\overline{u}(\Pi_hq^\dag)-\overline{u}_h(\Pi_hq^\dag)\|_n\\
 \leq  & \|\overline{u}_h(q_h^*) - \overline{u}_h(\Pi_hq^\dag)\|_{n}+ch^2(\|q_h^*\|_{H^1(\Omega)} + \|\Pi_hq^\dag\|_{H^1(\Omega)}).
    \end{split}
\end{equation}
This and the triangle inequality yield
\begin{align*}
 &\|\overline{u}(q_h^*) - \overline{u}(\Pi_hq^\dag)\|_{n}^{1- \frac{d}{6}}\big[(1+\rho_0)\|q_h^* - \Pi_hq^\dag\|_{H^1(\Omega)}\big]^{\frac{d}{6}}\\
\leq & c\|\overline{u}_h(q_h^*) - \overline{u}_h(\Pi_hq^\dag)\|_{n}^{1-\frac{d}{6}}\big[(1+\rho_0)\|q_h^* - \Pi_hq^\dag\|_{H^1(\Omega)}\big]^{\frac{d}{6}}  \\
& + ch^{2(1-\frac{d}{6})}(1+\rho_0)^{\frac{d}{6}}(\|q_h^*\|_{H^1(\Omega)} + \|\Pi_hq^\dag\|_{H^1(\Omega)}),
\end{align*}
in view of the elementary inequality $(a+b)^s<c(a^s+b^s)$ (for $a,b>0$ and $s>0$). Next we consider the following two sub-cases separately: (i) $\|q_h^*\|_{H^1\II} \leq \|\Pi_hq^\dag\|_{H^1\II} $ and (ii) $\|q_h^*\|_{H^1\II} \geq \|\Pi_hq^\dag\|_{H^1\II} $. If the event (i) holds, then we have $\gamma^{\frac12}\|q_h^*\|_{H^1\II} \leq c\gamma^{\frac12}\|q^\dag\|_{H^1\II}$ and by the $H^1(\Omega)$ stability of $\Pi_h$,
\begin{align*}
    &\|\overline{u}_h(q_h^*)-\overline{u}_h(\Pi_hq^\dagger)\|_n^2 \\
     \leq &\mathcal{O}_{p}\left(\sigma n^{-\frac12} (1+\rho_0)^{\frac{d}{6}}\right)\Big(\|\overline{u}_h(q_h^*) - \overline{u}_h(\Pi_hq^\dag)\|_{n}^{1-\frac{d}{6}}\|\Pi_hq^\dag\|_{H^1(\Omega)}^{\frac{d}{6}} + h^{2(1-\frac{d}{6})} \|\Pi_hq^\dag\|_{H^1(\Omega)}\Big) \\
     \leq &\mathcal{O}_{p}\left(\sigma n^{-\frac12}(1+\rho_0)^{\frac{d}{6}} \right)\Big(\|\overline{u}_h(q_h^*) - \overline{u}_h(\Pi_hq^\dag)\|_{n}^{1-\frac{d}{6}}\rho_0^{\frac{d}{6}} + h^{2(1-\frac{d}{6})}\rho_0\Big).
\end{align*}
Then one of the follow two events must hold
\begin{equation}
\left\{\begin{aligned}
		\|\overline{u}_h(q_h^*)-\overline{u}_h(\Pi_hq^\dag)\|_n^{1+\frac d6} &\leq \mathcal{O}_{p}(\sigma n^{-\frac12})(1+\rho_0)^{\frac{d}{6}}\rho_0^{\frac{d}{6}}, \\
		\|\overline{u}_h(q_h^*)-\overline{u}_h(\Pi_hq^\dag)\|^2_n &\leq \mathcal{O}_{p}(\sigma n^{-\frac12})h^{2(1-\frac{d}{6})}(1+\rho_0)^{\frac{d}{6}}\rho_0.
\end{aligned}\right.
\end{equation}
Thus with the choice $\gamma^{\frac12 + \frac{d}{12}} = O(\sigma n^{-\frac12}(1+\rho_0) ^{\frac d6}\rho_0^{-1})$, we have
\begin{align*}
 A  \leq& c\gamma^{\frac12}\rho_0 + c\big(\sigma n^{-\frac12}(1+\rho_0)^{\frac{d}{6}}\rho_0^{\frac d6}\big)^{\frac{6}{6+d}} + c\Big(\sigma n^{-\frac12}h^{2(1-\frac{d}{6})}(1+\rho_0)^{\frac{d}{6}}\rho_0\Big)^{\frac 12} \\
    \leq &c\gamma^{\frac12}\rho_0 + c\Big(\gamma^{\frac12+\frac{d}{12}}h^{4(\frac12-\frac{d}{12})}\Big)^{\frac 12}\rho_0.
\end{align*}
Alternatively, if the event (ii) holds, then similarly, we obtain
\begin{align*}
    &\|\overline{u}_h(q_h^*)-\overline{u}_h(\Pi_hq^\dagger)\|_n^2 + \gamma\|q_h^*\|_{H^1\II}^2 \\
    \leq &\mathcal{O}_{p}\left(\sigma n^{-\frac12} (1+\rho_0)^{\frac{d}{6}}\right)\Big(\|\overline{u}_h(q_h^*) - \overline{u}_h(\Pi_hq^\dag)\|_{n}^{1-\frac{d}{6}}\|q_h^*\|_{H^1(\Omega)}^{\frac{d}{6}} + h^{2(1-\frac{d}{6})}\|q_h^*\|_{H^1(\Omega)}\Big).
\end{align*}
Then one of the following events must hold
\begin{equation}\label{eqn:bdd-stoch0}
\left\{\begin{aligned}
\|\overline{u}_h(q_h^*)-\overline{u}_h(\Pi_hq^\dagger)\|_n^2 + \gamma\|q_h^*\|_{H^1\II}^2 \leq &\mathcal{O}_{p}\left(\sigma n^{-\frac12} (1+\rho_0)^{\frac{d}{6}}\right)\|q_h^*\|_{H^1(\Omega)}^{\frac{d}{6}}\\
&\times\|\overline{u}_h(q_h^*) - \overline{u}_h(\Pi_hq^\dag)\|_{n}^{1-\frac{d}{6}},\\
\|\overline{u}_h(q_h^*)-\overline{u}_h(\Pi_hq^\dagger)\|_n^2 + \gamma\|q_h^*\|_{H^1\II}^2 \leq& \mathcal{O}_{p}\left(\sigma n^{-\frac12} (1+\rho_0)^{\frac{d}{6}}\right)h^{2(1-\frac{d}{6})}\|q_h^*\|_{H^1(\Omega)}.
\end{aligned}\right.
\end{equation}
Upon using the first inequality in \eqref{eqn:bdd-stoch0}, it follows directly that
\begin{align*}
    \gamma\|q_h^*\|_{H^1\II}^{2-\frac d6} &\leq \mathcal{O}_{p}\left(\sigma n^{-\frac12} (1+\rho_0)^{\frac{d}{6}} \right)\|\overline{u}_h(q_h^*) - \overline{u}_h(\Pi_hq^\dag)\|_{n}^{1-\frac{d}{6}}, \\
    \|\overline{u}_h(q_h^*)-\overline{u}_h(\Pi_hq^\dagger)\|_n^{1+\frac d6} &\leq \mathcal{O}_{p}\left(\sigma n^{-\frac12} (1+\rho_0)^\frac{d}{6}\right) \|q_h^*\|_{H^1(\Omega)}^{\frac{d}{6}}.
\end{align*}
Then direct computation yields the estimate
$$
\|q_h^*\|_{H^1\II} \leq  \mathcal{O}_{p}\left(\sigma n^{-\frac12} (1+\rho_0)^{\frac{d}{6}}\right)\gamma^{-(\frac12 + \frac{d}{12})}.
$$
Thus, by the choice $\gamma^{\frac12 + \frac{d}{12}} = O(\sigma n^{-\frac12}(1+\rho_0) ^{\frac d6}\rho_0^{-1})$, we get the following estimates
\begin{equation*}
\bigg\|\|q_h^*\|_{H^1\II}\bigg\|_{\psi_2} \leq  c\rho_0\quad\mbox{and}\quad
 \bigg\|\|\overline{u}_h(q_h^*)-\overline{u}_h(\Pi_hq^\dagger)\|_n\bigg\|_{\psi_2}  \leq c\gamma^{\frac 12} \rho_0.
\end{equation*}
Similarly, solving the second inequality in \eqref{eqn:bdd-stoch0} leads to
\begin{equation*}
    A \leq c\sigma n^{-\frac12}h^{2(1-\frac{d}{6})}(1+\rho_0)^{\frac{d}{6}}\gamma^{-\frac12} \leq c\gamma^{\frac{d}{12}}h^{4(\frac12-\frac{d}{12})}\rho_0.
\end{equation*}
Combining the preceding discussions gives
\begin{align*}
A\leq c\gamma^{\frac12}\rho_0 + c\gamma^{\frac {d} {12}}h^{4(\frac12-\frac{d}{12})}\rho_0 +  c\Big(\gamma^{\frac12+\frac{d}{12}}h^{4(\frac12-\frac{d}{12})}\Big)^{\frac 12}\rho_0.
\end{align*}
\textbf{Step 3. Bound  $F_4$ and $F_5$.} These two terms can be bounded directly as
\begin{align*}
	&\big\|\overline{u}_h(q_h^*) - \overline{u}_h(\Pi_hq^\dagger)\big\|_n + \gamma^{\frac12}\big\|q_h^*\big\|_{H^1\II}
   \leq  ch^2\rho_0 +c\gamma^{\frac12}\rho_0.
\end{align*}
Combining the preceding five cases and Lemma \ref{lem:bdd-sub-Gaussian} yield the first estimate in the $\psi_2$-Orlicz norm, and the second estimate follows from direct computation.
\end{proof}

\begin{remark}
 Let $$\zeta:=\gamma^{\frac12} +h^2 + \gamma^{\frac{d}{12}}h^{4(\frac12-\frac{d}{8})}  + (\gamma^{\frac12+\frac{d}{12}}h^{4(\frac12-\frac{d}{8})})^{\frac12},$$
which for the choice $h=O(\gamma^{\frac{6-d}{6(4-d)}})$ gives $\zeta=O(\gamma^\frac12)$.
By combining Lemma \ref{lem:approx-sol-uppbound-qh} with Lemmas \ref{lem:bdd-sub-Gaussian} and \ref{lem:Orlicz-bound}, we deduce that equivalently, the following estimates hold  in high probability for $z>0$,
\begin{align*}
\mathbb{P}(\|u_h(q_h^*) - u_h(\Pi_hq^\dagger)\|_n>\zeta\rho_0z) &\leq 2\exp(-cz^2),\\ \mathbb{P}(\| q^*_h\|_{H^1\II}>\gamma^{-\frac12}\rho_0\zeta z) &\leq 2\exp(-cz^2).
\end{align*}
\end{remark}

\begin{remark}
In Lemma \ref{lem:approx-sol-uppbound-qh}, the mesh size $h$ and the regularization  parameter $\gamma$ are independent. With respect to the mesh size $h$, the convergence rate appears sub-optimal. This is due to the fact that the P1 FEM space $X_h$ does not make full use of the $H^3(\Omega)$ regularity of the state $u(q)$, cf. Lemma \ref{lem:reg-u}. Indeed, using the P1 FEM can yield at most an $O(h^2)$ rate which requires only the regularity $u(q)\in H^2\II$. In contrast, with  P2 elements, we obtain $O(h^3)$ convergence, which results in the optimal mesh size $h=O(\gamma^{\frac{1}{6}})$, and accordingly, the bound can be improved to
$\gamma^{\frac12} +h^2 +  \gamma^{\frac{d}{12}}h^{6(\frac12-\frac{d}{12})} + (\gamma^{\frac12+\frac{d}{12}}h^{4(\frac12-\frac{d}{12})})^{\frac12} +\gamma^\frac{d}{12}h^{4(\frac12-\frac{d}{12})}$, and the choice {$h=O(\gamma^{\frac{1}{4}})$} suffices achieving the model error $\gamma^{\frac12}  \rho_0$. The only difference in the analysis of the P2 elements lies in the first term $F_1$, cf. Remark \ref{h3-converge}.
\end{remark}
\begin{remark}
The following sharper version of the inequality \eqref{ineq:boundF2-pre-estimate} holds
    \begin{equation*}
    \|\overline{u}(q_h^*) - \overline{u}(\Pi_hq^\dag)\|_{n} \leq \|\overline{u}_h(q_h^*) - \overline{u}_h(\Pi_hq^\dag)\|_{n}+ch^2.
\end{equation*}
We have kept $\|q_h^*\|_{H^1(\Omega)}$ and $\|\Pi_h q^\dag\|_{H^1(\Omega)}$ for the alternative argument in Lemma \ref{lem:approx-sol-uppbound-qh}.
\end{remark}

Now we can give the proof of Theorem \ref{thm:err-potential}, i.e., the error bound on the discrete approximation $q_h^*$.
\begin{proof}[Proof of Theorem \ref{thm:err-potential}]
By Lemmas \ref{lem:point-err} and \ref{lem:approx-sol-uppbound-qh} and \eqref{ineq:proba-distrib-estimate}, the following two estimates hold with probability at
least $1-2\tau$
\begin{align*}
\|u_h(q_h^*)-u^\dag\|_n  &\leq \|u_h(q_h^*)-u_h(\Pi_hq^\dag)\|_n + \|u_h(\Pi_hq^\dag)-u^\dag\|_n \leq c\big(\eta\rho_0\ell_\tau + h^2\big),\\
\|q_h^*\|_{H^1\II} &\leq c\gamma^{-\frac12}\eta\rho_0\ell_\tau .
\end{align*}
For any $\varphi\in H_0^1\II$, the weak formulations of $u^\dag\equiv u(q^\dag)$ and $u_h(q_h^*)$ yield
\begin{equation*}
\begin{aligned}
	\big((q^\dagger-q_h^*)u^\dagger,\varphi \big)
	&=\big((q^\dagger-q_h^*)u^\dagger,\varphi-P_h\varphi\big) + \big( (q^\dagger-q_h^*)u^\dagger,P_h\varphi \big ) \\
	&=\big( (q^\dag-q_h^*)u^\dag,\varphi-P_h\varphi \big)+ \big(\nabla(u_h(q_h^*)-u^\dag),\nabla P_h\varphi\big) \\& \quad +\big(q_h^*(u_h(q_h^*)-u^\dag),P_h\varphi\big) =: {\rm I} + {\rm II} + {\rm III}.
\end{aligned}
\end{equation*}
Let $\varphi=(q^\dag-q_h^*)u^\dag$. Then direct computation gives
$\nabla\varphi = u^\dag\nabla( q^\dag -q_h^*)+(q^\dag-q^*_h)\nabla u^\dag.$
Using Assumption \ref{ass:reg} and the box constraint on $q^\dag$ and $q_h^*$, the uniform bound $\|u(q^\dag)\|_{H^2(\Omega)}\leq c$ in Lemma \ref{lem:reg-u} and the Sobolev embedding $H^2(\Omega)\hookrightarrow C(\overline{\Omega})$ from Lemma \ref{lem:Sob-embed}, we obtain $\|\varphi\|_{L^2\II} \leq c$ and
\begin{align*}	\|\nabla\varphi\|_{L^2\II} & \leq \|\nabla( q^\dag - q_h^*)\|_{L^2\II} \|u^\dag\|_{L^\infty\II} + \|q^\dag - q_h^*\|_{L^\infty\II} \|\nabla u^\dag\|_{L^2\II} \\ & \leq c\big( 1+\|\nabla q^\dag\|_{L^2\II}+\|\nabla q_h^*\|_{L^2\II} \big) \leq c\big(1+\gamma^{-\frac12}\eta\rho_0\ell_\tau  \big).
\end{align*}
This implies $\varphi\in H_0^1\II$. It suffices to bound the three terms ${\rm I}$, ${\rm II}$ and ${\rm III}$ separately.
By Lemma \ref{lem:approx-sol-uppbound-qh}, the approximation property \eqref{ineq:Ph-approx} of $P_h$ and the stability estimate \eqref{ineq:Pih-approx-2} of $\Pi_h$, the following bound holds
\begin{equation*}
	|{\rm I}| \leq  ch\|\nabla\varphi\|_{L^2\II} \leq ch(1+\|\nabla q_h^*\big\|\|_{L^2\II}) \leq ch\big(1+\gamma^{-\frac12}\eta\rho_0\ell_\tau  \big).
\end{equation*}
Next we bound the second term ${\rm II}$. Since $q_h^*\in\mathcal{A}_h$, there exists $c>0$ independent of $q_h^*$ such that $\|u(q_h^*)\|_{H^2\II} \leq c$.
Hence it follows directly that
\begin{equation*}
	\|\nabla \big(u(q_h^*)-u_h(q_h^*)\big)\|_{L^2\II} \leq ch.
\end{equation*}
Then Gagliardo--Nirenberg interpolation inequality \cite{BrezisMironescu:2018} (cf. the inequality \eqref{eqn:Gargliardo}) and Lemmas \ref{lem:connect-seminorm-stdnorm}, \ref{lem:reg-u} and \ref{lem:point-err} imply
\begin{align*}
\|\nabla \big(u^\dag-u(q_h^*)\big)\|_{L^2\II}  &\leq \|u^\dag-u(q_h^*)\|^\frac23_{L^2\II}\|u^\dagger-u(q_h^*)\|^{\frac13}_{H^3\II}\\& \leq c\big(\|u^\dag-u(q_h^*)\|_n + n^{-\frac3d}\|u^\dagger-u(q_h^*)\|_{H^3\II}\big)^\frac23\|u^\dagger-u(q_h^*)\|^{\frac13}_{H^3\II}.
\end{align*}
Meanwhile, by Lemma \ref{lem:reg-u}, we have $\|u^\dagger-u(q_h^*)\|_{H^3\II}\leq c(1+\rho_0)\|q^\dag-q_h^*\|_{H^1(\Omega)}$. These two estimates, Lemma \ref{lem:connect-seminorm-stdnorm} and the triangle inequality yield
\begin{equation*}
    \|\nabla \big(u^\dag-u(q_h^*)\big)\|_{L^2\II} \leq c\big(\eta + h^2 + n^{-\frac3d}\big)^\frac23 \big(1+\gamma^{-\frac12}\eta\rho_0\ell_\tau\big)^2.
\end{equation*}
Thus the Cauchy-Schwarz inequality,  the triangle inequality  and the $H^1(\Omega)$ stability  of $P_h$ (cf. \eqref{ineq:Ph-approx}) give
\begin{equation*}
\begin{aligned}
	|{\rm II}| &\leq \|\nabla(u^\dag-u_h(q_h^*))\|_{L^2\II}\|\nabla P_h\varphi\|_{L^2\II} \\
	&\leq  c\big(\|\nabla(u^\dag-u(q_h^*))\|_{L^2\II} + \|\nabla(u(q_h^*)-u_h(q_h^*))\|_{L^2\II}\big) \|\nabla \varphi\|_{L^2\II} \\
	&\leq  c\big[\big(\eta + h^2 + n^{-\frac3d}\big)^\frac23 +h\big]\big(1+\gamma^{-\frac12}\eta\rho_0\ell_\tau\big)^3.
\end{aligned}
\end{equation*}
Together with the \textit{a priori} estimate $\|\nabla(u^\dagger-u_h(q_h^*))\|_{L^2\II}\leq c$, this leads to
\begin{equation*}
	|{\rm II}| \leq  c\big[\min(1, \big(\eta + h^2 + n^{-\frac3d} \big)^\frac23 + h ) \big]\big(1+\gamma^{-\frac12}\eta\rho_0\ell_\tau  \big)^3.
\end{equation*}
Last, we bound the third term ${\rm III}$. By the standard duality argument and the \textit{a priori} regularity estimate  $\|u(q_h^*)\|_{H^2\II} \leq c$ from Lemma \ref{lem:reg-u}, we get
\begin{equation*}
	\|u(q_h^*)-u_h(q_h^*)\|_{L^2\II} \leq ch^2.
\end{equation*}
Thus the Cauchy-Schwarz inequality, the $L^2\II$ stability of $P_h$, Lemma \ref{lem:connect-seminorm-stdnorm}, and the triangle inequality imply
\begin{align*}
|{\rm III}| & \leq  c \| P_h\varphi  \|_{L^2\II} \|u^\dag-u_h(q_h^*)\|_{L^2\II} \\
&
\leq  c \|\varphi  \|_{L^2\II} \big(\|u^\dag - u (q_h^*)\|_{L^2\II} + \|u (q_h^*)-u_h(q_h^*)\|_{L^2\II}\big)\\
 & \leq c\big(\eta + h^2 + n^{-\frac3d}\big)\big(1+\gamma^{-\frac12}\eta\rho_0\ell_\tau  \big)^2.
\end{align*}
Combing the bounds on ${\rm I}$, ${\rm II}$ and ${\rm III}$ gives the bound on $\|u^\dag(q_h^*-q^\dag)\|_{L^2(\Omega)}$. Under Condition \ref{Cond:positivity}, repeating the argument in the proof of Theorem \ref{thm:con-err-potential} yields the desired $L^2\II$ estimate.
\end{proof}

\section{Numerical experiments and discussion}
\label{sec:numerics}
In this section we present numerical results to complement the theoretical results in Section \ref{sec:stoc-conver}.
\subsection{Adaptive parameter choice}
The suitable choice of the penalty parameter $\gamma$ in the formulation \eqref{eqn:dis-optim-pure-dis}--\eqref{eqn:dis-weak-potential} is one of the most important issues in practical inversion.
The theoretical result in Section \ref{sec:stoc-conver} motivates the following \textit{a priori} choice:
\begin{equation}\label{eqn:prior-selection-gamma}
	\gamma^{\frac12 + \frac{d}{12}} = O( \sigma n^{-\frac12}\big(\|q^\dagger\|_{H^1\II} + \sigma n^{-\frac12}\big)^{-1} ).
\end{equation}
However this rule requires an \textit{a priori} knowledge about $\|q^\dag\|_{H^1(\Omega)}$ and $\sigma$, which is generally unknown in practice, and thus it often cannot be applied directly. Following the strategy in \cite{chen2018stochastic}, one may employ a self-consistent iteration by estimating $\|q^\dag\|_{H^1(\Omega)}$ using $\|q_{\gamma_k,h}^*\|_{H^1(\Omega)}$ and $\sigma$ by $\|u_h(q_{\gamma_k,h}^*)-\boldsymbol{m}\|_n$. This leads to an adaptive algorithm for automatically choosing the regularization parameter $\gamma$; see Algorithm \ref{alg:adpa-gamma} for the detailed procedure (for $d=2$). In practice, the constant hidden in $O(\cdot)$ may be taken to be 1, as in the algorithm. This rule is heuristic in nature, since it does not use the knowledge of the noise level.

\begin{algorithm}[hbt!]
	\caption{Adaptive choice of the regularization parameter $\gamma$. \label{alg:adpa-gamma}}
	\begin{algorithmic}[1]
		\STATE Give an initial guess of $\gamma_0$ and set the maximum iteration number $K$.
		\FOR{$k=0,\ldots,K$}
		\STATE Find a minimizer $q_{\gamma_k,h}^*\in \mathcal{A}_h$ to problem \eqref{eqn:dis-optim-pure-dis}-\eqref{eqn:dis-weak-potential}
		with $\gamma = \gamma_k$.
		\STATE Update $\gamma_{k+1}$ by
	\begin{equation}\label{eqn:gamma} \gamma^{\frac23}_{k+1} =  O(n^{-\frac12} \|u_h(q_{\gamma_k,h}^*)-\boldsymbol{m}\|_n (\| q_{\gamma_{k},h}^*\|_{H^1(\Omega)} + n^{-1/2}{\|u_h(q_{\gamma_k,h}^*)-\boldsymbol{m}\|_n}  )^{-1} ).
 \end{equation}
		\ENDFOR
	\end{algorithmic}
\end{algorithm}

The next lemma shows the monotone convergence of the sequence generated by Algorithm \ref{alg:adpa-gamma}. Numerically we observe that the algorithm reaches convergence within ten iterations, and hence it is highly efficient.
\begin{lemma}\label{lem:gamma}
The sequence $(\gamma_k)_{k=0}^\infty$ given by Algorithm \ref{alg:adpa-gamma} is monotonically convergent.
\end{lemma}
\begin{proof}
We fix the constant in $O(\cdot)$ to be 1.
The minimizing property of $q_{\gamma,h}^*$ to $J_{\gamma,h}$ implies that
$\|u_h(q_{\gamma,h}^*)-\boldsymbol{m}\|_n$ monotonically increasing with $\gamma$, and $\|q_{\gamma,h}^*\|_{H^1(\Omega)}$ is monotonically decreasing with $\gamma$ \cite[Lemma 3.3, p. 42]{ItoJin:2015}. Now the condition $\gamma_{k-1}\leq \gamma_k$ implies
\begin{align*}
    \|u_h(q_{\gamma_{k-1},h}^*)-\boldsymbol{m}\|_n \leq \|u_h(q_{\gamma_k,h}^*)-\boldsymbol{m}\|_n\quad\mbox{and}\quad
    \| q_{\gamma_{k-1},h}^*\|_{H^1(\Omega)} \geq \| q_{\gamma_{k},h}^*\|_{H^1(\Omega)}.
\end{align*}
This, the update formula \eqref{eqn:gamma}, {and the monotonicity of the function $t\mapsto\frac{t}{c+t}$} in $t\geq0$ (for any $c\geq0$) imply
\begin{align*}
    \gamma_{k+1}^\frac{2}{3} - \gamma_k^\frac{2}{3} = &\frac{n^{-1/2} \|u_h(q_{\gamma_k,h}^*)-\boldsymbol{m}\|_n }{\| q_{\gamma_{k},h}^*\|_{H^1(\Omega)} + n^{-1/2} {\|u_h(q_{\gamma_k,h}^*)-\boldsymbol{m}\|_n}}\\
      &- \frac{n^{-1/2} \|u_h(q_{\gamma_{k-1},h}^*)-m\|_n}{\| q_{\gamma_{k-1},h}^*\|_{H^1(\Omega)} + n^{-1/2} {\|u_h(q_{\gamma_{k-1},h}^*)-\boldsymbol{m}\|_n}}\geq 0,
\end{align*}
i.e., $\gamma_{k+1}\geq \gamma_k$. The other case $\gamma_{k-1}\geq \gamma_k$ follows similarly. Thus the desired monotonicity relation holds.
The update formula \eqref{eqn:gamma} implies that the sequence $(\gamma_k)_{k=0}^\infty$ is uniformly bounded by $0 $ from below and $1$ from above, and hence it is monotonically convergent.
\end{proof}

\begin{remark}
Both Algorithm \ref{alg:adpa-gamma} and Lemma \ref{lem:gamma} are for the noisy data $\bm m$ with a fixed realization of the noise. Thus, they differ from the stochastic convergence results in Section \ref{sec:stoc-conver}. Note that the results in Section \ref{sec:stoc-conver} rely crucially on the a priori choice of $\gamma$, e.g., $\gamma^{{\frac12 + \frac{d}{12}}} = O(\sigma n^{-\frac12}(1+\rho_0)^{\frac d6}\rho_0^{-1})$ in Theorem \ref{thm:con-err-potential}. The stochastic convergence analysis under the heuristic choice of the regularization parameter $\gamma$ by Algorithm \ref{alg:adpa-gamma} remains open, for which a new analysis strategy and additional assumptions on the noise are likely needed.
\end{remark}

\subsection{Numerical results and discussions}

Now we present numerical experiments to complement the theoretical analysis. Unless otherwise stated, we take the  domain $\Omega$ to be the unit square $\Omega=(0,1)^2\subset\mathbb{R}^2$, and take the observation points $(x_i)_{i=1}^n$ to be uniformly distributed over the domain $\Omega$; see Fig. \ref{fig:quasiuniform-mesh-points}(a) for a schematic illustration. We generate the finite element mesh $\mathcal{T}_h$ in two steps: first divide $\Omega$ into $M\times M$ small squares (each has a side length $1/M$) and then connect the lower left and upper right vertices of each square. We obtain the exact data $u^\dag$ by solving problem \eqref{eqn:potential} on a finer grid with a mesh size $h = \frac{1}{200}$, and generate the noise data $\boldsymbol{m}$ by
\begin{equation*}
	m_i = u^\dagger(x_i) + \sigma\|u^\dagger\|_{L^\infty(\Omega)}\xi_i,\quad i=1,2,\dots,n,
\end{equation*}
where the random variables $\bm \xi = (\xi_i)_{i=1}^n$ follow the standard normal distribution and $\sigma$ denotes the noise strength. The resulting optimization problems are minimized using the conjugate gradient method; see Algorithm \ref{alg:cgm} in Appendix \ref{app:CG} for the implementation detail. All the experiments are run on a personal laptop. The lower and upper  bounds of the admissible set $\mathcal{A}$ are set to $c_0=1.0$ and $c_1 = 5.0$, respectively. We choose the regularization parameter $\gamma$ using the adaptive rule in Algorithm \ref{alg:adpa-gamma}. To measure the convergence behavior, we employ two error metrics: $e_q:=\|q^\dagger-q^*_h\|_{L^2\II}/\|q^\dagger\|_{L^2(\Omega)}$ and $e_u:= \|u^\dagger - u_h(q_h^*)\|_n$.

The first example is to illustrate the \textit{a priori} choice \eqref{eqn:prior-selection-gamma}.
\begin{example}\label{exam:opt-gamma}
Let $q^\dagger=3.0+\frac12\sin(\pi x)\sin(\pi y)$ and $f\equiv1$. Fix the number $n$ of  sampling points at $401^2$, and a mesh size $h=\frac{1}{36}$ small enough to neglect the FEM discretization error. Consider two noise levels: $\sigma=5.0\%$ and $\sigma=1.0\%$. Then we solve the regularized problem \eqref{eqn:conti-optim-prob-potential}-\eqref{eqn:conti-weak-potential} with $\gamma$ varying in the set $\{10^{-6},10^{-7},10^{-8},10^{-9},10^{-10},10^{-11},10^{-12}\}$.
\end{example}

\begin{figure}[hbt!]
	\centering
	\setlength{\tabcolsep}{0pt}
	\begin{tabular}{cc}
		\includegraphics[width=0.48\textwidth]{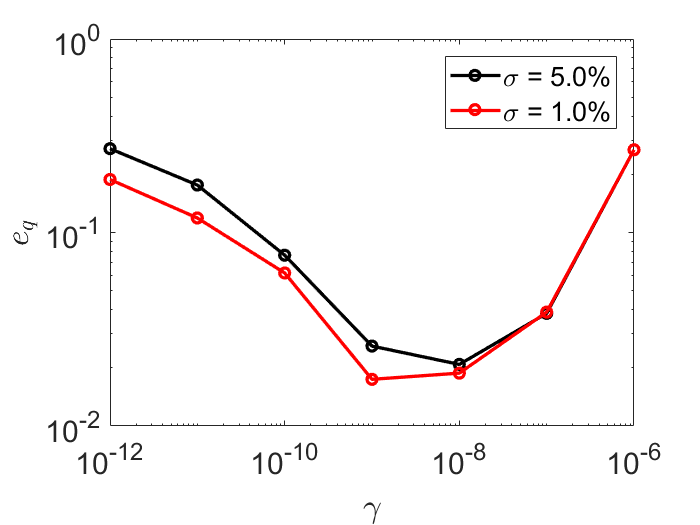} & \
		\includegraphics[width=0.48\textwidth]{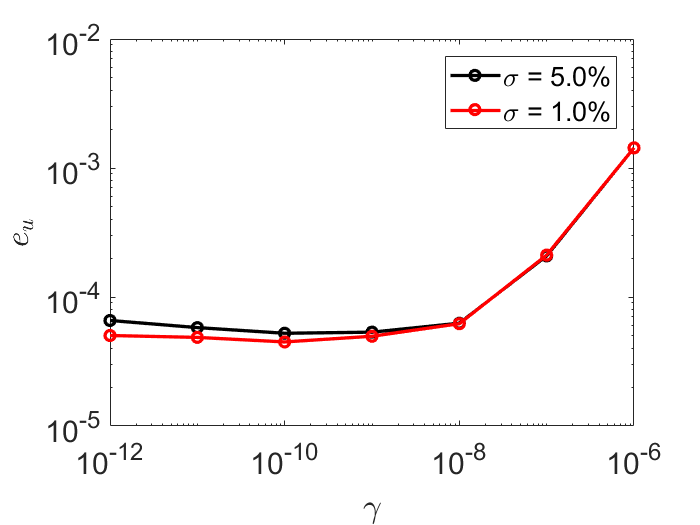}
	\end{tabular}
	\caption{The numerical results on the optimal regularization parameter $\gamma^*$ at two noise levels: the error $e_q$ versus the parameter $\gamma$ {\rm(}left{\rm)}; the error $e_u$ versus the parameter $\gamma$ {\rm(}right{\rm)}.}
	\label{fig:opt-gamma}
\end{figure}

Fig. \ref{fig:opt-gamma} shows the error $e_q$ of the minimizer $q_h^*$ and the state approximation $e_u:=\|u^\dagger - u_h(q_h^*)\|_n$. The optimal parameters are  found to be $10^{-8}$ ($\sigma=5.0\%$) and $10^{-9}$ ($\sigma=1.0\%$), which are close to the automatically determined $3.43$e-9 ($\sigma=5.0\%$) and $3.07$e-10 ($\sigma=1.0\%$). Thus, the rule \eqref{eqn:prior-selection-gamma} can provide a nearly optimal \textit{a priori} choice of $\gamma$.

The second example focuses on the adaptive algorithm to find a suitable $\gamma$.
\begin{example}\label{exam:adap}
Let $q^\dagger=3.0+\frac12\sin(\pi x)\sin(\pi y)$ and $f\equiv1$. Fix the number $n$ of sampling points at $n=501^2$, and a mesh size $h=\frac{1}{36}$. Then we solve the regularized problem \eqref{eqn:conti-optim-prob-potential}-\eqref{eqn:conti-weak-potential} with a noise level $\sigma = 5.0\%$ and the initialization $\gamma_0 = n^{-\frac34}$.
\end{example}

\begin{figure}[hbt!]
	\label{fig:adap}
	\centering
	\setlength{\tabcolsep}{0pt}
	\begin{tabular}{ccc}
		\includegraphics[width=0.33\textwidth]{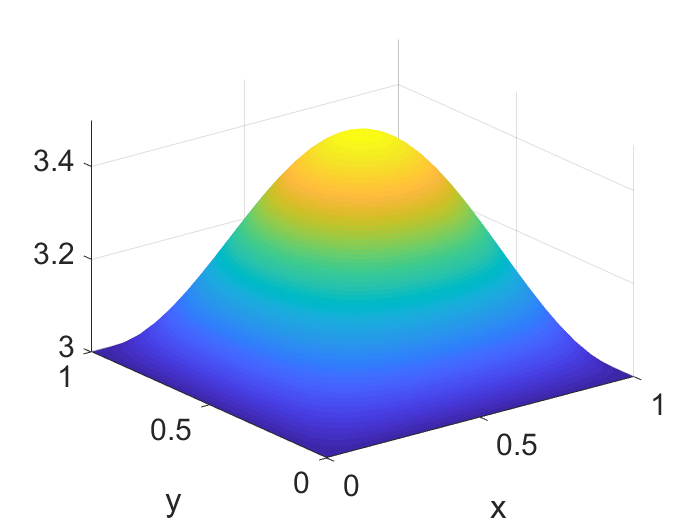} & \includegraphics[width=0.33\textwidth]{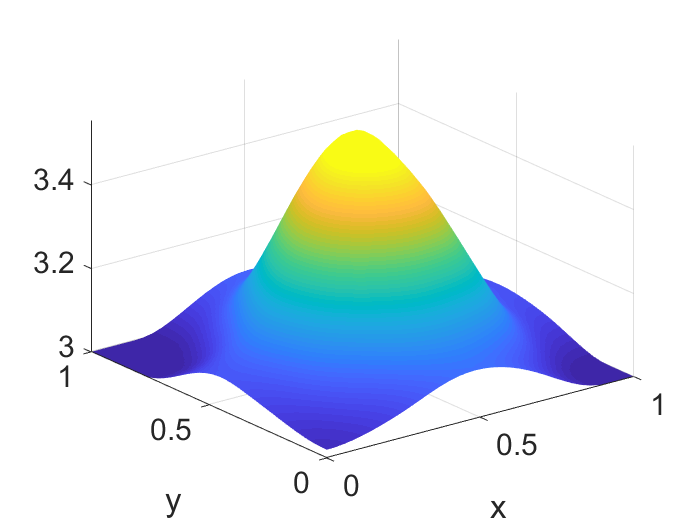} &
		\includegraphics[width=0.33\textwidth]{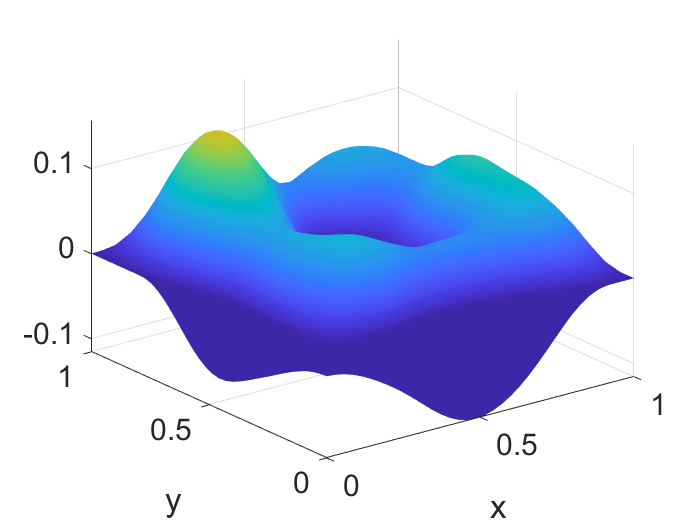} \\
        \includegraphics[width=0.33\textwidth]{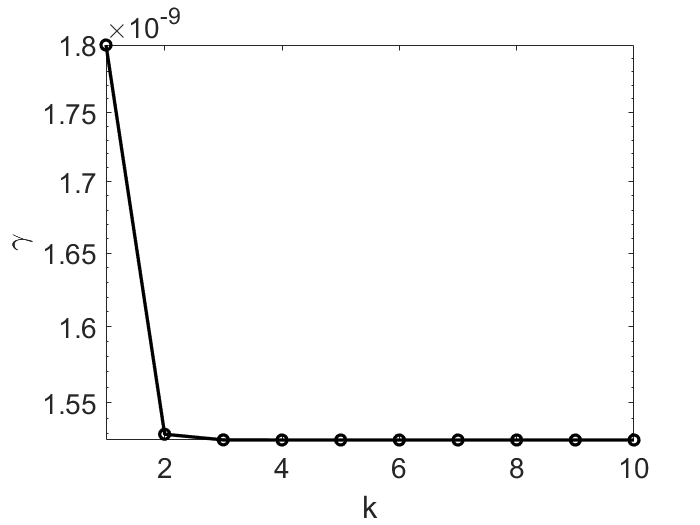} &
        \includegraphics[width=0.33\textwidth]{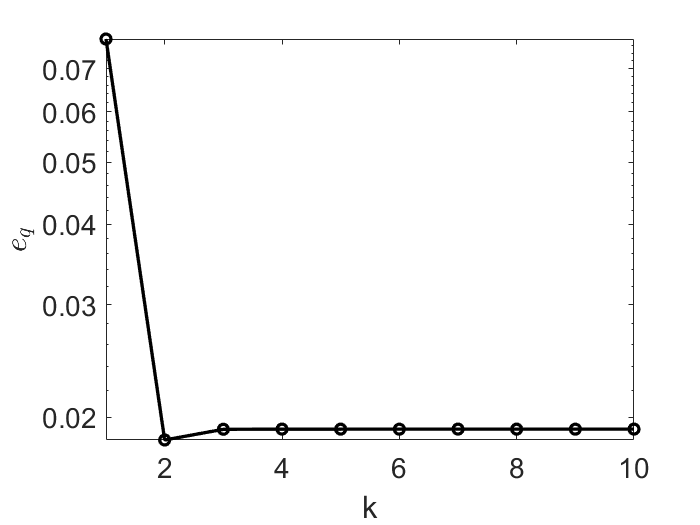} &
        \includegraphics[width=0.33\textwidth]{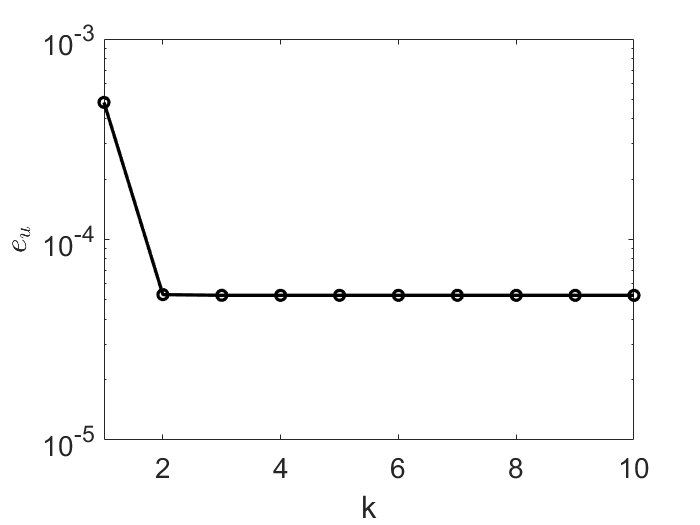}
	\end{tabular}
	\caption{The numerical results for Example \ref{exam:adap} with $\sigma=5.0\%$. The top row shows the exact potential $q^\dag$ (left), recovered potential $q_{h}^*$ from Algorithm \ref{alg:adpa-gamma} (middle), and the pointwise error $e:=q^\dag - q_{h}^*$ (right). The bottom row shows the value of $\gamma$ (left), the error $e_q$ (middle) and the error $e_u$ (right), all versus the iteration index $k$.}
\end{figure}

Example \ref{exam:adap} illustrates the effectiveness of Algorithm \ref{alg:adpa-gamma} for adaptively choosing the regularization parameter $\gamma$. In practice, the true potential $q^\dagger$ and the noise level $\sigma$ of the measurement data are typically unknown \textit{a priori}. To approximate these quantities, we utilize the numerical solution $q_h^*$ by setting $\|q^\dag\|_{H^1(\Omega)}\approx \|q_h^*\|_{H^1(\Omega)}$ and $\sigma\approx\|u_h(q_h^*) - m\|_n$. Fig. \ref{fig:adap} presents the results  by Algorithm \ref{alg:adpa-gamma} after several iterations. The reconstruction capture well the overall shape of the exact potential $q^\dag$, and the error near the boundary is dominant. The regularization parameter obtained by the adaptive algorithm converges to $\gamma = 1.53 \times 10^{-9}$,  which agree well with the theoretical prediction $2.46 \times 10^{-9}$.
This method converges in around ten steps, and hence it is computationally efficient. The monotone convergence of the algorithm is clearly observed also in the bottom row of Fig. \ref{fig:adap}.

The next example is to verify the stochastic convergence result in Theorem \ref{thm:err-potential}.

\begin{example}\label{exam:rate}
 Let $\Omega = (0,1)$ or $(0,1)^2$. Let $\gamma$ be given by \eqref{eqn:prior-selection-gamma} and $h = O(\gamma^\frac14)$. Consider the following four cases:
\begin{enumerate}
\item[{\rm(a)}] $\sigma = 5.0\%$, $q^\dagger = 1.0 + \frac12\sin^2(\pi x)$ and $f\equiv1$, with $n \in \{6.4, 14.4, 19.6, 36.0\}\times 10^5$;
\item[{\rm(b)}] $\sigma = 1.0\%$, $q^\dagger = 1.0 + \frac12\sin^2(\pi x)$ and $f\equiv1$, with $n \in\{6.4, 14.4, 19.6, 36.0\}\times 10^5$;
\item[{\rm(c)}] $\sigma = 2.0\%$, $q^\dagger = 3.0 + \frac12\sin(\pi x)\sin(\pi y)$ and $f\equiv1$, with $n \in \{101^2,201^2,301^2,401^2\}$;
\item[{\rm(d)}] $\sigma = 1.0\%$, $q^\dagger = 3.0 + \frac12\sin(\pi x)\sin(\pi y)$ and $f\equiv1$, with $n \in \{101^2,201^2,301^2,401^2\}$.
\end{enumerate}
\end{example}

\begin{figure}[hbt!]
	\label{fig:rate}
	\centering
	\setlength{\tabcolsep}{0pt}
    \begin{tabular}{ccc}
		\includegraphics[width=0.32\textwidth]{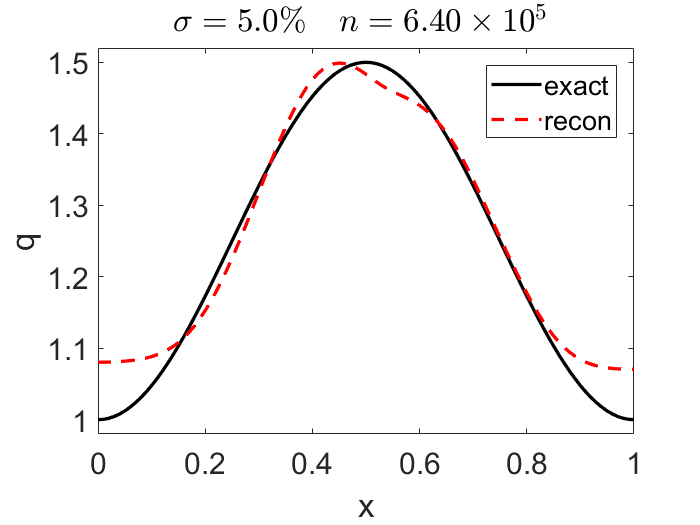} & \includegraphics[width=0.32\textwidth]{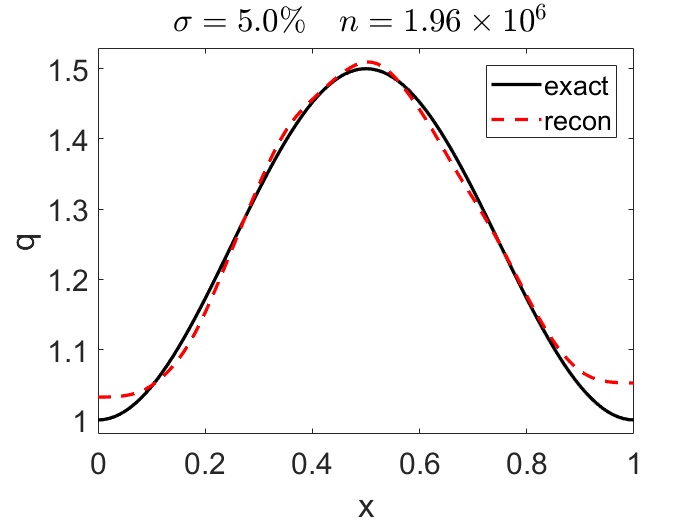} &
		\includegraphics[width=0.32\textwidth]{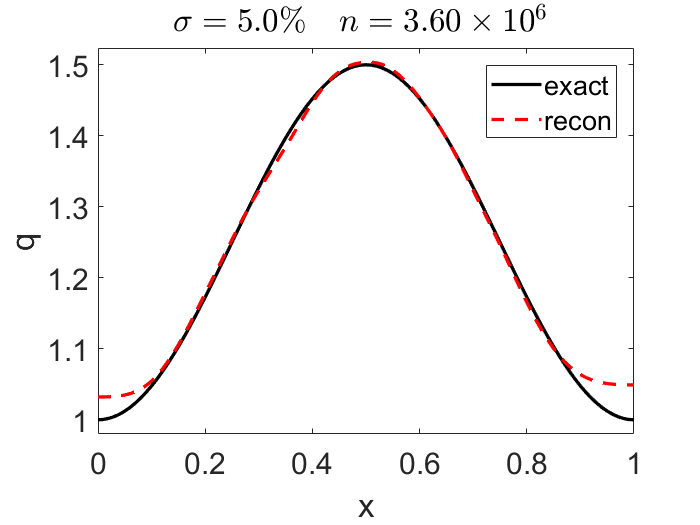} \\
		$\gamma = 5.72$e-10 & $\gamma = 2.19$e-10 & $\gamma = 1.19$e-10
	\end{tabular}
    \begin{tabular}{ccc}
		\includegraphics[width=0.32\textwidth]{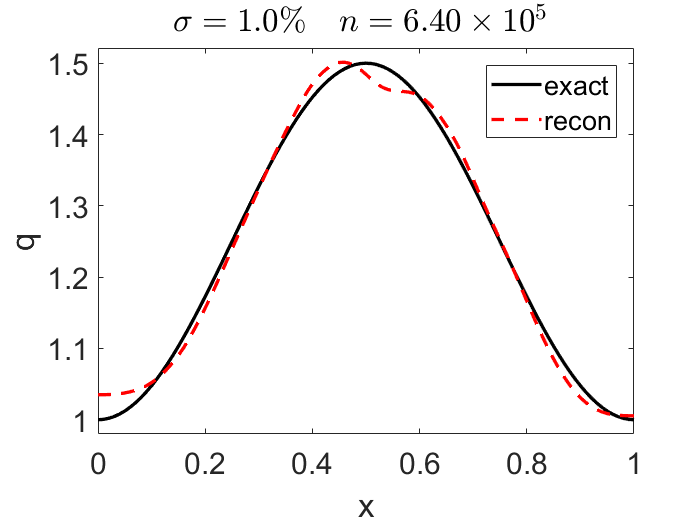} & \includegraphics[width=0.32\textwidth]{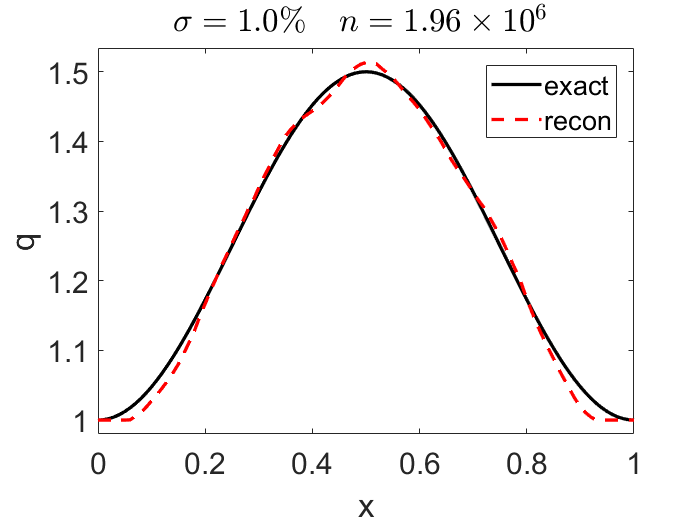} &
		\includegraphics[width=0.32\textwidth]{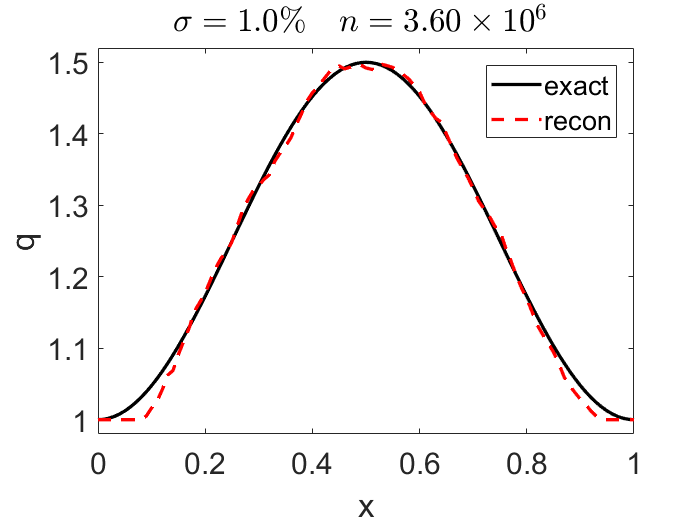} \\
		$\gamma = 3.62$e-11 & $\gamma = 1.39$e-11 & $\gamma = 8.24$e-12
	\end{tabular}
	\begin{tabular}{ccc}
		\includegraphics[width=0.32\textwidth]{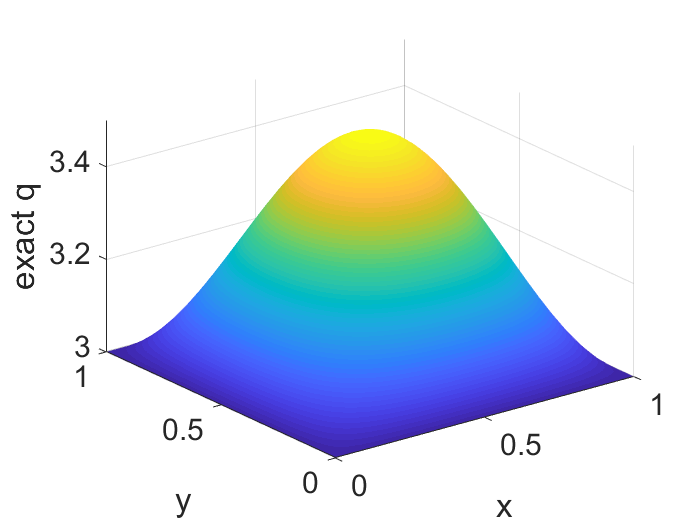} & \includegraphics[width=0.32\textwidth]{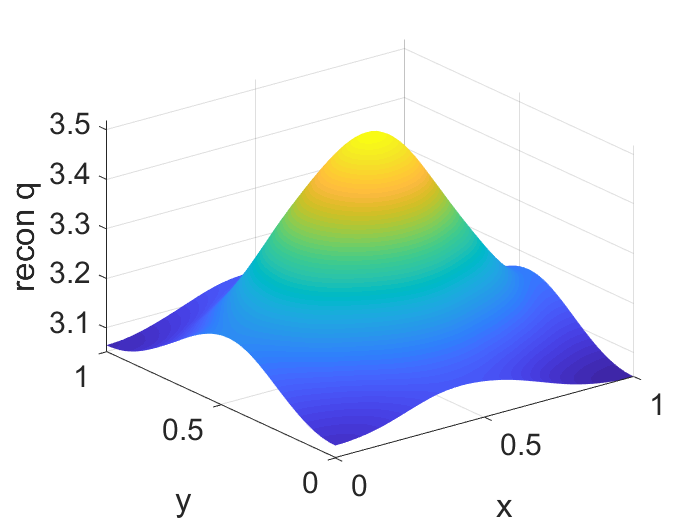} &
		\includegraphics[width=0.32\textwidth]{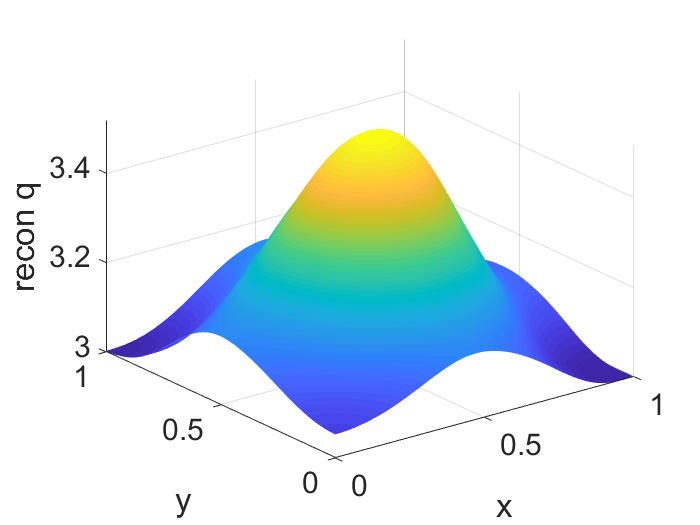} \\
		exact & $\gamma = 1.34$e-9 & $\gamma = 8.70$e-10
	\end{tabular}
    \begin{tabular}{ccc}
		\includegraphics[width=0.32\textwidth]{InvPotenfigs/exactq.png} & \includegraphics[width=0.32\textwidth]{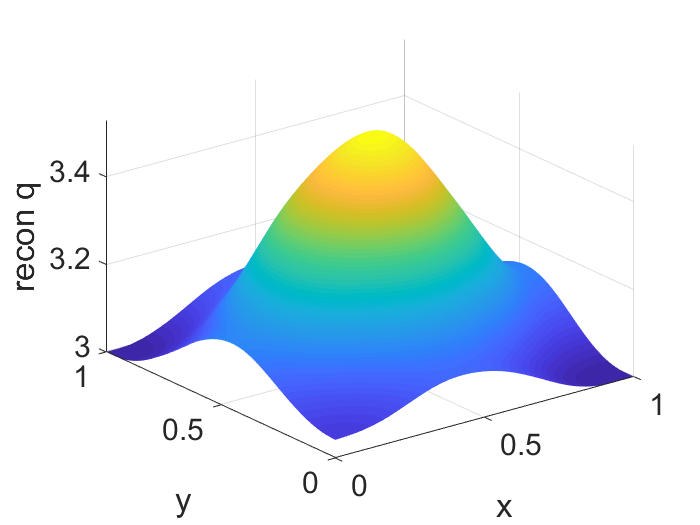} &
		\includegraphics[width=0.32\textwidth]{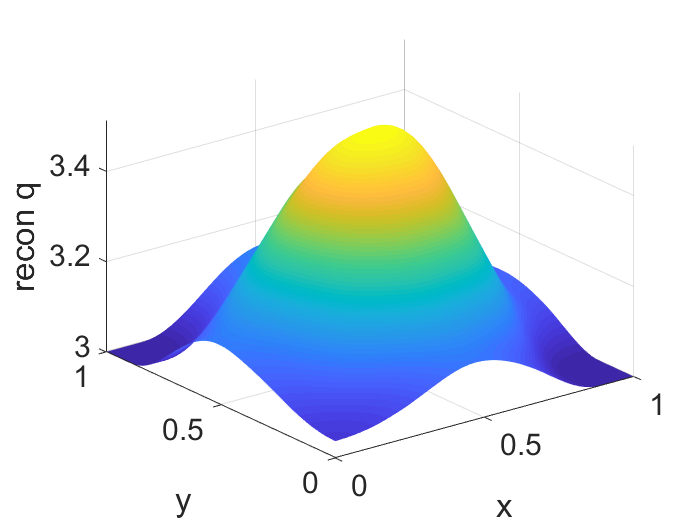} \\
		exact & $\gamma = 4.73$e-10 & $\gamma = 3.07$e-10
	\end{tabular}
	\caption{The reconstruction $q_{h}^*$ for Example \ref{exam:rate} (a)-(d): (a) and (b) at the first two rows; (c) and (d) at the last two rows.}
\end{figure}

In Table \ref{tab:rate}, we present the numerical results for the errors $e_q$ and $e_u$. The results show that the two error metrics $e_q$ and $e_u$ both decay to zero when $n$ tends to $\infty$ (see the final column), provided that the algorithmic parameters $h$ and $\gamma$ are chosen according to Theorem \ref{thm:err-potential}. Fig. \ref{fig:rate} shows the exemplary reconstructions for the one- and two-dimensional cases. Similar to Example \ref{exam:adap}, the pointwise  errors are concentrated near the boundaries.

\begin{table}[hbt!]
	\caption{Numerical results for Example \ref{exam:rate}}
	\centering
	\begin{threeparttable}
		\subfigure[case (a).]{
		\begin{tabular}{c|ccccc}
		\toprule
		 $\gamma$ & 5.72e-10 & 2.85e-10 & 2.19e-10 & 1.30e-10 & trend \\
		\midrule
	     $e_q$   & 2.55e-2 & 1.43e-2 & 1.42e-2 & 1.18e-2 & $\searrow$ \\
		 $e_u$   & 9.94e-6 & 6.79e-6 & 7.07e-6 & 6.07e-6 & $\searrow$ \\
		\bottomrule
		\end{tabular}}
        \subfigure[case (b).]{
		\begin{tabular}{c|ccccc}
		\toprule
		 $\gamma$ & 3.62e-11 & 1.80e-11 & 1.39e-11 & 8.24e-12 & trend \\
		\midrule
	     $e_q$    & 1.16e-2 & 1.05e-2 & 9.96e-3 & 9.11e-3 & $\searrow$ \\
		 $e_u$    & 4.96e-6 & 4.35e-6 & 4.21e-6 & 4.06e-6 & $\searrow$ \\
		\bottomrule
		\end{tabular}}
	\subfigure[case (c).]{
		\begin{tabular}{c|ccccc}
			\toprule
			$\gamma$ & 6.88e-9 & 2.45e-9 & 1.34e-9 & 8.70e-10 & trend \\
			\midrule
			 $e_q$   & 2.61e-2 & 1.85e-2 & 1.61e-2 & 1.41e-2 & $\searrow$ \\
			 $e_u$   & 3.99e-5 & 2.88e-5 & 2.81e-5 & 2.51e-5 & $\searrow$ \\
			\bottomrule
	\end{tabular}}
	\subfigure[case (d).]{
		\begin{tabular}{c|ccccc}
			\toprule
			$\gamma$ & 2.43e-9 & 8.66e-10 & 4.73e-10 & 3.07e-10 & trend \\
			\midrule
			 $e_q$   & 2.00e-2 & 1.47e-2 & 1.34e-2 & 1.31e-2 & $\searrow$ \\
			 $e_u$   & 3.21e-5 & 2.52e-5 & 2.50e-5 & 1.93e-5 & $\searrow$ \\
			\bottomrule
	\end{tabular}}
	\end{threeparttable}
    \label{tab:rate}
\end{table}

\section*{Acknowledgments} The authors are grateful to two anonymous referees and the associate editor for their constructive comments which have led a significant improvement in the quality of the paper.

\appendix
\section{Additional proofs}\label{app:tech-proof}

\subsection{Proof of Proposition \ref{prop:exist-cont}}
In the proof, fix any realization of the random measurement data  $\bm{m}$.
First, by Lemma \ref{lem:reg-u}, for $f\in L^2(\Omega)$ and for any $q\in \mathcal{A}$, the state $u(q)\in H^2(\Omega)$. This and the Sobolev embedding $H^2(\Omega)\hookrightarrow C(\overline{\Omega})$ (for $d=1,2,3$) in Lemma \ref{lem:Sob-embed} imply that the pointwise evaluation $u(q)(x_i)$ makes sense, and also the functional $J_\gamma(q)$ is indeed well defined over $q\in\mathcal{A}$.

The existence of a minimizer $q^*\in \mathcal{A}$ to problem \eqref{eqn:conti-optim-prob-potential}-\eqref{eqn:conti-weak-potential} follows by a standard argument in calculus of variation \cite{Dacorogna:2008}. Since the functional $J_\gamma(q)$ is nonnegative, there exists a minimizing sequence $(q^k)_{k=1}^\infty\subset\mathcal{A} $ such that
\begin{equation*}
    \lim_{k\to\infty} J_\gamma(q^k) = \inf_{q\in \mathcal{A}}J_\gamma(q).
\end{equation*}
Then the sequence $(q^k)_{k=1}^\infty\subset \mathcal{A}$ is uniformly bounded in $H^1(\Omega)$, and by Lemma \ref{lem:reg-u}, the sequence $(u(q^k) )_{k=1}^\infty \subset H_0^1(\Omega)\cap H^2(\Omega)$ is uniformly bounded in $H^2(\Omega)$. Hence, there exist two subsequences of $(q^k)_{k=1}^\infty$ and $(u(q^k))_{k=1}^\infty$, still denoted by $(q^k)_{k=1}^\infty$ and $( u(q^k))_{k=1}^\infty$, and some $q^* \in \mathcal{A}$ and $u^*\in H_0^1(\Omega)\cap H^2(\Omega)$ such that $q^k\to q^*$ weakly in $H^1(\Omega)$ and $u(q^k)\to u^*$ weakly in $H^2(\Omega)$. Then we claim that the following identity holds
\begin{equation}\label{eqn:weak-seq-close}
u^*=u(q^*).
\end{equation}
This claim and the compact embedding $H^2(\Omega)\hookrightarrow C(\overline{\Omega})$ (for $d=1,2,3$), cf. Lemma \ref{lem:Sob-embed} lead to
the convergence $u(q^k)\to u(q^*)$ in the discrete semi-norm $\|\cdot\|_n$.
This result and the weak lower semicontinuity of the $H^1(\Omega)$ norm imply
\begin{align*}
    J_\gamma(q^*) & = \|u(q^*)-\bm m\|_n^2 + \gamma \|q^*\|_{H^1(\Omega)}^2\\
     & \leq \lim_{k\to\infty}\|u(q^k)-\bm m\|_n^2 + \gamma \liminf_{k\to\infty}\|q^k\|_{H^1(\Omega)}^2\\
     & \leq \liminf_{k\to\infty}\big(\|u(q^k)-\bm m\|_n^2 + \gamma \|q^k\|_{H^1(\Omega)}^2\big)\\
     &= \liminf_{k\to\infty} J_\gamma(q^k) = \inf_{q\in \mathcal{A}}J_\gamma(q).
\end{align*}
That is, $q^*$ is a global minimizer to problem \eqref{eqn:conti-optim-prob-potential}-\eqref{eqn:conti-weak-potential}.
Next we show the claim \eqref{eqn:weak-seq-close}. By the  definition of $u(q^k)$, we have
\begin{equation}\label{eqn:weak-form-uqn}
    \int_\Omega \nabla u(q^k)\cdot\nabla \varphi +q^ku(q^k)\varphi \,\d x = \int_\Omega f\varphi \,\d x,\quad \forall \varphi \in H_0^1(\Omega).
\end{equation}
Fix any $\varphi\in H_0^1(\Omega)$. The weak convergence of $u(q^k)$ to $u^*$  in $H^2(\Omega)$ implies
\begin{equation}\label{eqn:limit1-app}
    \lim_{k\to \infty} \int_\Omega \nabla u(q^k)\cdot \nabla \varphi \,\d x = \int_\Omega \nabla u^*\cdot \nabla \varphi \,\d x.
\end{equation}
Next, by the weak convergence of $u(q^k)$ to $u^*$ in $H^2(\Omega)$, the compact embedding
$H^1(\Omega)\hookrightarrow \hookrightarrow L^2(\Omega)$ and the uniform boundedness of $(q^k)_{k=1}^\infty\subset\mathcal{A}$ in the $L^\infty(\Omega)$ norm, we have
\begin{align*}
   &\lim_{k\to\infty} \left|\int_\Omega q^k(u(q^k)-u^*)\varphi\,\d x\right| \leq \lim_{k\to\infty}\|q^k\|_{L^\infty(\Omega)}\|u(q^k)-u^*\|_{L^2(\Omega)}\|\varphi\|_{L^2(\Omega)} =0.
\end{align*}
Likewise, the weak convergence of $q^k$ in $H^1(\Omega)$, and $u^*\varphi\in L^2(\Omega)$ (due to the continuous embedding $H^1(\Omega)\hookrightarrow L^4(\Omega)$ for $d=1,2,3$, cf. Lemma \ref{lem:reg-u}) imply
\begin{align*}
 \lim_{k\to \infty}\int_\Omega q^ku^*\varphi \,\d x = \int_\Omega q^*u^*\varphi \,\d x.
\end{align*}
The last two identities immediately imply \begin{equation}\label{eqn:limit2-app}
\lim_{k\to\infty}\int_\Omega q^ku(q^k)\varphi \,\d x = \int_\Omega q^*u^*\varphi\,\d x.
\end{equation}
Thus, by passing to the limit $k\to \infty$ in the identity \eqref{eqn:weak-form-uqn} and using the identities \eqref{eqn:limit1-app} and \eqref{eqn:limit2-app}, we deduce that the limit $u^*\in H_0^1(\Omega)\cap H^2(\Omega)$ satisfies
\begin{equation*}
    \int_\Omega \nabla u^*\cdot\nabla \varphi\, \d x + \int_\Omega q^*u^*\varphi \,\d x = \int_\Omega f\varphi \,\d x,\quad \forall \varphi\in H^1_0(\Omega).
\end{equation*}
This and the definition of $u(q^*)$ yield $u^*=u(q^*)$, thereby showing the claim \eqref{eqn:weak-seq-close}. Since for each realization of the random data $\bm{m}$, we have the existence of a minimizer, the statement of Proposition \ref{prop:exist-cont} holds $\mathbb{P}$ almost surely. Moreover, by repeating the preceding argument, one can show that the minimizer $q^*$ is continuous with respect to the perturbation in the data $\bm m$ (see, e.g., \cite[Theorem 4.2, p. 109]{ItoJin:2015}).

\subsection{Proof of Proposition \ref{prop:exist-disc}}

Note that for any $q_h\in\mathcal{A}_h$, there exists a unique $u_h(q_h)\in X_h$ solving problem \eqref{eqn:dis-weak-potential}. Meanwhile, the finite element function $u_h(q_h)$ is continuous over $\overline\Omega$, so the pointwise evaluation $u_h(q_h)(x_i)$ is well defined. Thus the discrete semi-norm $\|u_h(q_h)-\bm m\|_n$ for each realization of the random data $\bm{m}$ and the functional $J_{\gamma,h}$ are both well defined over $\mathcal{A}_h$.

The existence of a minimizer $q_h^*\in \mathcal{A}_h$ follows similarly as Proposition \ref{prop:exist-cont} and relies on the finite-dimensionality of the discrete problem \eqref{eqn:dis-optim-pure-dis}-\eqref{eqn:dis-weak-potential}. Since the functional $J_{\gamma,h}(q_h)$ is bounded from below by zero, there exists a minimizing sequence $(q_h^k)_{k=1}^\infty\subset \mathcal{A}_h$:
\begin{equation*}
    \lim_{k\to\infty} J_{\gamma,h}(q_h^k) = \inf_{q_h\in\mathcal{A}_h}J_{\gamma,h}(q_h).
\end{equation*}
Meanwhile, for each $q_h\in \mathcal{A}_h$, by setting $\varphi_h=u_h(q_h)$ in \eqref{eqn:dis-weak-potential} and Poincar\'{e} inequality, we deduce that the state $u_h(q_h)$ can be bounded in $H^1(\Omega)$ independent of $q_h$. Thus, there exist two subsequences of $ (q_h^k)_{k=1}^\infty$ and $( u_h(q_h^k))_{k=1}^\infty$, still denoted by $( q_h^k)_{k=1}^\infty$ and $( u_h(q_h^k))_{k=1}^\infty$, and some $q_h^*\in \mathcal{A}_h\subset H^1(\Omega)$ and $u_h^* \in X_h\subset H_0^1(\Omega)$ such that $q_h^k\to q_h^*$ weakly in $H^1(\Omega)$ and $u_h(q_h^k)\to u_h^*$ weakly in $H^1(\Omega)$. Now due to the finite-dimensionality of the spaces $V_h$ and $X_h$ and the norm equivalence in a finite-dimensional space, we have  $q_h^k\to q^*$ in $H^1(\Omega)$ and $u_h(q_h^n) \to u_h^*$ in $C(\overline{\Omega})$. Furthermore, by repeating the argument of the claim \eqref{eqn:weak-seq-close}, we can deduce  $u_h^*=u_h(q^*)$. Then we can easily prove that the limit $q_h^*\in\mathcal{A}_h $ is indeed a minimizer of problem \eqref{eqn:dis-optim-pure-dis}--\eqref{eqn:dis-weak-potential}, and the statement of Proposition \ref{prop:exist-disc} holds $\mathbb{P}$ almost surely.

\section{Conjugate gradient algorithm}\label{app:CG}
Now we briefly describe the conjugate gradient algorithm
for minimizing problem \eqref{eqn:conti-optim-prob-potential}--\eqref{eqn:conti-weak-potential}. At each step, the algorithm has to compute the gradient $J'_\gamma(q)$ of $J_\gamma(q)$. This can
be achieved using the adjoint state: find $v(q)\in H_0^1$ such that
\begin{equation}\label{eqn:adj}
(\nabla v(q),\nabla\varphi)+(qv(q),\varphi) = 2(u(q)-\bm \xi,\varphi)_n, \quad\forall \varphi\in H^1_0\II.
\end{equation}
To evaluate $v_h = \sum_{i=1}^{N}v_i\varphi_i(x)$ at $(x_j)_{j=1}^n$, we form explicitly a sparse matrix $\mathbf{P}=[\varphi_i(x_j) ]\in\mathbb{R}^{N\times n}$.
Using the adjoint state $v(q)$, the $L^2(\Omega)$ gradient $J_\gamma'(q)$ is given by
\begin{equation}\label{eqn:grad}
	J_\gamma'(q) = -\gamma(\Delta+I) q + \nabla u(q)\cdot \nabla v(q) ,
\end{equation}
and the descent direction $g^k=-(-\Delta + I)^{-1}J_\gamma'(q^k)$ (with a zero Neumann boundary condition). The conjugate gradient direction $d_k$ is given by
\begin{equation}\label{eqn:cg-dir}
	d_k = g^k +\beta_k d_{k-1} , \quad \mbox{with } \beta_k = \|g^{k}\|_{L^2(\Omega)}^2/\|g^{k-1}\|_{L^2(\Omega)}^2,
\end{equation}
with $\beta_0=0$. The complete procedure of the method is listed in Algorithm \ref{alg:cgm}.

\begin{algorithm}[hbt!]
	\caption{Conjugate gradient method. \label{alg:cgm}}
	\begin{algorithmic}[1]
		\STATE Set the maximum iteration number $K$, and choose $q^0$.
		\FOR{$k=1,\ldots,K$}
		\STATE Solve for $u(q^k)$ the solution to problem \eqref{eqn:conti-weak-potential} with $q=q^k$;
		\STATE Solve for $v(q^k)$ the solution to the modified adjoint problem \eqref{eqn:adj} with $q=q^k$;
		\STATE Compute the gradient $J_\gamma'(q^k)$ via \eqref{eqn:grad}, and the descent direction $d_k$ via \eqref{eqn:cg-dir};
		\STATE Compute the step length $s_k$;
		\STATE Update the diffusion coefficient by $q^{k+1}=P_\mathcal{A}(q^k+s_kd_k)$;
		\STATE Check the stopping criterion.
		\ENDFOR
	\end{algorithmic}
\end{algorithm}

\bibliographystyle{siam}
\bibliography{reference}

\begin{thebibliography}{10}

\bibitem{Adams2003Sobolev}
{\sc R.~A. Adams and J.~J.~F. Fournier}, {\em {Sobolev Spaces}},
  Elsevier/Academic Press, Amsterdam, second~ed., 2003.

\bibitem{birman1967piecewise}
{\sc M.~S. Birman and M.~Z. Solomyak}, {\em Piecewise-polynomial approximations
  of functions of the classes {$W_p^\alpha$}}, Mat. Sb. (N.S.), 115 (1967),
  pp.~331--355.

\bibitem{Brenner2002The}
{\sc S.~C. Brenner and L.~R. Scott}, {\em {The Mathematical Theory of Finite
  Element Methods}}, Springer-Verlag, New York, second~ed., 2002.

\bibitem{Brezis2010FunctionalAS}
{\sc H.~Brezis}, {\em {Functional Analysis, {S}obolev Spaces and Partial
  Differential Equations}}, Springer, New York, 2011.

\bibitem{BrezisMironescu:2018}
{\sc H.~Brezis and P.~Mironescu}, {\em Gagliardo-{N}irenberg inequalities and
  non-inequalities: the full story}, Ann. Inst. H. Poincar\'e{} C Anal. Non
  Lin\'eaire, 35 (2018), pp.~1355--1376.

\bibitem{Chen2020Convergence}
{\sc D.-H. Chen, D.~Jiang, and J.~Zou}, {\em Convergence rates of {T}ikhonov
  regularizations for elliptic and parabolic inverse radiativity problems},
  Inverse Problems, 36 (2020), pp.~075001, 21 pp.

\bibitem{chen2018stochastic}
{\sc Z.~Chen, R.~Tuo, and W.~Zhang}, {\em Stochastic convergence of a
  nonconforming finite element method for the thin plate spline smoother for
  observational data}, SIAM J. Numer. Anal., 56 (2018), pp.~635--659.

\bibitem{chen2022stochastic}
{\sc Z.~Chen, W.~Zhang, and J.~Zou}, {\em Stochastic convergence of regularized
  solutions and their finite element approximations to inverse source
  problems}, SIAM J. Numer. Anal., 60 (2022), pp.~751--780.

\bibitem{choulli2021some}
{\sc M.~Choulli}, {\em Some stability inequalities for hybrid inverse
  problems}, C. R. Math. Acad. Sci. Paris, 359 (2021), pp.~1251--1265.

\bibitem{Dacorogna:2008}
{\sc B.~Dacorogna}, {\em {Direct Methods in the Calculus of Variations}},
  Springer, New York, second~ed., 2008.

\bibitem{Engl1996Regularization}
{\sc H.~W. Engl, M.~Hanke, and A.~Neubauer}, {\em {Regularization of Inverse
  Problems}}, Kluwer Academic Publishers Group, Dordrecht, 1996.

\bibitem{EnglKunischNeubauer:1989}
{\sc H.~W. Engl, K.~Kunisch, and A.~Neubauer}, {\em Convergence rates for
  {T}ikhonov regularisation of nonlinear ill-posed problems}, Inverse Problems,
  5 (1989), pp.~523--540.

\bibitem{GineNickl:2016}
{\sc E.~Gin\'{e} and R.~Nickl}, {\em {Mathematical Foundations of
  Infinite-Dimensional Statistical Models}}, Cambridge University Press,
  Cambridge, 2016.

\bibitem{Grisvard:2011}
{\sc P.~Grisvard}, {\em {Elliptic Problems in Nonsmooth Domains}}, SIAM,
  Philadelphia, 2011.

\bibitem{hao2010convergence}
{\sc D.~N. H{\`a}o and T.~N.~T. Quyen}, {\em Convergence rates for {T}ikhonov
  regularization of coefficient identification problems in {L}aplace-type
  equations}, Inverse Problems, 26 (2010), p.~125014.

\bibitem{ItoJin:2015}
{\sc K.~Ito and B.~Jin}, {\em Inverse {P}roblems: Tikhonov {T}heory and
  {A}lgorithms}, World Scientific Publishing Co. Pte. Ltd., Hackensack, NJ,
  2015.

\bibitem{jin2022convergence}
{\sc B.~Jin, X.~Lu, Q.~Quan, and Z.~Zhou}, {\em Convergence rate analysis of
  {G}alerkin approximation of inverse potential problem}, Inverse Problems, 39
  (2022), pp.~015008, 26 pp.

\bibitem{Kekkonen:2022}
{\sc H.~Kekkonen}, {\em Consistency of {B}ayesian inference with {G}aussian
  process priors for a parabolic inverse problem}, Inverse Problems, 38 (2022),
  pp.~035002, 29 pp.

\bibitem{Koers:2024}
{\sc G.~Koers, B.~Szab\'{o}, and A.~van~der Vaart}, {\em Linear methods for
  non-linear inverse problems}.
\newblock Preprint, arXiv:2411.19797v1, 2024.

\bibitem{larsson2003partial}
{\sc S.~Larsson and V.~Thom{\'e}e}, {\em {Partial Differential Equations with
  Numerical Methods}}, Springer, Berlin, 2003.

\bibitem{Nickl:2020}
{\sc R.~Nickl}, {\em Bernstein--von {M}ises theorems for statistical inverse
  problems {I}: {S}chr\"odinger equation}, J. Eur. Math. Soc. (JEMS), 22
  (2020), pp.~2697--2750.

\bibitem{Nickl:2023}
\leavevmode\vrule height 2pt depth -1.6pt width 23pt, {\em {Bayesian Non-Linear
  Statistical Inverse Problems}}, EMS Press, Berlin, 2023.

\bibitem{NicklWang:2020}
{\sc R.~Nickl, S.~van~de Geer, and S.~Wang}, {\em Convergence rates for
  penalized least squares estimators in {PDE} constrained regression problems},
  SIAM/ASA J. Uncertain. Quantif., 8 (2020), pp.~374--413.

\bibitem{NickWang:2024}
{\sc R.~Nickl and S.~Wang}, {\em On polynomial-time computation of
  high-dimensional posterior measures by {L}angevin-type algorithms}, J. Eur.
  Math. Soc. (JEMS), 26 (2024), pp.~1031--1112.

\bibitem{pennes1948analysis}
{\sc H.~H. Pennes}, {\em Analysis of tissue and arterial blood temperatures in
  the resting human forearm}, J. Appl. Physiol., 1 (1948), pp.~93--122.

\bibitem{Siebel:2024}
{\sc M.~Siebel}, {\em Convergence rates for the maximum a posteriori estimator
  in {PDE}-regression models with random design}.
\newblock Preprint, arXiv:2409.03417v2, 2024.

\bibitem{Stone:1982}
{\sc C.~J. Stone}, {\em Optimal global rates of convergence for nonparametric
  regression}, Ann. Statist., 10 (1982), pp.~1040--1053.

\bibitem{thomee2007galerkin}
{\sc V.~Thom{\'e}e}, {\em {Galerkin Finite Element Methods for Parabolic
  Problems}}, Springer, New York, 2006.

\bibitem{trucu2010space}
{\sc D.~Trucu, D.~B. Ingham, and D.~Lesnic}, {\em Space-dependent perfusion
  coefficient identification in the transient bio-heat equation}, J. Eng.
  Math., 67 (2010), pp.~307--315.

\bibitem{utreras1988convergence}
{\sc F.~I. Utreras}, {\em Convergence rates for multivariate smoothing spline
  functions}, J. Approx. Theory, 52 (1988), pp.~1--27.

\bibitem{geer2000empirical}
{\sc S.~A. van~de Geer}, {\em Empirical Processes in M-estimation}, Cambridge
  university press, Cambridge, 2000.

\bibitem{vanDerVaart:1996}
{\sc A.~W. van~der Vaart and J.~A. Wellner}, {\em {Weak Convergence and
  Empirical Processes}}, Springer-Verlag, New York, 1996.
\newblock With applications to statistics.

\bibitem{yamamoto2001simultaneous}
{\sc M.~Yamamoto and J.~Zou}, {\em Simultaneous reconstruction of the initial
  temperature and heat radiative coefficient}, Inverse Problems, 17 (2001),
  pp.~1181--1202.

\bibitem{yang2008inverse}
{\sc L.~Yang, J.-N. Yu, and Z.-C. Deng}, {\em An inverse problem of identifying
  the coefficient of parabolic equation}, Appl. Math. Model., 32 (2008),
  pp.~1984--1995.

\bibitem{zhang2009family}
{\sc S.~Zhang}, {\em A family of {3D} continuously differentiable finite
  elements on tetrahedral grids}, Appl. Numer. Math., 59 (2009), pp.~219--233.

\bibitem{zhang2022identification}
{\sc Z.~Zhang, Z.~Zhang, and Z.~Zhou}, {\em Identification of potential in
  diffusion equations from terminal observation: analysis and discrete
  approximation}, SIAM J. Numer. Anal., 60 (2022), pp.~2834--2865.

\end{thebibliography}
\end{document}